\documentclass[11pt]{article}

\usepackage{amscd,amsmath, amssymb, fancyhdr, epsfig, color, tocloft, mathtools}
\usepackage{dutchcal, extarrows, enumerate, xypic}

\usepackage[backref=page]{hyperref}
\renewcommand*{\backrefalt}[4]{%
	\ifcase #1 (Not cited.)%
	\or        (Cited on page~#2.)%
	\else      (Cited on pages~#2.)%
	\fi}

\hypersetup{
	colorlinks   = true,
	citecolor    = magenta,
	linkcolor    = blue,
	urlcolor     = magenta	
}

\newcommand{\version}{version 1.0,\ \ June 4, 2026}
\makeatletter
\def\x@arrow{\DOTSB\Relbar}
\def\xlongequalsignfill@{\arrowfill@\x@arrow\Relbar\x@arrow}
\providecommand{\xlongequal}[2][]{%
	\ext@arrow 0099\xlongequalsignfill@{#1}{#2}}
\def\xlongrightarrowfill@{\arrowfill@\relbar\relbar\longrightarrow}
\makeatother

\numberwithin{equation}{section}

\def\eqref#1{(\ref{#1})}
\newcommand{\goth}{\mathfrak}
\newcommand{\g}{{\mathfrak g}}

\newcommand{\Z}{{\mathbb Z}}
\newcommand{\C}{{\mathbb C}}
\newcommand{\CP}{{\mathbb C}P}
\newcommand{\R}{{\mathbb R}}

\renewcommand{\to}{\longrightarrow}
\def\1{\sqrt{-1}\:}
\newcommand{\restrict}[1]{{\left|_{{\phantom{|}\!\!}_{#1}}\right.}}
\newcommand{\cntrct}                
{\hspace{2pt}\raisebox{1pt}{\text{$\lrcorner$}}\hspace{2pt}}

\newcommand{\calo}{{\cal O}}
\newcommand{\cac}{{\cal C}}

\renewcommand{\bar}{\overline}
\renewcommand{\phi}{\varphi}
\renewcommand{\epsilon}{\varepsilon}
\renewcommand{\geq}{\geqslant}
\renewcommand{\leq}{\leqslant}

\newcommand{\Mor}{\operatorname{\mathcal{Mor}}}

\newcommand{\End}{\operatorname{End}}
\newcommand{\Tot}{\operatorname{Tot}}
\newcommand{\St}{\operatorname{St}}
\newcommand{\Id}{\operatorname{Id}}

\newcommand{\const}{\operatorname{\text{\sf const}}}

\newcommand{\Coh}{\operatorname{{\cal Coh}}}
\newcommand{\HHom}{\operatorname{{\cal Hom}}}

\newcommand{\Aut}{\operatorname{Aut}}

\newcommand{\Pic}{\operatorname{Pic}}

\newcommand{\codim}{\operatorname{codim}}

\newcommand{\Sym}{\operatorname{Sym}}
\newcommand{\Proj}{\operatorname{Proj}}

\newcommand{\Spec}{\operatorname{Spec}}
\newcommand{\GL}{\operatorname{GL}}

\newcounter{Mycounter}[section]
\newcounter{lemma}[section]
\renewcommand{\thelemma}{{Lemma \thesection.\arabic{lemma}}}
\newcommand{\lemma}{%
	\setcounter{lemma}{\value{Mycounter}}
	\refstepcounter{lemma}
	\stepcounter{Mycounter}
	{\noindent \bf \thelemma:\ }}

\newcounter{claim}[section]
\renewcommand{\theclaim}{{Claim \thesection.\arabic{claim}}}
\newcommand{\claim}{%
	\setcounter{claim}{\value{Mycounter}}
	\refstepcounter{claim}
\stepcounter{Mycounter}
	{\noindent \bf \theclaim:\ }}

\newcounter{sublemma}[section]
\newcounter{corollary}[section]
\renewcommand{\thecorollary}{{Corollary \thesection.\arabic{corollary}}}
\newcommand{\corollary}{%
	\setcounter{corollary}{\value{Mycounter}}
	\refstepcounter{corollary}
	\stepcounter{Mycounter}
	{\noindent \bf \thecorollary:\ }}

\newcounter{theorem}[section]
\renewcommand{\thetheorem}{{Theorem \thesection.\arabic{theorem}}}
\newcommand{\theorem}{%
	\setcounter{theorem}{\value{Mycounter}}
	\refstepcounter{theorem}
	\stepcounter{Mycounter}
	{\noindent \bf \thetheorem:\ }}

\newcounter{conjecture}[section]
\newcounter{proposition}[section]
\renewcommand{\theproposition} {{Proposition \thesection.\arabic{proposition}}}
\newcommand{\proposition}{%
	\setcounter{proposition}{\value{Mycounter}}
	\refstepcounter{proposition}
	\stepcounter{Mycounter}
	{\noindent \bf \theproposition:\ }}

\newcounter{definition}[section]
\renewcommand{\thedefinition} {{Definition~\thesection.\arabic{definition}}}
\newcommand{\definition}{%
	\setcounter{definition}{\value{Mycounter}}
	\refstepcounter{definition}
	\stepcounter{Mycounter}
	{\noindent \bf \thedefinition:\ }}

\newcounter{example}[section]
\renewcommand{\theexample}{{Example \thesection.\arabic{example}}}
\newcommand{\example}{%
	\setcounter{example}{\value{Mycounter}}
	\refstepcounter{example}
	\stepcounter{Mycounter}
	{\noindent \bf \theexample:\ }}

\newcounter{remark}[section]
\renewcommand{\theremark}{{Remark \thesection.\arabic{remark}}}
\newcommand{\remark}{%
	\setcounter{remark}{\value{Mycounter}}
	\refstepcounter{remark}
	\stepcounter{Mycounter}
	{\noindent \bf \theremark:\ }}

\newcounter{problem}[section]
\newcounter{question}[section]
\makeatletter

\@addtoreset{equation}{section}
\makeatother

\def\blacksquare{\hbox{\vrule width 5pt height 5pt depth 0pt}}
\def\endproof{\blacksquare}

\newcommand{\proof}{{\bf Proof: \ }}
\newcommand{\pstep}{{\bf Proof. Step 1: \ }}

\begin{document}
	
\begin{center}
{\Large\bf  Coherent sheaves on subvarieties \\[3mm] in Hopf manifolds}\\[5mm]

{
Liviu Ornea\footnote{Liviu Ornea is  partially supported by the PNRR-III-C9-2023-I8 grant CF 149/31.07.2023 Conformal
Aspects of Geometry and Dynamics.},  
Misha Verbitsky\footnote{Misha Verbitsky is partially supported by
FAPERJ SEI-260003/000410/2023 and CNPq - Process 310952/2021-2.\\[1mm]
}\\[4mm]

}

\end{center}

{\scriptsize
\hspace{0.13\linewidth}
\begin{minipage}[t]{0.87\linewidth}
{\bf Abstract:} \\
We prove a version of GAGA
theorem for a normal complex analytic variety $X$ equipped with
an invertible holomorphic contraction $\gamma$ with 
center in $x$. We show that $X$ admits a natural structure of an affine variety,
and any $\gamma$-equivariant complex analytic reflexive coherent
sheaf on $X$ admits a natural algebraic structure. We prove a
structure theorem for $X_0:=X\backslash x$, showing that it 
admits a proper action of ${\Bbb C}^*$, and is isomorphic
to the space of non-zero vectors in the total space of
an ample line bundle over the projective variety 
$Z:= X_0/{\Bbb C}^*$ equipped
with an orbifold structure. We show that the quotient
$M:=X_0/\gamma$ admits a holomorphic embedding
to a Hopf manifold, and, conversely, any 
normal subvariety $M$ in a Hopf manifold is obtained this way. 
We prove a form of structure theorem, showing that
any reflexive coherent sheaf on $M$, $\dim M > 2$, 
admits a filtration such that its associated graded 
subquotients, tensored with an appropriate
 line bundle, are obtained as pullbacks of
coherent sheaves on the projective variety 
$Z=X_0/{\Bbb C}^*$. This is used to show that 
any reflexive coherent sheaf on $M$ is 
filtrable, that is, admits a filtration with 
associated graded quotients of rank $\leq 1$.
\end{minipage} 
}

{\scriptsize
\tableofcontents
}

\section{Introduction}


\subsection{Cone varieties and subvarieties in Hopf manifolds}

In this paper we study 
germs of complex analytic varieties
equipped with a holomorphic invertible contraction.
Such objects appeared in several geometric
situations, often independently. 

A formal definition of a contraction of $\check X$ to $x\in \check X$ is given 
in \ref{_contraction_cone_Definition_};
it is a map $\gamma:\; \check X \to \check X$ such 
that $\gamma^n$ converges to the constant
map $\check X \mapsto x$ uniformly on compacts, when  $n\to \infty$.

The aim of this paper is to show that
a contraction $\gamma$ on a complex analytic variety $X$
provides it with a natural algebraic structure.
This algebraic structure is used to treat the
$\gamma$-equivariant coherent sheaves on $X$.
In the end, we obtain results about filtrability
of coherent sheaves on complex subvarieties
in a Hopf manifold.

H. Poincar\'e (\cite{_Poincare:Thesis_}, 1879) and H. Dulac
(\cite{_Dulac_}, 1912) were the first to realize that
a contraction vector field has a polynomial normal form.
Later, S. Sternberg (\cite{_Sternberg_contraction_}, 1957) has 
generalized the Poincar\'e-Dulac theorem to a 
contraction acting on a manifold. His results 
were used by K. Kodaira in his classification
of Hopf surfaces (\cite{_Kodaira_Structure_III_}).

In his dissertation, Poincar\'e proved the general case of his theorem,
showing that any contraction (smooth or holomorphic)
is linear in appropriate coordinates, assuming
the eigenvalues of the differential of the 
contraction are general (the so-called ``non-resonant case'').
After that, Dulac found a polynomial normal form also 
for the resonant case. 

Curiously enough, the methods used in this
paper and other works on locally conformally K\"ahler
geometry can be applied to obtain a geometric
proof of Poincar\'e-Dulac theorem. We first noticed
this in \cite{_OV:Mall_}, where we have shown that
the Poincar\'e-Dulac theorem follows from the
vanishing of the cohomology of a certain holomorphic
vector bundle on a Hopf manifold. This way we
proved the non-resonant case of Poincar\'e-Dulac;
the resonant case was proven in \cite{_Madera:Hopf_}, 
using the same approach.

In the present paper, we study the holomorphic
contractions acting on a variety which is possibly 
singular, and generalize the
Poincar\'e-Dulac theorem in the following form: if $\gamma$ is a
contraction acting on a complex analytic variety $\check X$ 
then $\check X$ admits an algebraic structure, 
such that $\gamma$ is algebraic. This is
similar to Poincar\'e-Dulac constructing
a polynomial normal form for a contraction.

We show that $(\check X, \gamma)$ is 
a germ of a Stein variety (\ref{_closed_cone_Stein_Theorem_}), which is uniquely
determined by the action of $\gamma$, and admits a canonical algebraic structure
making it an affine variety  
(\ref{_algebraic_str_on_cone_uniqueness_Theorem_},
\ref{_normalization_closed_cone_algebraic_Remark_}).
In the proof, we assume normality of $\check X$,
to avoid unnecessary complications.

This research started in 2000-ies, when we 
were exploring the geometry of complex submanifolds
of a Hopf manifold. In \cite{_OV_lckpot_} 
we discovered that the submanifolds of Hopf manifolds
have an intrinsic characterization in terms of
a special class of Hermitian metrics; we call
them ``LCK metrics with potential''.

An LCK (locally conformally K\"ahler)
manifold is a complex Hermitian manifold $(M, I, \omega)$
such that $d\omega =\theta \wedge \omega$, where $d\theta=0$. The definition 
was proposed by I. Vaisman who later proved (\cite{_Vaisman_trans_}) that on
compact K\"ahler manifolds, any LCK metric is in fact
globally conformally K\"ahler (that is, $\theta$ is exact,
$\theta=df$, and $e^{-f}\omega$ is K\"ahler).

An equivalent definition says that a complex manifold
$(M,I)$ is LCK if and only if it admits a K\"ahler cover
$(\tilde M, I, \tilde\omega)$ on which the deck group
$\Gamma$ acts by holomorphic homotheties:
$\gamma^*\tilde\omega=c_\gamma\tilde\omega$, where
$c_\gamma\in\R^{>0}$, for all $\gamma\in \Gamma$ (see
\cite[\S 3.4]{_OV_book_}. When the metric $\tilde\omega$
has a global potential,
$\tilde\omega=dd^c\phi$, on which the deck group acts by
scalar multiplication,
$\gamma^*\phi=c_\gamma\phi$, $(M,I)$ is called LCK with
potential\footnote{Such a potential is called
  ``automorphic''.}. Unlike the general LCK metrics, LCK
metrics with potential are stable at small deformations of
the complex structure (see \cite{_OV_lckpot_},
\cite[Chapter 12]{_OV_book_}). 

The first example of LCK metrics with potential appeared
implicitly in \cite{_GO_Hopf_surfaces_}, on Hopf
surfaces. Later on, we proved that {\em all} Hopf
manifolds, i.e. quotients of $\C^n\backslash 0$ by a
cyclic group generated by an invertible holomorphic
contraction at 0, linear or non-linear, admit LCK metrics
with potential
(\cite{_OV_lckpot_,_OV_non_linear})\footnote{If the
  contraction is linear, we speak about ``linear Hopf
  manifolds''.}.

Clearly, all complex submanifolds of
an LCK manifold with potential inherit
this structure. Conversely, we proved in
\cite{_OV_lckpot_} (see also \cite[Chapter 13]{_OV_book_})
that a compact LCK manifold with potential of dimension at
least three is biholomorphic to a submanifold of a
linear Hopf manifold. All in all, we know that a compact
LCK manifold of dimension at least three is LCK with
potential if and only if it is a submanifold of a linear
Hopf manifold.

The dimensional bound is due to 
the technicalities of our proof.
We deduce the results on LCK manifolds
with potential from the  fundamental theorems
of complex analysis which bear the same dimensional restriction.
The proof of the embedding result 
for LCK manifolds with potential was based on 
Cartan's theorem A and B and Grauert's solution
of the Levi problem, (and, more specifically,
on Rossi and Andreotti-Siu's theorem, 
\cite{_Andreotti_Siu:embeddings_}) 
used to fill a pseudoconcave hole
in a complex manifold with a Stein 
domain with a compact boundary.
We don't know whether our results
on LCK manifolds with potential
need this dimensional restriction.

An LCK manifold  with potential
always admits a holomorphic embedding to
a Hopf manifold, $M \hookrightarrow \frac{\C^n\backslash
  0}{\langle A \rangle}$. Consequently,
$M$ admits a $\Z$-covering which is equipped with a 
$\Z$-equivariant embedding to $\C^{n}\backslash 0$. At this point
the study becomes focused on $\Z$-invariant
submanifolds $X_0\subset \C^{n}\backslash 0$,
where $\Z$ acts on $\C^n$ by linear contractions.
We have shown that $X_0$
admits a proper $\C^*$-action, containing a holomorphic
contraction, with $Z:=X_0/\C^*$ a projective orbifold
(\cite{_OV_Algebraic_Cones_}). In \cite{_OV_Algebraic_Cones_},
we defined {\bf an open algebraic cone} as a complex
manifold equipped with a $\C^*$-action which
contains a contraction and admitting a $\C^*$-invariant
embedding to $\C^n$, where the $\C^*$-action on $\C^n$ 
is linear and has eigenvalues $|\alpha_i| <1$.

By Remmert-Stein theorem,
the closure of an open algebraic cone in $\C^n$
is a complex subvariety; we called it
{\bf a closed algebraic cone}.
In \cite{_OV_Algebraic_Cones_},
we have shown that a complex analytic variety
equipped with a holomorphic contraction
is biholomorphic to a closed algebraic cone if its
only singularity is the attraction
point of this contraction (called
{\em the apex} or {\it the origin}
later on). 

Moreover, in \cite{_OV_Algebraic_Cones_}
we have shown that an open algebraic cone
$X_0$ is biholomorphic to the total space
$\Tot^\circ(L)$ of non-zero vectors in 
an ample line bundle on the projective orbifold $Z$. 
This argument is used to obtain a structure
theorem for LCK manifolds $M$ with potential,
$\dim_\C M>2$.\footnote{The same result is true
for complex surfaces if we assume the so-called 
GSS conjecture, \cite{ov_indam}.}
This argument also 
implies that any LCK manifold with 
potential of dimension $>2$ is biholomorphic to $\Tot^\circ(L)/\langle\gamma\rangle$,
where $\gamma$ denotes an automorphism of $Z$ which
equivariantly acts on the ample line bundle $L$
by Hermitian homotheties.

In the present paper we ignore the metric
aspect of this theory, and focus on its complex analytic
and algebro-geometric aspects. Among other things,
this allows us to  greatly relax
the assumptions.

The central notion of the present
paper is {\bf the cone variety}
(\ref{_contraction_cone_Definition_}).
Just like it was with the algebraic
cones, it comes in two flavors: 
{\bf a closed cone variety} $X$ is
a complex analytic variety of dimension $>1$ equipped with
a holomorphic contraction $\gamma$, and
{\bf an open cone variety} $X_0$ is
a closed cone variety without the origin.

Clearly, the contraction $\gamma$ acts
on $X_0$ properly discontinuously.
We prove that the quotient $X_0/\langle\gamma\rangle$
admits a holomorphic embedding to a Hopf manifold.
Conversely, any complex subvariety $M$ in a Hopf
manifold, $\dim_\C M >0$, is obtained as 
 $X_0/\langle\gamma\rangle$ for some open cone
variety $(X_0, \gamma)$ (\ref{_cone_embedded_Theorem_}).

To simplify the statements, we generally
assume that our cone varieties are normal.
This assumption can be avoided, sacrificing 
the clarity.

Similarly to algebraic cones,
a cone variety is biholomorphic to the
total space of non-zero vectors in 
an ample line bundle over a projective
orbivariety. 

The orbivarieties, which are known in algebraic
geometry as Deligne-Mumford stacks (\cite{_Calle:Stacks_}), are
(roughly speaking) complex analytic varieties which are locally obtained
as finite quotients of complex varieties,
with the finite group action still attached.
We provide a rather comprehensive introduction to this
theory in Section
\ref{_orbivarieties_equiv_sheaves_Section_}.

For the present purpose, this notion cannot
be avoided, because the ample line bundle $L$ on an
orbivariety $Z$, used
to identify an open cone variety $X_0$ with $\Tot^\circ(L)$,
is an orbibundle. Indeed, for each open chart $U/G\subset Z$
where $U$ is a complex variety,  the line bundle on $U/G$
is a $G$-equivariant bundle on $U$, and the stabilizer 
$\St_x(G)$ of $x\in U$ can act non-trivially on the fiber
of $L$ in $x$. Then the total space $\Tot(L)$
is identified locally with 
the quotient of the total space of a line bundle on $U$ by
$G$, and it is not locally trivial over $U/G$,
because its general fiber is $\C^*$ and
its special fiber in $x$ is $\C^*/\St_x(G)$.

This problem is purely technical, and at no 
point we use anything specific to orbivarieties
or orbibundles. Indeed, for most part, we
prove our results ignoring the orbivariety
aspect, to make the exposition more comprehensible.

\subsection{GAGA theorems for cone varieties}

Throughout this paper, ``variety'' refers
to a complex analytic variety, unless indicated otherwise.
However, in most cases the complex analytic varieties
we consider are equipped with an algebraic structure,
or obtained as $\Z$-quotients of complex algebraic varieties.

In the present 
paper, we prove several versions of the Serre's 
GAGA for the
complex varieties equipped
with a contraction.

Results which are similar to ours appeared
in the literature, usually motivated by 
research in singularities. In singularity theory,
a variety equipped with an action of $\C^*$
which contains a contraction is called {\bf homogeneous}
or {\bf quasi-homogeneous}. Since the category
of representations of $\C^*$ is equivalent to
the category of graded vector spaces, 
a quasi-homogeneous singularity is defined
by polynomial equations for an appropriate 
coordinate system. This result was first
formulated by P. Orlik and Ph. Wagreich
(\cite[Proposition 1.1.2]{_Orlik_Wagreich_}).
Since then, ``quasi-homogeneous singularity'' was 
understood alternatively as the one defined by 
quasi-homogeneous polynomial equation 
or as a $\C^*$-invariant singularity.\footnote{
A quasi-homogeneous polynomial of weight
$w_1, ..., w_n$, where $w_i \in \Z^{>0}$, 
is a polynomial $P(t_1, ..., t_n)$ obtained as a sum
of monomials of degrees $d_1, ..., d_n$
such that $\sum w_i d_i=\const$.}

The earliest research on
quasi-homogeneous singularities is due
to E. Brieskorn and to K. Saito, who was
Brieskorn's graduate student at the moment.
In his first paper \cite{_Saito:quasihomo_},
Saito gave an intrinsic characterization
of a quasi-homogeneous hypersurface singularity, 
showing that an isolated singularity
is quasi-homogeneous if and only if
its Tyurina number is equal to its Milnor number.
This result also gives an algebraization theorem,
because the $\C^*$-action allows to define the
singularity using the polynomial equations 
(\cite[Proposition 1.1.2]{_Orlik_Wagreich_}).
Saito's theorem, proven for hypersurfaces,
 was generalized to all isolated
 quasi-homogeneous singularities by H. Vosegaard
(\cite{_Vosegaard_}).

Since quasi-homogeneity implies algebraicity,
it is natural to expect that the existence
of a contraction implies the quasi-homogeneity,
and this would give the algebraization theorem.

The current argument goes in this direction,
but our original motivation was different.
Motivated by differential geometry, we have
been studying a variety $X$ with an isolated singularity
admitting a holomorphic contraction $\gamma$,
showing that $(X, \gamma)$ admits an equivariant
embedding to $(\C^n, A)$, where $A$ is a linear
contraction. In \cite{_OV_lckpot_} (see also 
\cite{ov_indam}), it was shown that
such an embedding can always be constructed,
when $\gamma$ acts by homotheties on a K\"ahler
metric, and in \cite{_OV_pams_} it was shown that
such a K\"ahler metric always exists,
making the quotient manifolds 
$\frac{\C^n\backslash 0}{\langle A \rangle}$ and 
$\frac{X\backslash x}{\langle \gamma \rangle}$
locally conformally K\"ahler.

It is not hard to show that an equivariant embedding
to $(\C^n, A)$ (or, equivalently, a holomorphic
embedding from the LCK manifold
$\frac{X\backslash x}{\langle \gamma \rangle}$
to a linear Hopf manifold) equips $(X, \gamma)$
with an algebraic structure. Indeed, a holomorphic
function $f$ on $\C^n$ is polynomial if and only if
the space $\langle f, A^*(f), A^*A^*(f), A^*A^*A^*(F),
... \rangle$ is finite-dimensional 
(\cite[Lemma 6.1]{_OV_Algebraic_Cones_}). 
This proves the existence of an algebraic structure
on $X$ compatible with $\gamma$. Its uniqueness
follows from \ref{_algebraic_str_on_cone_uniqueness_Theorem_},
though a weaker version of this result was
shown in \cite{_OV_Algebraic_Cones_}.

In the literature, results of a similar nature
appeared from a completely different angle.
In \cite{_Favre_Ruggiero_}, C. Favre and M. Ruggiero
prove an algebraization result in dimension 2,
showing that any 2-dimensional variety admitting 
a holomorphic contraction is also quasi-homogeneous,
that is, admits a $\C^*$-action with positive weights.
Using Kodaira's classification of compact complex surfaces,
Favre-Ruggiero also obtained a classification theorem
for 2-dimensional varieties admitting a contraction, 
showing that the corresponding complex surface is Hopf or
elliptic (in \cite{_ovv:surf_}, the same classification result
was obtained using the LCK geometry).

Another algebraization theorem was obtained by K. Morvan
in \cite{_Morvan_}. Morvan has shown that any complex
singularity $X$ admits a local embedding to $\C^n$ such that
any automorphism of $X$ is extended to a local
automorphism of $\C^n$. Then he applies the
Poincar\'e-Dulac theorem to show that a local 
contraction of $X$ is polynomial in appropriate
coordinates.

In \cite{_Serre:GAGA_}, J.-P. Serre coined the term ``GAGA
principle'' referring to the whole spectrum of
algebraization results, showing that the complex
subvarieties of projective varieties are algebraic,
complex morphisms of projective varieties are algebraic,
and coherent sheaves on projective varieties are
algebraic, and, moreover, the algebraic structure
is unique and functorial.

In the present paper we aim to prove the same
for complex varieties equipped with a contraction,
that is, for the cone varieties.

In fact, there are three categories which are relevant
for our present work. Let $(X, \gamma, x)$ be a (closed) cone variety,
where $x\in X$ denotes the attraction point of $\gamma$.
Recall that the complement $X_0:=X \backslash x$ is called 
an open cone variety. Clearly, $\gamma$ acts on
 $X_0$ properly discontinuously, and the quotient
$\frac{X \backslash x}{\langle \gamma\rangle}$
is a complex analytic variety. Such a variety admits
a holomorphic embedding to a Hopf manifold
(\ref{_cone_embedded_Theorem_}), and, conversely, any subvariety
of a Hopf manifold is obtained this way. 

The three categories we referred to were 
the category of closed cone varieties, 
open cone varieties, and varieties admitting
an holomorphic embedding to a Hopf manifold.

The natural functor $(X,\gamma, x) \mapsto (X_0, \gamma)$
is not an equivalence, unless $X$ is normal.
The geometry which is hidden behind this
observation is rich and complicated 
(\cite{_OV_Algebraic_Cones_}), and 
we will ignore it for the present paper.

The functor taking $(X_0, \gamma)$
to $M:=\frac{X_0}{\langle \gamma\rangle}$
also loses information, unless $\pi_1(M)=\Z$.
However, we can recover $X_0$ if we fix a choice of an infinite cyclic
subgroup in $\pi_1(M)$.

The variety $M$ should be understood as a 
singular version of a locally conformally K\"ahler manifold
with potential (indeed, if $M$ is smooth, it is LCK with
potential). Sometimes such varieties are called ``LCK
varieties with potential''.

The GAGA principle can be formulated for open and closed
cones; in both cases it states the existence and
uniqueness of an algebraic structure compatible
with the action of $\gamma$. We prove it
for open cone varieties in
\ref{_algebraic_str_on_cone_uniqueness_Theorem_}
(by functoriality of the normalization,
this also takes care for GAGA for normal 
closed cone varieties).

We did not go further in this direction,
leaving the non-normal varieties firmly out
of the scope of the present paper.

It is hard to interpret any of our results
as GAGA or algebraization of the subvarieties
in a Hopf manifold. Indeed, such varieties are
not algebraic in any conventional sense
of this word. However, the GAGA principle
for coherent sheaves or holomorphic vector bundles,
is manifest in all three categories. Indeed,
our principal motivation was to understand
the geometry of holomorphic vector bundles
on the LCK manifolds with potential.

The category of coherent sheaves on 
$M:=\frac{X_0}{\langle \gamma\rangle}$
is tautologically equivalent to the
category of $\langle\gamma\rangle$-equivariant
coherent sheaves on the open cone $X_0$.
Speaking of GAGA-type or algebraization
theorems for coherent sheaves on cone
varieties, we always talk about
$\langle\gamma\rangle$-equivariant
coherent sheaves (for the definition
and introduction to equivariant
coherent sheaves, see Subsection 
\ref{_equivariant_coherent_Subsection_}).

The GAGA theorem for coherent sheaves
on the closed cones is the most straightforward.
Let $(X,\gamma)$ be a closed cone variety,
$\Coh^{an}_\gamma(X)$ the category of 
complex analytic $\langle\gamma\rangle$-equivariant
coherent sheaves, and $\Coh^{alg}_\gamma(X)$ the category of 
algebraic $\langle\gamma\rangle$-equivariant
coherent sheaves. Then the identity functor
$\Coh^{alg}_\gamma(X)\to \Coh^{an}_\gamma(X)$
is an equivalence of categories (\ref{_coherent_GAGA_Theorem_}).

However, for our applications we are more
interested in the category of $\langle\gamma\rangle$-equivariant
coherent sheaves on the open cone variety $(X_0, \gamma)$,
and here the equivalence becomes more subtle.
The main problem is dimension 2: even a vector
bundle ($\langle\gamma\rangle$-equivariant or not)
on a 2-dimensional cone $(X_0, \gamma)$
such as $(\C^2\backslash 0, v\mapsto 2v)$
does not necessarily extend to a coherent sheaf
on the closed cone.

When the extension exists, the equivariance
is not always imposed automatically either.
The reflexivity and the normality of a sheaf are the key notions here.

A coherent sheaf ${\cal F}$ on $M$ is called {\bf reflexive}
if the natural map ${\cal F}\to {\cal F}^{**}$
is an isomorphism, where ${\cal F}^*:= \HHom({\cal F}, \calo_M)$.

A coherent sheaf on a normal variety $M$
is called {\bf normal} if for any subvariety
$Z\subset M$ of codimension $\geq 2$,
the natural map ${\cal F} \to j_* j^* {\cal F}$
is an isomorphism, where $j:\; (M\backslash Z) \to M$
is the open embedding map. 
If $j_* {\cal F}$ is coherent, we call it {\bf the standard extension}
of a coherent sheaf ${\cal F}$ from $M\backslash Z$ to $M$.
When $M$ is normal itself,
normality is automatic for any vector bundle
(this is essentially the Hartogs theorem).
Indeed, the normality in this situation is
equivalent to reflexivity (\ref{_normal_is_reflexive_Theorem_}).

From this observation it is clear that 
the standard extension of a reflexive
$\langle\gamma\rangle$-equivariant coherent sheaf from 
$(X_0, \gamma)$ to $X$ is also normal and 
$\langle\gamma\rangle$-equivariant.

When $\dim X \geq 3$, the extension 
of any reflexive sheaf from $X_0:= X \backslash x$
to $X$ exists and is reflexive
 (\ref{_Siu_over_point_Corollary_}). This gives the equivalence
of categories between the reflexive 
$\langle\gamma\rangle$-equivariant coherent sheaves on $X$ and $X_0$
(\ref{_categories_on_cone_Proposition_})
and brings the GAGA theorem to the category of
$\gamma$-equivariant reflexive sheaves on $X_0$.
For $\dim X \leq 2$, this equivalence fails,
as we explain in Subsection \ref{_filtrability_Subsection_}.
This puts a stop on attempts to produce such an
equivalence for coherent sheaves which are not necessarily
reflexive,  even when
$\dim X \geq 3$, because the category
of $\langle\gamma\rangle$-equivariant coherent
sheaves supported on 2-dimensional $\gamma$-invariant subvariety
is naturally embedded to $\Coh_\gamma(X)$.

The GAGA-type algebraization results in formal setting
were explored in Sections \ref{_algebraiz_for_rings_Section_} 
and \ref{_modules_over_complete_Section_},
which are somewhat isolated from the rest of this paper.
In Section \ref{_algebraiz_for_rings_Section_}, we prove
that a complete local ring $R$ equipped with a contraction
automorphism $\gamma$ admits a natural algebraic structure
(that is, a $\gamma$-invariant finitely generated subring $R_{alg}$  which has
the same completion). In Section \ref{_modules_over_complete_Section_},
we prove GAGA for $\langle\gamma\rangle$-equivariant 
finitely generated $R$-modules, showing that
every such module is obtained from a finitely generated $R_{alg}$-module
by tensoring with $R$. These results are more
elementary and straightforward, but they are
necessary for the proof of the GAGA theorem for
coherent sheaves (\ref{_formal_used_here_Remark_}).

\subsection{Coherent sheaves on complex manifolds and filtrability}
\label{_filtrability_Subsection_}

This paper explores a phenomenon which is specific
for dimension $\geq 3$. Indeed, there is much literature
about vector bundles on complex surfaces, where
it is shown that they behave in entirely different
manner from those on Hopf manifolds in dimension $\geq 3$
and their subvarieties.

The key notion here is {\em filtrability} of coherent sheaves.
A coherent sheaf ${\cal F}$ is called {\bf filtrable} if it
admits a filtration 
$0={\cal F}_0 \subset {\cal F}_1 \subset ...\subset {\cal F}_n={\cal F}$ 
by coherent subsheaves, with
all subquotients ${\cal F}_i/{\cal F}_{i-1}$ of rank $\leq 1$.
Filtrability is automatic for coherent sheaves
on projective varieties (the argument used in
\ref{_filtrability_for_G_sheaves_Claim_} can be applied).
It is clearly false for non-algebraic K\"ahler
manifolds; indeed, for a general deformation
of a hyperk\"ahler manifold $M$ with maximal holonomy, any non-trivial
subsheaf of its tangent bundle $TM$
has maximal rank (see {\em e.g.} \cite{_KV1_}).

For non-algebraic complex surfaces, filtrability 
is known to be false. Indeed, most examples of non-algebraic
complex surfaces are known to admit {\bf irreducible vector
bundles,} that is, holomorphic vector bundles without non-trivial
coherent subsheaves of smaller rank.

The earliest result in this direction is probably
due G. Elencwajg and O. Forster
(\cite{_Elencwajg_Forster_}), who proved the existence
of irreducible holomorphic vector bundles on 
2-dimensional complex tori with trivial N\'eron-Severi
group. A few years later, C. B\u anic\u a and 
J. Le Potier (\cite{_Banica_Le_Potier_})
generalized their work, proving the 
existence of irreducible holomorphic vector bundles
on any complex surface of algebraic dimension 0,
and on Kodaira surfaces.

The classification of non-K\"ahler surfaces is based
on the so-called  GSS conjecture, which is widely accepted
by many mathematicians. It claims that a minimal surface $M$
of Kodaira dimension $-\infty$, $b_2>0$ and $b_1=1$ admits
a real 3-sphere $S \subset M$ such that
$M\backslash S$ is connected, and a neighbourhood
of $S$ is biholomorphic to a neighbourhood of 
the standard 3-sphere in $\C^2$. Such a sphere
is called {\bf a global spherical shell} (GSS).
Surfaces admitting a global spherical shell are
called Kato surfaces, after Masahide Kato; also, 
the GSS conjecture is often called Kato's conjecture
(\cite{_Dloussky:Kato_}).

In \cite{_ovv:surf_}, it was shown that, assuming
the GSS conjecture, all non-K\"ahler
surfaces except one class of Inoue surfaces and its blow-ups
are locally conformally K\"ahler.

The minimal models of
all of these surfaces, except Inoue and Kato, are LCK
with potential and admit a holomorphic embedding 
to a Hopf manifold (\cite{_OV_book_}). 
Therefore, they can be obtained as a 
quotient of a smooth open cone variety $(X_0, \gamma)$
by $\langle \gamma\rangle$. Coherent sheaves
on $\frac{X_0}{\langle \gamma\rangle}$ are
the same as $\gamma$-equivariant coherent sheaves
on $X_0$; however, these sheaves do not always
admit an extension to the corresponding closed
cone $X$, and this is where the rest of our
results fail and cannot be applied anymore.

Indeed, in contrast with the results of the present paper
(which are valid in dimension $\geq 3$),
filtrability is known to fail for all 
compact complex surfaces which admit an 
LCK metric with potential, that is,
for Hopf surfaces and non-K\"ahler
elliptic surfaces.

All non-K\"ahler surfaces are either
elliptic or have algebraic dimension 0 (\cite{_BHPV_}).
For Hopf surfaces of algebraic dimension 0 and
Kodaira surfaces, the existence of non-filtrable bundles is already covered
by B\u anic\u a and Le Potier (\cite{_Banica_Le_Potier_}).

Existence of non-filtrable holomorphic bundles
on elliptic surfaces with the base of genus $>1$
seems to be mostly open even now, though 
in \cite{_Brinzanescu_Moraru_}, V. Br\^inz\u anescu and R. Moraru 
address this question and construct 
infinitely many non-filtrable holomorphic bundles
on elliptic surfaces such that their relative
Jacobian admits an irreducible bisection
(\cite[Remark after Proposition 4.5]{_Brinzanescu_Moraru_}).
In \cite{_Aprodu_Toma_}, M. Aprodu and M. Toma have shown that
all non-filtrable bundles of rank 2 on elliptic
surfaces are obtained from this construction.

In dimension $>2$, the situation is rather different.
Most research in this direction is based on the notion
of stability, which is defined for non-K\"ahler complex
manifolds in terms of the Gauduchon metrics (\cite{_Lubke_Teleman:Book_}).
Any coherent sheaf admits the Harder-Narasimhan
filtration, with all associated graded subquotients
semistable (\cite{_Wang:Harder-Narasimhan_}).
Any semistable sheaf admits the Jordan-H\"older
filtration with all associated graded subquotients
stable of the same slope (\cite{_Clarke_Tipler_}).
This reduces the question of filtrability 
to stable bundles and sheaves.

For stable bundles on Hopf manifolds of
Vaisman type (that is, the diagonal Hopf manifold,
which are obtained as quotients of $\C^n$ by
a diagonalizable contraction), the filtrability is
proven in \cite{_Verbitsky:Filtrable_} using the Kobayashi-Hitchin
correspondence associating to each stable bundle
(or coherent sheaf) a certain special connection.
The same argument was applied to stable bundles on
manifolds of dimension $> 2$ admitting positive elliptic fibrations
 (\cite{_Verbitsky:Stable_}).\footnote{
Consider an elliptic fibration $\pi:\; M \to X$.
It is called {\bf positive} if
$\pi^*\omega$ is exact for a K\"ahler form $\omega$ on $X$.}
This covers elliptic Hopf manifolds, Calabi-Eckmann
manifolds,as well as some
of invariant complex structures on Lie groups, such as 
$\mathrm{SU}(3)$.

This was the starting point of the present paper, though
in its current version we do not use either stability
or Yang-Mills connections. Instead, we prove the same
conclusion as obtained in \cite{_Verbitsky:Stable_},
using the GAGA principle in place of the Kobayashi-Hitchin
correspondence.

We prove the filtrability of reflexive coherent sheaves
on any normal variety $M$, $\dim_\C M\geq 3$, which admits
a holomorphic embedding in a Hopf manifold 
(by \ref{_cone_Hopf_embedding_Theorem_}, this is equivalent
to being biholomorphic to $\frac{X_0}{\langle\gamma\rangle}$,
where $(X_0, \gamma)$ is an open cone variety).
We use the equivalence between the category
of coherent sheaves on $M$ and the category
of $\gamma$-equivariant coherent sheaves
on $X_0$. Applying Siu's extension theorem,
we show that the category of $\gamma$-equivariant coherent sheaves
on $X_0$ is equivalent to the category of $\gamma$-equivariant coherent sheaves
on the corresponding closed cone variety $(X, x, \gamma)$
(\ref{_categories_on_cone_Proposition_}).

The action of $\gamma$ on $X$
is {\bf strongly algebraic} (\ref{_strongly_alg_Definition_}), that is, factorized
through the action of an algebraic group ${\cal G}$,
which can be identified with the Zariski closure
of $\langle \gamma\rangle$ (\ref{_algebraic_str_on_cone_uniqueness_Theorem_}).
Since ${\cal G}$ is an abelian algebraic
group which acts on the Zariski tangent space 
$T_x X$ by contractions, it contains a subgroup
 $\C^*\subset {\cal G}$ acting on $T_x X$ 
by contractions. In \ref{_closed_cone_alg_cone_Theorem_}, 
we use this $\C^*$-action
to identify $X$ with an affine cone,
that is, the spectrum of the graded ring
$\bigoplus_i H^0(Z,L^{\otimes i})$, where
$Z:=X_0/\C^*$ is a projective orbivariety, and $L$
an ample orbibundle of rank 1 on $Z$.
The choice of $\C^*\subset {\cal G}$
is not unique; indeed, it is possible to end
up with non-isomorphic $Z$ for different choices of $\C^*$.
We fix the choice of $\C^*\subset {\cal G}$
attached to $Z$, and denote this subgroup by $\C^*_Z$.

In this situation, the GAGA theorem  shows
that the category of $\gamma$-equivariant analytic 
coherent sheaves on $X$ is equivalent to the category
of $\gamma$-equivariant algebraic
coherent sheaves on $X$. 
In \ref{_G_F_equivariant_on_F_Proposition_} we show that
the $\gamma$-action on an equivariant
coherent sheaf ${\cal F}$ on $X$ is naturally extended
to the action of an algebraic group ${\cal G}_{\cal F}$,
which is projected to ${\cal G}$ surjectively.
If this projection is an isomorphism,
the sheaf ${\cal F}$ is $\C^*_Z$-equivariant,
and we can identify it with the pullback
of the coherent sheaf ${\cal F}^{\C^*_Z}$ on $Z$
(\ref{_Roberts_invariant_sections_Theorem_});
this would insure the filtrability, because
$Z$ is projective, and ${\cal F}^{\C^*_Z}$ is 
a coherent sheaf on a projective orbivariety.

However, not every coherent sheaf is obtained this
way. This phenomenon was first discovered 
on principal elliptic fibrations  
(\cite{_Verbitsky:Stable_}, 
and extended to diagonal Hopf manifolds in
\cite{_Verbitsky:Filtrable_}). For elliptic Hopf manifolds, which are
elliptically fibered over a weighted projective
space, $\pi:\;M \to W \C P^n$,
it can be described as follows. For these examples,
the group $\C^*_Z$ acts fiberwise on the elliptic fibers
of $\pi$, and its action is factorized through the principal
action of the elliptic curve on the fibers of the
elliptic fibration. Therefore, the sheaves
which are lifted from $Z=W \C P^n$ are
trivial on the elliptic curves. However,
it is not hard to construct a coherent sheaf
(even a line bundle) which is non-trivial
and non-torsion on each of the elliptic
fibers of $\pi$.

In \cite{_Verbitsky:Stable_,_Verbitsky:Filtrable_},
this problem was resolved using the Kobayashi-Hitchin
correspondence. The Yang-Mills connection which 
is obtained from the Donaldson-Uhlenbeck-Yau theorem
is flat on the leaves of the elliptic fibration $\pi$
and has scalar monodromy. Tensoring a stable coherent sheaf
with a line bundle $L^{-1}$, where $L$ has the same monodromy
on the fibers of $\pi$, leads to a coherent sheaf which
is trivial on the fibers of $\pi$, hence lifted from $\Z$.

In the present paper we discover that the same
conclusion can be reached without using the Kobayashi-Hitchin
correspondence (which cannot be applied to the singular
open cone variety anyway). We prove that any reflexive
$\gamma$-equivariant coherent sheaf ${\cal F}$ on a closed cone variety
$X$ admits a filtration with associated graded quotients
${\cal F}_i$ which become $\C^*_Z$-equivariant after
tensoring with a line bundle 
(\ref{_filtra_with_subquot_lifted_up_to_L_i_Theorem_}).
This proves the filtrability of ${\cal F}$, because $\C^*_Z$-equivariant 
sheaves are lifted from $Z$, which is projective
(\ref{_filtra_with_subquot_lifted_up_to_L_i_Theorem_}).

In most of the paper, we assume that the dimension
of the cone variety is $\geq 3$ and the coherent sheaves
are reflexive. However, there are cases when these assumptions are
excessive and can be avoided.

The main obstruction in extending these results to 
dimension 2 is the failure of Siu extension theorem
(\ref{_Siu_extending_Theorem_}) in this dimension. 
Indeed, there are non-filtrable vector bundles on 
Vaisman and Hopf surfaces, as explained earlier
in this subsection. A more subtle difficulty
of the same nature is the requirement of reflexivity:
if $K$ is the kernel of a surjective map
$F \to F_1$, where $F_1$ has support of dimension 2,
it is not clear how to extend $K$, even if $F$ is extensible.
We do not have counterexamples in this situation, 
and theoretically our results could be proven
for some (or all) torsion-free sheaves.

We also use the reflexivity to extend
the Poincar\'e-Dulac theorem to construct an
equivariant structure on coherent sheaves
in formal situation
(\ref{_algebraization_R_modules_Theorem_}). 
With the exception of this result (which is
used only once), everything related to closed
cone varieties does not need either reflexivity
or the dimension $\geq 2$ assumption.

\subsection{Main results}

This subsection serves as a reference list
collecting the main results of this paper.

\hfill

In Sections \ref{_algebraiz_for_rings_Section_} 
and \ref{_modules_over_complete_Section_}, we explore
complete local rings equipped with contractions
and prove the GAGA-type theorems. 
Let ${\goth m}$ be a maximal ideal in a local ring $R$.
{\bf A contraction} on a local ring $R$
is an automorphism which acts on ${\goth m}/{\goth m}^2$
with eigenvalues $|\alpha_i|<1$.
Let $\gamma$ be an endomorphism of a vector space $V$. A vector
$v\in V$ is called {\bf $\gamma$-finite} if
it is contained in a finite-dimensional $\gamma$-invariant
subspace of $V$.

\hfill

\theorem
Let $R$ be a complete local ring of finite type, 
and $\tau\in \Aut(R)$ a contraction. 
Then the ring $R^\tau$ of $\tau$-finite elements in $R$ is dense
in $R$ and finitely generated.

\proof \ref{_algebraization_Theorem_}. \endproof

\hfill

The next theorem deals with an algebraic analogue
of equivariant coherent sheaves. Let $\tau$ be an
automorphism of a ring $R$. A {\bf $\tau$-equivariant
structure} on an $R$-module $V$ is a map $\tilde\tau:\; V \to V$
such that $\tilde\tau(av) = \tau(a) \tilde \tau (v)$
for any $v\in V, a \in R$.

\hfill

\theorem
Let $R$ be a complete local ring of finite type without zero divisors,
$\tau\in \Aut(R)$ a contraction, and $V$ a finitely generated
reflexive $R$-module. Let $\tilde\tau:\; V \to V$ be the  
structure map of a $\langle\tau\rangle$-equivariant module on $V$.
Denote by $R^\tau$, $V^{\tilde\tau}$ the corresponding
spaces of $\tau$-finite and $\tilde\tau$-finite elements.
Then $V^{\tilde\tau}$ is dense in $V$. Moreover, $V^{\tilde\tau}$
is finitely generated as an  $R^\tau$-module.

\proof \ref{_algebraization_R_modules_Theorem_}. \endproof

\hfill

In the next two sections, \ref{_cone_in_Hopf_Section_} 
and \ref{_GAGA_for_cones_Section_},
we prove the GAGA theorem for cone varieties.
We start with an embedding theorem, which 
parallels the results of \cite{_OV_lckpot_}
where the same theorem was proven 
for LCK manifolds with potential.

\hfill

\theorem
Let $H$ be a Hopf manifold, and $M \subset H$ a 
complex subvariety. Denote by $\tilde H = \C^n\backslash 0$
the $\Z$-covering of $H$, and let $\tilde M \subset \tilde H$
be the corresponding $\Z$-covering of $M$.
Then the closure of $\tilde M$ in 
$\C^n \supset \tilde H=\C^n\backslash 0$
is a complex analytic subvariety of $\C^n$, equipped
with the structure of a closed cone variety.
Moreover, all closed cone varieties are obtained this way.

\proof \ref{_cone_to_Hopf_Theorem_}. \endproof

\hfill

Let $\gamma$ be an automorphism 
of an algebraic variety $X$. We say that
$\gamma$ is {\bf strongly algebraic} if
it belongs to an algebraic group ${\cal G}$
acting on $X$. This notion is non-trivial.
Indeed, it is not hard to produce 
automorphisms which are not strongly algebraic
(Subsection \ref{_Henon_automo_Subsection_}).

\hfill

Let $A\in \End(\C^n)$ be a 
linear contraction in a complex vector space.
We already cited
\cite[Lemma 6.1]{_OV_Algebraic_Cones_}, claiming that 
an $A$-finite holomorphic function on $\C^n$ is polynomial.
The GAGA theorem follows from this observation and
the embedding theorem (\ref{_cone_to_Hopf_Theorem_}).

\hfill

\theorem
Let $(X_0, \gamma)$ be an open cone variety.
Then $X_0$ admits a unique algebraic structure 
such that the map $\gamma$ is strongly algebraic.
Moreover, in this algebraic structure, $\gamma$-finite functions
are regular (that is, algebraic), 
and, conversely, all regular sections of $H^0(X, \calo_X)$
are $\gamma$-finite.

\proof \ref{_algebraic_str_on_cone_uniqueness_Theorem_}.
\endproof

\hfill

In Section \ref{_affine_cone_Section_},
we analyze this algebraic structure on open cone varieties,
showing that they admit a proper $\C^*$-action
with projective quotients. Later on, this
statement will be transformed in an isomorphism
between an open cone variety and the space of non-zero
vectors in an ample line bundle over a projective
orbivariety.

\hfill

\theorem
Let $(X, \gamma)$ be a closed cone variety, and
$X_0$ the corresponding open cone variety.
Then $X_0$ is equipped with 
a locally free $\C^*$-action, commuting with the
action of $\gamma$, such that $Z:=X_0/\C^*$ is a projective
variety.

\proof  \ref{_cone_str_on_open_cone_variety_Theorem_}. \endproof

\hfill

We prove the GAGA theorem for equivariant coherent sheaves
on closed cone varieties in
Section \ref{_GAGA_coherent_Section_}.

\hfill

\theorem
Let $(X, \gamma)$ be a closed cone variety.
Using \ref{_algebraic_str_on_cone_uniqueness_Theorem_}, 
we consider $X$ as an affine variety.
Consider a $\langle \gamma\rangle$-equivariant analytic reflexive coherent 
sheaf ${\cal F}$ on $X$. Then ${\cal F}$
admits an algebraic structure, which is 
uniquely determined by the 
$\gamma$-equivariant structure.

\proof \ref{_coherent_GAGA_Theorem_}. \endproof

\hfill

In Section \ref{_orbivarieties_equiv_sheaves_Section_}, 
we introduce the orbivarieties. We prove that
for any cone variety $(X, \gamma)$,
there exists a sequence of contractions $\gamma_i$
converging to $\gamma$ such that 
each $\gamma_i$ is contained in 
a subgroup $\C^*\subset \Aut_\gamma(X)$
(\ref{_cone_str_on_open_cone_variety_Theorem_}).%
\footnote{$\Aut_\gamma(X)$ denotes the group of 
automorphisms commuting with $\gamma$.}
This implies that any closed algebraic
cone is isomorphic to an a affine cone over
a projective orbivariety.

\hfill

\theorem
Let $(X, x, \gamma_0)$ be a closed cone variety such that
the Zariski closure of $\gamma_0$ is isomorphic
to $\C^*$, and $X_0:= X\backslash x$ the
corresponding open cone variety. 
In this situation we denote the 
Zariski closure of $\gamma_0$ by $\C^*_Z$. Then 
the quotient $Z:=X_0/\C^*_Z$ is equipped
with a natural structure of a projective
orbivariety and an ample orbifold line bundle $L$ on $Z$, such that
$X_0$ is biholomorphic to $\Tot^0(L)$.

\proof \ref{_C^*_action_orbispace_Theorem_}. \endproof

\hfill

In the last section of this paper we analyse the
category of equivariant reflexive coherent sheaves on an
open cone variety. We prove that any such sheaf
admits a filtration with associated graded terms
${\cal F}_i$ which are isomorphic, as $\gamma$-equivariant
sheaves, to $\pi^*(\pi_*({\cal F}_i\otimes L_i^{-1}))\otimes L_i$,
where $L_i$ is an equivariant line bundle on $X_0$.
This is used to prove the filtrability theorem:

\hfill

\theorem
Let $(X_0, \gamma)$ be a cone variety
and ${\cal F}$ a $\gamma$-equivariant
reflexive sheaf on $X_0$, or, equivalently,
a coherent sheaf on $X_0/\langle \gamma\rangle$.
Assume that $\dim_\C X_0 \geq 3$. Then 
${\cal F}$ admits a $\gamma$-equivariant 
filtration with all associated graded subquotients
$\gamma$-equivariant coherent sheaves 
of rank $\leq 1$.

\proof \ref{_filtrable_Corollary_}. \endproof


\section{Cone varieties}
\label{_cone_Section_}



In this section we define the cone varieties, open and closed,
and explain how they are related.

\hfill

\definition\label{_contraction_cone_Definition_}
Let $\gamma:\; X \to X$ be a 
holomorphic automorphism of a connected complex analytic variety
preserving a point $x\in X$, called {\bf the origin}, or {\bf the apex}.
We say that $\gamma$ is {\bf a holomorphic contraction}
of $X$ to $x$, or {\bf a holomorphic contraction with origin in $x$}
if for any open $U\ni x$ and any compact $K\subset X$,
a sufficiently big power of $\gamma$ takes $K$ to $U$.
Consider a connected
complex analytic variety $X$, normal outside $x\in X$.
Assume that $X$ is equipped with an invertible holomorphic contraction
with origin in $x$. Then $X$ called {\bf a closed cone variety}, and
the normal variety 
$X \backslash x$ is called {\bf an open cone variety}.
This terminology is related to the notion of an affine
cone of a projective variety. Indeed, a closed cone variety
is always biholomorphic to an affine cone over a projective
(orbi)-variety (Subsection \ref{_all_are_cones_Subsection_}).

\hfill

\remark 
By \ref{_closed_cone_Stein_Theorem_}, a closed cone variety is always Stein.

\hfill

\remark
Clearly, an orbit of a holomorphic contraction acting on the closed cone
variety has a unique limit point in the origin. This means
that $\Z=\langle\gamma\rangle$ acts on an open
cone variety properly discontinuously.

\hfill

\example
Let $\gamma$ be an invertible holomorphic contraction of
$\C^n$ with apex 0. Then  $\frac{\C^n\setminus 0}{\langle
\gamma\rangle}$ is called {\bf a Hopf manifold}. If
$\gamma\in \GL(n,\C)$ is a linear contraction at 0, then
$\frac{\C^n\setminus 0}{\langle \gamma\rangle}$ is called {\bf a
  linear Hopf manifold}.

\hfill

\theorem\label{_cone_Hopf_embedding_Theorem_}
Let $H$ be a Hopf manifold, $Z\subset H$
a closed subvariety, $\tilde H = \C^n \backslash 0$ 
the universal cover of $H$, and $\tilde Z$ the preimage
of $Z$ in $\tilde H$. Then $\tilde Z$ is an open cone variety,
and any open cone variety can be obtained this way.

\proof \ref{_cone_embedded_Theorem_}.
\endproof

\hfill

\theorem
Every open cone variety 
is isomorphic to the space of non-zero vectors in an ample
bundle on a projective variety.

\proof \ref{_cone_str_on_open_cone_variety_Theorem_}. \endproof

%
%

\hfill

The passage from an open cone variety to the closed cone
variety is not without pitfalls. Recall that a 
complex analytic variety $X$ is called {\bf normal}
if any locally bounded meromorphic function on $X$ 
is holomorphic.
Unless we assume the normality of $X$, 
there are several non-isomorphic closed cones $(X,x,\gamma)$
with isomorphic open cone variety $(X_0, \gamma)$.
These non-isomorphic closed cone varieties have
geometric significance, because they can be obtained
as closures of $X_0$ holomorphically and $\gamma$-equivariantly
embedded to $\C^n$: different embeddings can produce
different closed cone varieties.
We discuss the details of this arrangements in 
\cite{_OV_Algebraic_Cones_}. For the present
purposes, it is sufficient to assume that $X$ is normal.

\hfill

\definition\label{_Stein_completion_Definition_}
Let $X$ be a Stein variety, normal
outside of $x\in X$, and $X_0:= X \backslash x$.
{\bf A Stein completion} of $X_0$ is
a Stein variety homeomorphic to $X$,
such that the natural embedding $X_0 \to X$
is a holomorphic embedding.

\hfill

\remark As explained in \cite{_OV_Algebraic_Cones_},
the Stein completion is not unique; however, its
normalization is the unique normal Stein completion
(\cite[Corollary 3.2]{_Andreotti_Siu:embeddings_}).

\hfill

\theorem\label{_closed_cone_from_open_Theorem_}
Let $(X_0, \gamma_0)$ be an open cone variety.
Then there exists a unique normal closed cone variety
$(X, \gamma)$ with origin in $x$ such that $X_0= X \backslash x$.

\proof By definition of an open cone, the closed cone exists;
it is only the uniqueness of $X$ which is an issue. By 
\ref{_closed_cone_Stein_Theorem_} below, $X$ is Stein. 
Then, by definition, $X$ is isomorphic to the 
(unique) normal Stein completion of $X_0$.
\endproof


\section{Algebraization for complete local rings}
\label{_algebraiz_for_rings_Section_}


In algebraic geometry, the ``algebraization''
alludes to passage from formal rings or formal schemes to 
finitely generated rings or schemes, originally
performed by Michael Artin \cite{_Artin:Approximation_}.
In this section we produce an algebraization result
for a complete local ring equipped with a contraction
(\ref{_Contraction_for_local_rings_Definition_}), 
showing that such a ring is always
a completion of a finitely generated ring. 

\subsection{Contractions of complete local ring and finite vectors}

In this subsection, we state the main result of this
section, the algebraization theorem.
The ``rings'' in this paper are 
commutative $\C$-algebras with unity. We express this by saying
``$R$ is a commutative ring over $\C$''. All automorphisms
of rings over $\C$ are assumed to be $\C$-linear.

\hfill

\definition\label{_F_finite_Definition_}
Let $F\in \End(V)$ be an endomorphism of a vector space.
A vector $v\in V$ is called {\bf $F$-finite}
if the space generated by $v, F(v), F(F(v)),  ...$
is finite-dimensional. 

\hfill

\definition
Let $R$ be a local ring, and $\goth m$
its maximal ideal. We say that $R$ is {\bf complete}
if it is complete in the $\goth m$-adic topology,
and {\bf finite type} if $\frac{\goth  m}{\goth m^2}$
is finite-dimensional.

\hfill

\definition
Let $V=\C^n$. A {\bf linear contraction} is an
endomorphism $A\in \End(V)$ such that all its
eigenvalues $\alpha_i$ satisfy $|\alpha_i| < 1$.

\hfill

\definition\label{_Contraction_for_local_rings_Definition_}
Let $R$ be a local ring, $\goth m$ its maximal ideal, 
and $\tau\in \Aut(R)$. We say that $\tau$ is a {\bf contraction}
if $\tau$ acts on $\frac{\goth m}{\goth m^2}$ as a linear contraction.

\hfill

The main result of this section is the following theorem.

\hfill

\theorem\label{_algebraization_Theorem_}
Let $R$ be a complete local ring of finite type, 
and $\tau\in \Aut(R)$ a contraction. 
Then the ring $R^\tau$ of $\tau$-finite elements in $R$ is dense
and finitely generated.

\hfill

We prove \ref{_algebraization_Theorem_}
in Subsection \ref{_Algebraiza_Subsection_}.

\subsection{Poincar\'e-Dulac theorem for complete local rings}

We start with the following formal power series version
of the Poincar\'e-Dulac theorem.

\hfill

\theorem\label{_Poinc_Dulac_formal_Theorem_}
Let $\tau\in \Aut(\C[[t_1, ..., t_n]]$ be an
automorphism of the ring of power series, and 
$\goth m\subset \C[[t_1, ..., t_n]]$ the maximal ideal.
Assume that $\tau$ acts on $\frac{\goth m}{\goth m^2}$
as a linear contraction. Then there exists
an automorphism $P\in \Aut(\C[[t_1, ..., t_n]])$
such that map $T:=P\tau P^{-1}$ is polynomial.
Moreover, for an appropriate choice of
$t_1, ..., t_n$, $T(t_k) = \lambda_k t_k + Q_k(t_1, ..., t_{k-1})$,
where $Q$ is a polynomial map.

\hfill

\proof
A complex analytic version of this theorem
 is proven in \cite[Theorem 4.6]{_Berteloot:PD_}.
Its proof is based on \cite[Lemma 4.4]{_Berteloot:PD_},
which works for power series, and the proof
of \cite[Theorem 4.6]{_Berteloot:PD_}
 carries over
to the power series ring (see also \cite{_Boureau:Hopf_}).
The last statement follows from the normal form
of a contraction in the presence of a resonance
(\cite{_Arnold:ODE+_}).
\endproof

\hfill

\proposition\label{_finite_PD_Proposition_}
Let $\tau\in \Aut(\C[[t_1, ..., t_n]])$ 
be an automorphism of the power series ring,
written in the Poincar\'e-Dulac normal form
as above, $\tau(t_k) = \lambda_k t_k + Q_k(t_1, ..., t_{k-1})$.
Then a power series $f\in \C[[t_1, ..., t_n]]$
is $\tau$-finite if and only if it is polynomial.

\hfill

\proof
Since $\tau(t_k)=\lambda_k t_k$ modulo $\C[t_1,..., t_{k-1}]\subset \C[[t_1, ..., t_n]]$, 
the space $\langle f, \tau(f), \tau(\tau(f)),... \rangle$
is finite-dimensional when $f$ is a polynomial.
For the converse, we generalize the argument
used in \cite[Lemma 6.1]{_OV_Algebraic_Cones_}
for the case of a linear $\tau$.

Define a lexicographic order on monomials
$t_1^{a_1}t_2^{a_2}... t_n ^{a_n}$,
with $t_1^{a_1}t_2^{a_2}... t_n ^{a_n}\prec t_1^{b_1}t_2^{b_2}... t_n ^{b_n}$
if $a_n < b_n$, or if $a_n = b_n$ and $a_{n-1} < b_{n-1}$,
or if $a_n = b_n$ and $a_{n-1} = b_{n-1}$ and
$a_{n-2} < b_{n-2}$ etc.
Then $\tau$ acts linearly on monomials
up to lower order terms. This implies
that all eigenvalues of $\tau$ are of form
$\lambda_1^{a_1}\lambda_2^{a_2}... \lambda_n^{a_n}$.

Since $\tau$ is a contraction, all $\lambda_i$
satisfy $|\lambda_i| <1$. Therefore,
any sequence $\{\lambda_1^{a_1}\lambda_2^{a_2}... \lambda_n^{a_n}\}$
converges to 0 as $\sum a_i$ goes to infinity. 
We obtain that every given number can be
realized as an eigenvalue of $\tau$ on
homogeneous polynomials of degree $d$ 
for finitely many choices of $d$ only.
Therefore, any root vector of $\tau$
is a finite sum of homogeneous polynomials.
This implies that
the Taylor decomposition of a $\gamma$-finite function $f$ can   have 
only finitely many components, otherwise the eigenspace
decomposition of $f$ with respect to the
action of $\tau$ is infinite. 
\endproof

\subsection{Algebraization of complete rings}
\label{_Algebraiza_Subsection_}

In this subsection, we prove \ref{_algebraization_Theorem_}.
Poincar\'e-Dulac theorem
shows that any formal contraction of the power series
ring becomes polynomial after an appropriate change
of coordinates. We are going to extend this result
to complete local rings of finite type. We need the 
following elementary lemma, proven in 
\cite[Lemma 5.15]{_OV_Algebraic_Cones_}
for Banach spaces. 

\hfill

\definition\label{_compact_type_Jordan_Definition_}
Let $K:\; V \to V$ be a
morphism of topological vector spaces.
We say that $K$ {\bf has a compact type Jordan form}
if $V$ has a dense $K$-invariant subspace $V'$,
and the restriction $K\restrict{V'}$ admits a Jordan block decomposition,
with each Jordan block associated with an eigenvalue
$\mu_i$ finite-dimensional, and $\lim_i \mu_i =0$.

\hfill

This lemma repeats \cite[Lemma 5.15]{_OV_Algebraic_Cones_}.

\hfill

\lemma\label{_RS_surj_Lemma_}
Consider a commutative square of continuous 
operators of topological vector spaces
\[ \begin{CD}
W_1 @> K_1 >> W_1\\
@V \phi VV @VV \phi V \\
W_2 @> K_2 >> W_2
\end{CD}
\]
with $\phi$ a surjective map. Suppose that
$K_i$, $i=1,2$, admit compact type Jordan forms.
Restricting $\phi$ to the space $W_i^{K_i}$ of $K_i$-finite vectors, we
obtain a map $\phi_f:\; W_1^{K_1}\to  W_2^{K_2}$.
Then $\phi_f$ is surjective.

\hfill

\proof
Let $U_2\subset W_2$ be a finite-dimensional space
of $K_2$-finite vectors, and $U_1\subset W_1$ its preimage.
Since $\phi$ is surjective, the natural map
$U_1 \to U_2$ is surjective. Since $U_1^{K_1}$ is dense in $U_1$,
and $U_2$ is finite-dimensional, the restriction
of $\phi$ to $U_1^{K_1}$ is also surjective.
\endproof

\hfill

We finish by proving the main result of this section
(\ref{_algebraization_Theorem_}):

\hfill

\theorem\label{_Poincar_Dulac_arb_ring_Theorem_}
Let $R$ be a complete local ring of finite type, 
and $\tau\in \Aut(R)$ a contraction. 
Then the ring $R^\tau$ of $\tau$-finite elements in $R$ is dense
and finitely generated.

\hfill

\pstep
Let $t_1,..., t_n\in {\goth m}$ be a set of elements
whose classes in  $\frac{\goth m}{\goth m^2}$ define a basis.
We are going to show that the natural map 
$j:\; \C[[z_1, ..., z_n]]\to R$ taking $z_i$ to $t_i$ is surjective.
Indeed, by \cite[Lemma 10.23]{_Atiyah_MacDonald_}, 
a morphism of complete local rings is surjective if
and only if the corresponding morphism of
associated graded rings is surjective.
On the other hand, by Nakayama lemma, every finitely generated 
module $W$ over a local Noetherian ring $A$ is generated by
any collection of elements $w_1, ..., w_n$
which generate $W/{\goth m}_A W$, where ${\goth m}_A$
is the maximal ideal of $A$. Therefore, ${\goth m}$
is generated by $t_1,..., t_n$ as an $R$-module.
Similarly, ${\goth m}^k$ is generated by homogeneous
polynomials of degree $k$ in $t_1, ..., t_n$,
which implies that the associated graded map
of $j$ is also surjective.

\hfill

{\bf Step 2:} 
Any homomorphism from $\C[z_1, ..., z_n]$ to
$\C[[z_1, ..., z_n]]$ taking $z_1, ..., z_n$
to elements generating the maximal ideal of $\C[[z_1, ..., z_n]]$
can be extended to a homomorphism of
$\C[[z_1, ..., z_n]]$, because $\C[z_1, ..., z_n]$
is dense in $\C[[z_1, ..., z_n]]$.
This homomorphism is surjective because
its associated graded map is surjective
(Step 1).

\hfill

{\bf Step 3:}
We lift the map $\tau\in\Aut(R)$
to an automorphism of $\C[[z_1, ..., z_n]]$
by taking each $z_i$ to any element $z_i'$
which satisfies $j(z_i')=\tau(z_i)$, and
extending this homomorphism to
$\C[[z_1, ..., z_n]]$ (Step 2). 
The power series version of Poincar\'e-Dulac
theorem (\ref{_Poinc_Dulac_formal_Theorem_}) implies that there exists
a set of generators ${\cal z}_1, ..., {\cal z}_n\in {\goth m}$ 
of the maximal ideal ${\goth m}$ of $\C[[z_1, ..., z_n]]$ such that
$\tau$ defines a automorphism of the polynomial ring
$\C[{\cal z}_1, ..., {\cal z}_n]\subset \C[z_1, ...,  z_n]$.
By \ref{_finite_PD_Proposition_}, 
an element $v\in \C[[{\cal z}_1, ..., {\cal z}_n]]$
is $\tau$-finite if and only if it is polynomial.
By \ref{_RS_surj_Lemma_}, the map $j$ is surjective
on $\tau$-finite vectors, hence $R^\tau$ is
a quotient of a polynomial ring.
Finally, $\C[{\cal z}_1, ..., {\cal z}_n]$
is dense in $\C[[{\cal z}_1, ..., {\cal z}_n]]$,
hence its image is dense in $R$.
\endproof


\section{Algebraization for modules over complete local rings}
\label{_modules_over_complete_Section_}

In this section, we extend the algebraization theorem,
which is central in Section \ref{_algebraiz_for_rings_Section_},
to equivariant modules over complete local rings.

\subsection{Equivariant modules over a ring}

\hfill

\definition
Let $V$ be a module over a ring $R$, and $\tau\in \Aut(R)$ 
an automorphism. We say that $V$ is {\bf $\langle\tau\rangle$-equivariant}
if there exist a vector space automorphism $\tilde \tau:\; V\to V$
such that for any $a\in R$, $v\in V$ one has 
$\tilde\tau(a v)=\tau(a)  \tilde \tau(v)$.

\hfill

\remark\label{_Serre_Swan_Remark_}
This definition is compatible with the notion
of $\langle \tau\rangle$-equivariant coherent sheaves (see
Subsection \ref{_equivariant_coherent_Subsection_} for the
definition).
By Serre-Swan theorem (\cite{_Serre_Faisceaux_,_Swan_}), 
the category of coherent sheaves
over a spectrum of a Noetherian ring $R$ is equivalent to the category
of finitely generated $R$-modules. Under this equivalence,
$\tau$-equivariant $R$-modules correspond to $\tau$-equivariant
coherent sheaves over $\Spec(R)$.

\hfill

\remark
Let $R$ be a complete local ring.
In the sequel, we will consider the topology on
finitely generated $R$-modules induced by the adic topology
on $R$. The base of this topology on an $R$-module $V$ is
$a+ {\goth m}^n V$, where $a\in V$ is an arbitrary element,
$n\in \Z^{\geq 0}$, and ${\goth m}$  the
maximal ideal of $R$.

\hfill

Recall that an $R$-module $V$ is called {\bf reflexive}
if the natural map $V\to V^{**}$ is an isomorphism, where
$V^*:=\HHom_R(V, R)$.
The main result of this section is the following 
algebraization theorem for reflexive equivariant $R$-modules.

\hfill

\theorem\label{_algebraization_R_modules_Theorem_}
Let $R$ be a complete local ring of finite type without zero divisors,
$\tau\in \Aut(R)$ a contraction, and $V$ a finitely generated
reflexive $R$-module. Let $\tilde\tau:\; V \to V$ be the  
structure map of a $\langle\tau\rangle$-equivariant module on $V$.
Denote by $R^\tau$, $V^{\tilde\tau}$ the corresponding
spaces of $\tau$-finite and $\tilde\tau$-finite elements.
Then $V^{\tilde\tau}$ is dense in $V$. Moreover, $V^{\tilde\tau}$
is finitely generated as an  $R^\tau$-module.

\hfill

We prove \ref{_algebraization_R_modules_Theorem_}
in Subsection \ref{_proof_finite_gen_modules_Subsection_}.

\subsection{The total space of a reflexive coherent sheaf}

Heuristically, we always think of $R$-modules as of coherent sheaves
over the corresponding affine variety
(\ref{_Serre_Swan_Remark_}). 
A torsion-free coherent sheaf 
is understood as a vector bundle with possible singularities, and 
vector bundles have total spaces. To prove 
\ref{_algebraization_R_modules_Theorem_},
we are going to define
a ``total space'' of an $R$-module, using this analogy, and
apply the algebraization theorem to this total space.
In this section, we always assume that
$R$ is a complete local ring of finite type
and without zero divizors, and $\goth m$ is its maximal ideal.

\hfill

\remark
Let $X$ be an affine variety.
Invoking again the Serre-Swan theorem,
we interpret vector bundles on $X$ as
projective modules over the ring $A:=H^0(\calo_X)$.
The fiberwise polynomial functions on the 
total space of a bundle $B$ 
are sections of $\bigoplus_i \Sym^i_A(H^0(X, B^*))$,
therefore $\Tot(B)= \Spec(\Sym^i_A(H^0(X, B^*))$.
This motivates the following definition.

\hfill

\definition
Denote the group of permutations of the set
$\{1, 2, ..., i\}$ by $\Sigma_i$.
Let $V$ be a finitely generated reflexive $R$-module.
Its {\bf total space} is the $I$-completion of the
ring $\bigoplus_{i=0}^\infty \Sym^i_R(V^*)$, where
$\Sym^i_R(V^*)$ is the space of $\Sigma_i$-invariant
vectors in $(V^*)^{\otimes i}$. Here $I$ denotes 
the maximal ideal in $\bigoplus_{i=0}^\infty \Sym^i_R(V^*)$ generated
by $\goth m$ and $\Sym^i_R(V^*)$, $i>0$.

\hfill

\remark
The total space of $V$, denoted as $\Tot(V)$,
is a complete local ring of finite type. 
Indeed, let ${\goth m}_V= I + \goth m$ be its
maximal ideal. Then 
\[ 
\frac{{\goth m}_V}{{\goth m}_V^2}= \frac{V^*}{\goth m V^*} 
\oplus \frac{\goth m}{{\goth m}^2},
\]
and this space is finite-dimensional because $V^*$ is
a finitely generated $R$-module.

\hfill

The contraction $\tau$ is extended to the total
space of $V$ using the equivariant structure operator $\tilde \tau$;
however, this new automorphism is not a contraction, because
it is not necessarily a contraction ``fiberwise''.
Fix a number $\lambda\in \C$, $|\lambda|<1$, and
define a new contraction automorphism $\tau_\lambda$ on 
the total space $\Tot(V)$ as follows:
$\tau_\lambda \restrict {\Sym^i_R(V^*)}= \lambda^i \Sym(\tilde \tau)$,
where $\Sym(\tilde \tau)$ denotes the extension of $\tilde\tau$
to $\Sym^i_R(V^*)$ acting on all tensor factors, as
$\Sym(\tilde\tau)(x_1\otimes... \otimes x_i)= 
\tilde\tau(x_1)\otimes... \otimes \tilde\tau(x_i).$

In the next subsection, we are going to apply
the algebraization theorem (\ref{_algebraization_Theorem_})
to the ring $(\bigoplus_i\Sym^i_R(V^*), \tau_\lambda)$ and interpret the
ring of $\tau_\lambda$-finite elements in $\bigoplus_i\Sym^i_R(V^*)$ 
as $\tilde\tau$-finite elements of $V$,
proving \ref{_algebraization_R_modules_Theorem_}.

\subsection{$\tau_\lambda$-finite elements in $\bigoplus_i\Sym^i_R(V^*)$}
\label{_proof_finite_gen_modules_Subsection_}

{\bf The proof of \ref{_algebraization_R_modules_Theorem_}:}
In the assumptions of \ref{_algebraization_R_modules_Theorem_},
consider $V=\Sym^1(V^{**})$ as a subspace of  the ring of
functions on $\Tot(V^*)$.
Then $\tau_\lambda$ acts on $V$ as $\lambda\tilde\tau$, hence
$\tau_\lambda$-finite vectors in $\Sym^1(V^{**})$ 
are the same as $\tau$-finite vectors in $V$.
By \ref{_algebraization_Theorem_}, the set 
of $\tau_\lambda$-finite 
vectors in $\bigoplus_i \Sym^i(V^{**})$ is dense in 
$\bigoplus_i \Sym^i(V^{**})$. Then 
$\Sym^1(V^{**})^{\tau_\lambda}$ is dense in $\Sym^1(V^{**})$.
This implies that the space $V^{\tilde\tau}$
is dense in $V$. This takes care of the first claim on
\ref{_algebraization_R_modules_Theorem_}.

To prove \ref{_algebraization_R_modules_Theorem_},
it remains to show that $V^{\tilde \tau}$ is
finitely generated as an $R^\tau$-module. However,
the ring $\bigoplus_i \Sym^i(V^{**})^{\tau_\lambda}$ is finitely 
generated by \ref{_algebraization_Theorem_}.
Consider the filtration on $\bigoplus_i \Sym^i(V^{**})$ with
$k$-th term generated by $\bigoplus_{i=k}^\infty \Sym^i_R(V^{**})$.
Its associated graded algebra is
$\bigoplus_{i=0}^\infty \Sym^i_R(V^{**})$. 
The associated graded algebra of $\bigoplus_{i=0}^\infty \Sym^i_R(V^{**})$ 
induced by this filtration is 
$\bigoplus_{i=0}^\infty \Sym^i_R(V^{**})^{\Sym(\tilde \tau)}$,
and it has to be finitely generated, because
the ring $\bigoplus_{i=0}^\infty \Sym^i_R(V^{**})^{\tau_\lambda}$ is finitely generated.

Clearly, any set of generators of
$\bigoplus_{i=0}^\infty \Sym^i_R(V^{**})^{\Sym(\tilde \tau)}$
intersected with $\Sym^1({V^{**}})^{\tau_\lambda}$
generates $\Sym^1({V^{**}})^{\tau_\lambda}=V^{\tilde\tau}$
over $R^\tau=\Sym^0({V^{**}})^{\tau_\lambda}$. Therefore, 
$V^{\tilde\tau}$ is finitely generated as an $R^\tau$-module.
We proved \ref{_algebraization_R_modules_Theorem_}.
\endproof


\section{Cone varieties and subvarieties in Hopf manifolds}
\label{_cone_in_Hopf_Section_}


We have already introduced the main characters of this paper,
the open and closed cone varieties
in Section \ref{_cone_Section_} (\ref{_contraction_cone_Definition_}).
Now we are going to prove the basic properties of 
the open and closed cone varieties and give examples.

\subsection{Closed cone varieties are Stein}

Recall that {\bf a closed cone variety}
is a complex analytic variety $X$ equipped with an invertible 
holomorphic contraction $\gamma$ (\ref{_contraction_cone_Definition_}).

\hfill

\theorem\label{_closed_cone_Stein_Theorem_}
A closed cone variety $(X, \gamma)$ is Stein.

\hfill

\noindent
\proof
By Behnke-Stein theorem (\cite[Proposition 12.39]{_OV_book_};
see also \cite{_Behnke_Stein_,_Stein:58_,_Andreotti_Narasimhan_}),
a union of an increasing sequence 
$U_1\Subset U_2 \Subset ... \Subset U_n \Subset ...$
of Stein varieties is Stein, if each $U_i$ is relatively
compact in $U_{i+1}$. Let $V\Subset X$ be a relatively
compact Stein neighbourhood of the origin. By definition
of a holomorphic contraction (\ref{_contraction_cone_Definition_})
for a sufficiently big $N\in \Z^{>0}$, we have
$\gamma^N(V) \Subset V$. Then 
$X= \bigcup U_i$, where $U_1:=V$ and
$U_i:= \gamma^{-Ni}(V)$. Clearly,
$U_1\Subset U_2 \Subset ... \Subset U_n \Subset ...$,
hence $X$ is Stein by Behnke-Stein theorem.
\endproof

\subsection{Subvarieties in Hopf manifolds}

%
%
%

For the sequel, we are going to use the celebrated
Remmert-Stein theorem (\cite{_Gunning_Rossi_},
\cite[Chapter 2, \S 8.2]{_Demailly:AG_}).

\hfill

\theorem\label{_Remmert_Stein_Theorem_} 
({\bf Remmert--Stein}) \\
Let $M$ be a complex analytic variety, $A\subset M$ a closed complex
subvariety, $\dim A\leq k-1$ and $M_0=M\backslash A$.
Consider a subvariety $Z\subset M_0$
with all components of dimension $\geq k$.
Then its closure in $M$ is complex analytic.
\endproof

\hfill 

The main result of the present section is the following theorem.

\hfill

\theorem\label{_cone_to_Hopf_Theorem_}
Let $H$ be a Hopf manifold, and $M \subset H$ a 
complex subvariety. Denote by $\tilde H = \C^n\backslash 0$
the $\Z$-covering of $H$, and let $\tilde M \subset \tilde H$
be the corresponding $\Z$-covering of $M$.
Then the closure of $\tilde M$ in 
$\C^n \supset \tilde H=\C^n\backslash 0$
is a complex analytic subvariety of $\C^n$, equipped
with the structure of a closed cone variety.
Moreover, all closed cone varieties are obtained this way.

\hfill

\pstep 
By the Remmert-Stein theorem (\ref{_Remmert_Stein_Theorem_}),
the closure $X$ of $\tilde M$ in $\C^n$ 
is complex analytic. The deck transform group
$\Z$ acts on $\C^n$ by holomorphic contractions,
hence $X$ is a closed cone variety.

\hfill

{\bf Step 2:} It remains to show that for any
closed cone variety $(X, \gamma)$ with apex in $x\in X$,
the quotient\footnote{The quotient 
$(X\backslash x)/\langle\gamma\rangle$ is well defined,
because $\gamma$ acts on $X\backslash x$ properly
discontinuously.}
$M:= (X\backslash x)/\langle\gamma\rangle$
admits a holomorphic embedding to a Hopf manifold.
This is \ref{_cone_embedded_Theorem_} below.
\endproof

\subsection{Montel theorem and $\gamma$-finite functions}

\definition \label{normal_family}
Let $M$ be a complex analytic variety,
and ${\cal F}\subset H^0(\calo_M)$ a 
family of holomorphic functions. We call ${\cal F}$ 
{\bf a normal family} if for each compact
$K\subset M$ there exists a constant $C_K>0$ such that
for each $f\in {\cal F}$, $\sup_K |f| \leq C_K$.

\hfill

\definition 
The {\bf $C^0$-topology} on the space
of functions on $M$ is the topology
of uniform convergence on compacts.

\hfill

\theorem \label{_Montel_Theorem_}
({\bf Montel})\\
Let ${\cal F}\subset H^0(\calo_M)$
be a normal family of functions,
and $\bar {\cal F}\subset H^0(\calo_M)$ its closure
in the $C^0$-topology. Then $\bar {\cal F}$ is compact
in the $C^0$-topology.

\proof 
\cite[Lemma 1.4]{_Wu:Montel_}. \endproof

\hfill

\theorem\label{_bounded_functions_Banach_Theorem_}
Let $Z$ be a complex analytic variety, $U\Subset Z$
a relatively compact open subset,
and $W:= H^0_b(\calo_U)$ the space of bounded
holomorphic functions, admitting an extension
to a neighbourhood of the closure $\bar U$.
We consider $W$ as a normed space equipped with the 
sup norm $\|f\|:=\sup_{z\in \bar U}|f(z)|.$
Then $W$ is Banach.

\hfill

\proof By \ref{_Montel_Theorem_}, any bounded subset in
$W$ is a normal family, hence every bounded sequence
of holomorphic functions in $W$ has a convergent subsequence.
Applying this to a Cauchy sequence, we obtain that $W$
is complete.
\endproof

\hfill

\corollary\label{_gamma_on_cone_compact_Corollary_}
Let $(X, \gamma)$ be a closed cone variety,
and $U\subset X$ a precompact open subset such that
$\gamma(U) \Subset U$. Then the map
$\gamma^*:\; H^0_b(U, \calo_U)\to H^0_b(U, \calo_U)$
is compact.

\hfill

\proof 
Since $\gamma(U) \Subset U$, for any bounded set 
$A\subset H^0(U, \calo_U)$, its pullback
$\gamma^*(A)\restrict U$ is a normal family by \ref{_Montel_Theorem_}.
Therefore, $\gamma^*$ maps a ball in the
Banach space $H^0_b(U, \calo_U)$ to a precompact
subset.
 \endproof

\hfill

We use the following corollary of the 
Riesz-Schauder theorem, see \cite{_OV_Algebraic_Cones_}
or \cite[Section 5.2]{friedman}.
Let $F:\; V \to V$ be a endomorphism of a vector space.
Recall that a vector $v\in V$ is called {\bf $F$-finite}
if the space $\langle v, F(v), F(F(v)), ... \rangle$
is finite-dimensional. An operator $F$ is called
{\bf quasi-nilpotent} if its spectral radius is 0,
that is, if $\lim_n \|F^n\|^{1/n}=0$. {\bf The quasi-nilpotent
part} of $(V, F)$ is the maximal $F$-invariant subspace
$W\subset V$ such that $F\restrict W$ is quasi-nilpotent.

\hfill

\theorem\label{_Riesz_Schauder_Theorem_}
Let  $F:\; V \to V$ be a compact operator on a Banach space
with finite-dimensional quasi-nilpotent part.
Then the space of $F$-finite vectors is dense in $V$.
\endproof

\hfill

The following criterion for $(V,F)$ having
no quasi-nilpotent part holds in all cases we consider.

\hfill

\claim\label{_bounds_by_C^n_Claim_}
Let $F:\; V \to V$  be a compact operator on a Banach space.
Assume that for any non-zero $v\in V$, 
there exists $C>0$ such that $\lim\sup_n C^n\| F^n(v)\|>0$.
Then $(V,F)$ has no quasi-nilpotent part.

\hfill

\proof
Let $W\subset V$ be the quasi-nilpotent part, 
and $w\in W$ a non-zero vector. Since $\lim_n \|F^n(w)\|^{1/n}=0$,
for any $A>0$ one has
\[
\lim \|F^n(w)\|^{1/n}=0 \Rightarrow \lim(A^n\|F^n(w)\|)^{1/n}=0
 \Rightarrow  \lim\sup_n A^n\|F^n(w)\| =0,
\]
contradicting $\lim\sup_n C^n\| F^n(v)\|>0$ when $A=C$.
\endproof

\hfill

\proposition\label{_Banach_no_quasinilpotent_Proposition_}
Let $(X, \gamma, x)$ be a closed cone variety with the apex $x$,
and $\gamma^*:\; H^0_b(U, \calo_U)\to H^0_b(U, \calo_U)$ 
the compact morphism of Banach spaces
defined in \ref{_gamma_on_cone_compact_Corollary_}.
Then $\gamma^*$ satisfies the assumptions of 
\ref{_bounds_by_C^n_Claim_}.

\hfill

\pstep
Let $\goth m$ be the maximal ideal of $x$.
Then {\bf the $k$-jet of a function} is its class
in $\calo_{U}/{\goth m}^{k+1}$. Since the ring of 
germs of holomorphic functions is Noetherian,
Krull lemma implies that $\bigcap_n {\goth m}^n=0$.
Therefore, any non-zero function has non-zero $k$-jet
for some $k>0$.

\hfill

{\bf Step 2:}
Consider the map associating to a function $f\in H^0_b(U, \calo_U)$
its $k$-jet $f_k$ in the origin $x$ of the cone variety 
$X$. Since the derivatives of a function can be expressed
using the Cauchy formula, this map is continuous.
We are going to prove that for some $C$ the sequence 
$C^n (\gamma^*)^nf$ does not
converge to 0. For this purpose, we 
show that $\lim\sup_n |C^n (\gamma^*)^nf_k|> 0$
for appropriate $C$ and $k$.

\hfill

{\bf Step 3:}
Let $|\mu|$ be the minimum of the absolute values
of the eigenvalues of $\gamma^*$
on the Zariski cotangent space ${\goth m}/{\goth m}^2$.
Clearly, the smallest absolute value of the eigenvalues of $\gamma^*$
acting on the space of $k$-jets is $|\mu|^k$.
Therefore, $\lim\sup C^n|(\gamma^*)^nf_k| >0$ for $C> |\mu|^{-k}$,
and the sequence $C^n (\gamma^*)^nf$ cannot converge to 0
in $H^0_b(U, \calo_U)$, when $k$ is the smallest integer
such that $f_k \neq 0$ (Step 1).
\endproof

\subsection{$\gamma^*$-finite functions extended to cone varieties}

\proposition\label{_finite_function_extends_Proposition_}
Let $(X, \gamma)$ be a closed cone variety, with 
the apex $x\in X$, and $f$ a germ of a holomorphic
function in $x$. Assume that $f$ is $\gamma^*$-finite.
Then $f$ can be holomorphically extended to the whole of $X$.

\hfill

\proof
Let $W=\langle f_1, ..., f_n\rangle$ 
be a $\gamma^*$-invariant finite-dimensional space
of germs containing $f=f_1$. Denote by $A\in \End(W)$
the restriction of $\gamma^*$ to $W$.
Replacing $\gamma$ by $\gamma^n$ if necessary,
we may assume that there exists a Stein neighbourhood
$U\ni x$ such that $\gamma(U) \Subset U$.
Replacing $U$ by $\gamma^k U$, $k\gg 0$,
we can also assume that all $f_i$ are defined on $U$.
By definition, $A^{-m}(\gamma^*)^m f_i= f_i$. However,
the function $(\gamma^*)^m f_i$ is defined on $\gamma^{-m}(U)$.
Since $X= \bigcup_m \gamma^{-m}(U)$, the function $f=f_1$
naturally extends to $X$.
\endproof

\hfill

\theorem\label{_finite_restiction_dense_Theorem_}
Let $(X, \gamma)$ be a closed cone variety,
and $U\subset X$ a precompact open subset such that
$\gamma(U) \Subset U$. Consider the Banach space
$H^0_b(U, \calo_U)$.
Then the restrictions of $\gamma^*$-finite functions from $X$ 
to $U$ are dense in $H^0_b(U, \calo_U)$.

\hfill

\proof
Applying 
\ref{_Banach_no_quasinilpotent_Proposition_}
and \ref{_Riesz_Schauder_Theorem_}, we obtain that 
$\gamma^*$-finite functions are dense in
$H^0_b(U, \calo_U)$. By \ref{_finite_function_extends_Proposition_},
any $\gamma^*$-finite function on $U$ can be
extended to all of $X$.
\endproof

\subsection{Embedding the cone varieties in Hopf manifolds}

The following theorem is known for smooth open cone varieties
(\cite{_OV_book_}). Its proof for general cone varieties is not
much different from the one given in \cite{_OV_book_},
and follows the same lines as the standard proof of
Kodaira embedding theorem.

\hfill

\theorem\label{_cone_embedded_Theorem_}
Let $(X, \gamma)$ be a closed cone variety, with 
the apex $x\in X$, and $X_0:= X \backslash x$ the
corresponding open cone variety. Then the quotient
$H_X:=X_0/\langle \gamma\rangle$ admits a holomorphic
embedding to a linear Hopf manifold.

\hfill

\pstep
Let $W \subset H^0(\calo_X)$ be a $\gamma^*$-invariant
subspace of $H^0(\calo_X)$. Substracting a constant,
we assume that all $f\in W$  vanish in the origin 
$x\in X$. Since $\gamma$ is a contraction, 
$(\gamma^*)^n(f)$ converges to a constant function 
taking the value $f(x)$ everywhere. Since $f(x)=0$,
the map $\gamma^*$ acts on $W$ as a linear contraction;
in particular, all 
eigenvalues of $\gamma^*$ on $W$ satisfy $|\alpha_i| <1$.
Therefore, $\frac{W\backslash 0}{\langle \gamma^* \rangle}$
is a Hopf manifold. In Step 3 we are going to show that we can choose
$W$ in such a way that the functions
in $W$ have no common zeros outside of the origin. Then 
the evaluation map $X\backslash x\stackrel e\to W\backslash 0$
is compatible with the $\gamma$-action and defines
a commutative diagram
\[
\begin{CD}
	X\backslash x @> e >> W\backslash 0 \\
	@V \gamma VV @VV \gamma^* V \\
	X\backslash x @> e >> W\backslash 0.
\end{CD}
\]
Therefore, it defines a holomorphic map 
$H_X\to  \frac{W\backslash 0}{\langle \gamma^* \rangle}$
to a Hopf manifold. It remains to check that $W$ can be
chosen without common zeros, so that 
this map is well defined. In Step 5, we are going to show that 
for an appropriate choice of $W$, the map $e$ is well defined and 
$X\backslash x\stackrel e\to W\backslash 0$ is
a holomorphic embedding.

\hfill

{\bf Step 2:}
By \ref{_finite_restiction_dense_Theorem_},
the space of $\gamma^*$-finite functions on $X$
is dense in the space of holomorphic functions on
any relatively compact subset $U\Subset X$.
In particular, for any two distinct 
points $a, b\in U\backslash x$,
there exists a $\gamma^*$-finite function
$f$ such that $f(a) \neq f(b)$,
and also $f(a) \neq 0$.
Similarly, for any 1-jet $\xi$ of a function
on $X$, there exists a $\gamma^*$-finite 
function with 1-jet $\xi$.

\hfill

{\bf Step 3:} 
Let $Z_W\subset X\backslash x$ be the set of
common zeros of $W$. Since $Z_W$ is $\gamma$-invariant,
$Z_W$ is the preimage of a subvariety $Y_W$ of a compact
variety $H_X:=\frac{X\backslash x}{\langle \gamma \rangle}$.
Using Step 2, for any $a\in Y_W$, we can always enlarge $W$ 
to obtain a space of $\gamma^*$-finite functions
$W_1\supset W$ such that $Y_W \supset Y_{W_1}$ and $a\notin Y_{W_1}$.
Repeating this process, we obtain a sequence of subvarieties of 
$H_X$,
$Y_W \supsetneq Y_{W_1} \supsetneq Y_{W_2}\supsetneq \cdots$.
Since $H_X$
is compact, and the Zariski topology is Noetherian,
this sequence stabilizes, giving $Y_{W_k}=\emptyset$.
We have shown that there exists a $\gamma^*$-invariant finite-dimensional
space $W$ of functions on $X$ having no common zeros except $x$.

\hfill

{\bf Step 4:} From Step 1 and Step 3, we obtain 
a holomorphic map $\Psi_W:\;H_X\to 
\frac{W\backslash 0}{\langle \gamma^* \rangle}$.
It remains to show that it is an embedding for
a sufficiently large $W$. 
Let $A_W\subset H_X \times H_X$ be the set of all
pairs of points $(a,b)$ such that $\Psi_W(a)=\Psi_W(b)$.
Using the same argument as in 
Step 2, we can always enlarge $W$ to $W_1\supset W$
in such a way that $A_{W_1} \subsetneq A_W$, unless $A_W=\emptyset$.
Passing to a sequence
$W\subsetneq W_1 \subsetneq W_2 \subsetneq ...$
as in Step 3, we get a sequence
$A_w\supsetneq A_{W_1} \supsetneq ...$
of complex analytic subsets of $H_X \times H_X$.
Since $H_X \times H_X$ is compact, this sequence
terminates at $W_k$ such that $A_{W_k}$
is empty and $\Psi_{W_k}$ is injective.

\hfill

{\bf Step 5:} 
In Step 4, we have shown that 
$\Psi_W:\; H_X \to \frac{W\backslash 0}{\langle \gamma^* \rangle}$
is injective. To prove 
\ref{_cone_embedded_Theorem_} it remains to
show that $\Psi_W$ is a closed embedding,
that is, every germ of a holomorphic function
on $H_X$ is extended to a germ of a holomorphic
function on $\frac{W\backslash 0}{\langle \gamma^* \rangle}$.
As shown in Step 2, for an appropriate choice of
$W$, the corresponding pullback map is surjective
on 1-jets. Then $\Psi_W$ is a closed embedding
by \cite[Proposition 45.9]{_Kaup_}.
\endproof

\subsection{A $\C^*$-equivariant embedding to a Hopf manifold}

Let $G$ be a connected complex algebraic group, obtained
as the complexification of a real Lie group $G_0\subset G$.
The Weyl theorem known as ``the Weyl unitarian trick'' 
(\cite{_Humphreys_}) states that every complex representation of $G_0$
extends to a complex representation of $G$, and this
correspondence defines an equivalence of categories.
This observation becomes obvious if we
notice that the Lie algebra action of $G_0$
is automatically extended to the action
of its complexification.

We are going to extend the
$\langle \gamma \rangle$-equivariant structure on the 
embedding $X \hookrightarrow \C^n$
of a cone variety to $\C^n$ to a $\C^*$-equivariant 
structure, provided that the $\gamma$-action
on this cone variety extends to a $\C^*$-action.

\hfill

\theorem\label{_C^*_on_cone_embedding_Theorem_}
Let $(X_0, \gamma)$ be an open cone variety. 
Suppose that the action of $\gamma\in \Aut(X_0)$
can be extended to a $\C^*$-action $\rho:\; \C^* \to \Aut(X_0)$ on $X_0$.
Consider the embedding 
$H_X\to  \frac{W\backslash 0}{\langle \gamma^* \rangle}$
constructed in  \ref{_cone_embedded_Theorem_}.
Then $W\subset H^0(\calo_X)^*$ is $\rho$-invariant.

\hfill

\pstep
When a group $G$ acts on a vector space $V$, we will
say that a vector or a subspace is ``$G$-finite''
if it is contained in a finite-dimensional
representation of $G$.\footnote{When $G$ is cyclic, $G=\langle g\rangle$, we shall also say ``$g$-finite'' instead of $G$-finite.}
Using the Fourier series  decomposition, we obtain
that the space of $S^1$-finite vectors is dense in 
any smooth Banach representation of $S^1$. 
Then Weyl's unitarian trick implies that 
the $\C^*$-finite vectors are dense
in any Banach representation of $\C^*$
where $\C^*$ acts complex analytically.

\hfill

{\bf Step 2:}
Let $(X, \gamma, x)$ be the closed
cone associated with $(X_0, \gamma)$ (\ref{_closed_cone_from_open_Theorem_}). 
As in the proof of 
\ref{_cone_embedded_Theorem_}, consider 
an open, relatively compact neighbourhood
$U$ of $x\in X$.
Let $H^0_b(\calo_U)$ denote the space
of holomorphic functions on $U$
extensible to a neighbourhood of the
closure $\bar U$. Since the $\C^*$-action
commutes with $\gamma$, it preserves
the space $H^0(\calo_X)^\gamma$ of $\gamma$-finite
functions on $X$.\footnote{When $\gamma$ is an automorphism of $X$,  it acts linearly on the functions on $X$ through $\gamma^*$, hence $\gamma$-finite is just a shortcut for $\gamma^*$-finite.} We endow this space with the norm defined as $\|f\|:= \sup_{z\in U} |f(z)|$.
By \ref{_finite_restiction_dense_Theorem_}, 
the completion of $H^0(\calo_X)^\gamma$
with respect to this norm is equal to 
$H^0_b(\calo_U)$. For any $t\in \C^*$,
by definition of a contraction, 
there exists a sufficiently big power
of $\gamma$ such that 
$\rho(t)\circ \gamma^n(U)\Subset U$;
this implies that $\rho(t)$ is 
continuous.\footnote{This is where we use that $\gamma$ 
commutes with the $\C^*$-action: we twist $\rho(t)$
with $\gamma^n$ and show that the result is continuous.} Also, clearly, the corresponding
operator $\rho(t)\in \End(H^0(\calo_X)^\gamma)$ 
complex-analytically depends on $t$.
We have proven that $\rho$
defines a complex analytic 
representation of $\C^*$ on
the Banach space $H^0_b(\calo_U)$.

\hfill

{\bf Step 3:} 
Using the argument in Step 1,
we show that the Banach space
$H^0_b(\calo_U)$ is the completion of a 
direct sum of $\C^*$-finite subspaces.
Clearly, $\C^*$-finiteness implies
$\gamma$-finiteness. The converse implication is
also true (Step 4), and this implies the statement of
\ref{_C^*_on_cone_embedding_Theorem_}.

\hfill

{\bf Step 4:}
Consider
a $\gamma$-finite representation $V\subset H^0(\calo_X)$, 
and let $\tilde V$ be the smallest $\rho$-invariant
subspace in $H^0(\calo_X)$ containing $V$.
Since $\rho(\C^*)$ commutes with $\gamma$,
the action of $\gamma$ on $\tilde V$ has the same eigenvalues 
as on $V$. Since $\gamma$ has compact-type Jordan form
(\ref{_compact_type_Jordan_Definition_}),
the space $\tilde V$ is finite-dimensional.
Now, $\rho(\C^*)\subset \GL(\tilde V)$ is obtained
as the Zariski closure of $\langle \gamma\rangle$
(it is the smallest-dimensional connected
complex algebraic group containing $\langle \gamma\rangle$).
Therefore, any subspace $V'\subset \tilde V$ 
which is $\langle \gamma\rangle$-invariant is also
$\rho(\C^*)$-invariant. We have shown that $V$ 
is $\rho(\C^*)$-invariant. 
\endproof

\hfill

\remark\label{_C^*_on_cone_embedding_Remark_}
For the space $W$ constructed in \ref{_C^*_on_cone_embedding_Theorem_}, 
the Hopf manifold $\frac{W\backslash 0}{\langle \gamma^* \rangle}$
is equipped with a $\C^*$-action which is compatible
with the $\C^*$-action $\rho:\; \C^*\to \Aut(X_0)$ on $X_0$.


\section{GAGA theorem for cone varieties}
\label{_GAGA_for_cones_Section_}


We start this section with a few remarks about algebraic groups,
define strong algebraic automorphisms, 
and explain how the notion of a strong algebraic automorphism
is used to construct a canonical algebraic structure on a cone variety.

\subsection{Algebraic groups and automorphisms of affine varieties}

\definition
{\bf An algebraic group} $G$ is a group object in the
category of affine varieties.

\hfill

\remark
In other words, $G$ is equipped with a
multiplication map $\mu:\; G \times G\to G$, which
is a morphism of affine varieties, and
an inverse $G \to G$ which is also a
morphism of affine varieties, and these
two maps satisfy the axioms for group 
operations.

\hfill

\remark
As shown in \cite{_Borel:alg_groups_},
an algebraic group can always be realised
as a subvariety in $\GL(n, \C)$, closed under the
group operations.

\hfill

\definition
Let $M$ be an algebraic variety,
and $G$ an algebraic group. A morphism
$a:\; G\times M \to M$ is called {\bf an action 
of $G$ on $M$} if the composition $G\times G\times M\to M$
is associative: 
$(\mu\times \Id)\circ a = (\Id\times a) \circ a$.

\hfill

\remark
Clearly, in these assumptions, the group 
$G$ acts on $H^0(M, \calo_M)$.
A regular function on $M$ is called {\bf $G$-finite}
(\ref{_C^*_on_cone_embedding_Theorem_}, Step 1)
if it belongs to a finite-dimensional $G$-representation.

\hfill

\claim\label{_regular_G_finite_Claim_}
Let $M$ be a variety equipped with an action
of an algebraic group $G$. Then all regular functions
on $M$ are $G$-finite.

\hfill

\pstep
Consider $G$ equipped with the left action of $G$ 
on itself. Any regular function on $G$ is $G$-finite,
because $G$ is a subgroup of $\GL(n, \C)$, and the
action of $\GL(n, \C)$ on itself preserves the degree
of a polynomial, hence any polynomial on $\GL(n, \C)$
belongs to a finite-dimensional subrepresentation.

\hfill

{\bf Step 2:}
Consider the multiplication map $G \times M\stackrel \Psi \to M$,
and let $f\in H^0(M, \calo_M)$ be a regular function.
Consider the action of $G$ on $G \times M$ by translations
along the first variable. As follows from Step 1, 
$\Psi^* f$ is $G$-finite with respect to this action.
By associativity of the multiplication, the map 
$G \times M\stackrel \Psi \to M$ is compatible with the
$G$-action, hence the $G$-finiteness of $\Psi^* f$
implies the $G$-finiteness of $f$.
\endproof

\subsection{Strongly algebraic automorphisms}
\label{_Henon_automo_Subsection_}

\remark
Let $X_0$ be an open cone variety. 
The closed cone variety is obtained as a Stein
completion (\ref{_Stein_completion_Definition_}) 
of $X_0$, which is not unique, 
as we have explained in \cite{_OV_Algebraic_Cones_}.
This non-uniqueness is a subtle matter which
lies out of the scope of the present paper.
This is why most statements in the present
section are made for open cone varieties.

\hfill

\definition
An inverse limit of algebraic groups is called {\bf a pro-al\-ge\-b\-raic
group}  (\cite{_Deligne:Tannakian_}). It is often infinite-dimensional.
Given a subroup $\Gamma \in \Aut(\calo_X)$, consider
the Zariski closure $\hat \Gamma_n$ of the subgroup generated by
its image in $\Aut(\calo_X/{\goth m}^n)$.
The {\bf proalgebraic Zariski closure} 
$\hat \Gamma$ of $\Gamma$ is 
$\lim\limits_{\leftarrow} \hat \Gamma_n$.

\hfill

\remark
Let $\goth G$ be the group of automorphisms
of the localization $\calo_{x, X}$ of $\calo_{X}$ in $x$,
where $X$ is the corresponding closed cone variety.
Clearly, $\goth G$ acts on $\frac{\calo_{X}}{\goth m^n}$,
where $\goth m$ is the maximal ideal of $x$.
This implies that $\goth G$ is embedded
to the inverse limit of algebraic groups,
$\goth G\subset \hat{\goth G}:=\lim\limits_{\leftarrow} \Aut\left(\frac{\calo_{X}}{\goth  m^n}\right)$.
Generally speaking, the proalgebraic Zariski closure of
a cyclic subgroup of $\Aut(\calo_{x, X})$ in $\hat{\goth G}$
is  proalgebraic, but sometimes the corresponding infinite sequence
of algebraic groups terminates. This is the subject
of interest in this subsection.

\hfill

\definition\label{_strongly_alg_Definition_}
Let $X$ be a quasiprojective variety,
and $\gamma\in \Aut(X)$ its automorphism.
We say that $\gamma$ is {\bf stronly algebraic}
if it belongs to an algebraic group acting on $X$.

\hfill

\remark\label{_strong_alg_degree_Remark_}
Let $\phi:\; \C^n \to \C^n$ be
an automorphism. Clearly, if $\phi$ is strongly
algebraic, all its iterates belong to a finite-dimensional
variety of automorphisms. Therefore, the degree of
$\phi^n$ is bounded for all $n$. Conversely, if
$\deg \phi^n$ is bounded, the linear span
of its iterates is finite-dimensional, hence
the Zariski closure of $\langle \phi\rangle$
is also finite dimensional. This characterizes
the strong algebraicity of $\phi\in \Aut(\C^n)$
in terms of the growth of the degree of its iterations.

\hfill

From now until the end of this subsection,
we are going to discuss this notion, giving 
a few examples and basic results.
We are not going to use these examples 
or results further in this paper, but the strong algebraicity is 
essential in \ref{_algebraic_str_on_cone_uniqueness_Theorem_}, and hence 
it is important to clarify its meaning.

\hfill

\definition
A {\bf generalized H\'enon map}
is an automorphism $\phi:\; \C[x, y]\to \C[x, y]$
defined by $(x, y) \to (\alpha y, P(y)-\beta x)$, 
where $\alpha, \beta\neq 0$ are constants, and $P$ is a polynomial with 
$\deg P >1$. 

\hfill

\remark
In \cite{_Friedland_Milnor:Henon_},
Friedland and Milnor prove that any polynomial
automorphism $\phi\in \Aut(\C^2)$ preserving zero is  conjugate to either an
affine map, to the map $(x, y) \mapsto (\alpha x+ P(y), \beta y)$,
called {\bf an elementary automorphism}, 
or to a composition of generalized H\'enon maps (such automorphisms
are called {\bf loxodromic}).
For the first two types, the degree of $\phi^n$
is constant for all $n\neq 0$, 
hence the Zariski closure of $\langle \phi\rangle$
is an algebraic group. Such automorphisms are strongly
algebraic. The degree of $\phi^n$ grows exponentially
when $\phi$ is loxodromic, hence a loxodromic automorphism
is not strongly algebraic. However, by a theorem
of Cantat and Dujardin
(\cite{_Cantat_Dujardin_}), any two polynomial
automorphisms of $\C^2$ which are holomorphically
conjugate are also polynomially conjugate.
Since any contraction is holomorphically
conjugate to an elementary automorphism
(Poincar\'e-Dulac theorem), the 
Cantat-Dujardin theorem implies that
the degree of a contraction in $\C^2$ does not
grow when we iterate, and hence it is strongly algebraic.
This proves the following result, which is due to
D. Korshunov (\cite{_Korshunov:MO_,_Korshunov:contractions_}).

\hfill

\proposition
Let $\phi\in \Aut(\C^2)$ be a polynomial contraction.
Then $\phi$ is strongly algebraic.
\endproof

\hfill

For $n\geq 3$, not every polynomial contraction of $\C^n$ is
strongly algebraic. The following example is also due
to Dmitri Korshunov, \cite{_Korshunov:MO_,_Korshunov:contractions_}.

\hfill

\example
Let $(x,y,z) \stackrel F \mapsto 
\left(\frac{y+xz}{2},\frac{x}{2},\frac{z}{2}\right)$. 
This map has linear degree growth. On the other hand,
denote $F^n(x,y,z)$ by $(x_n, y_n, z_n)$. Then
\[ z_n=\frac{z_0}{2^n}, \ y_n=\frac{x_{n-1}}{2},\ \text{ and}\ 
x_n=\frac{\frac{x_{n-2}}{2}+x_{n-1}\frac{z_0}{2^{n-1}}}{2}.
\]
Clearly, $x_n \leq \frac 12 \frac{x_{n-1}+x_{n-2}}{2}$ once
$z_0 \leq 2^{n-1}$. Therefore, for $n$ sufficiently big,
$\frac{x_{n-1}+x_{n-2}}{2}\leq \max(x_{n-1}, x_{n-2})$.
This implies that $x_n \leq \frac 12 \max(x_{n-1}, x_{n-2})$,
and this quantity converges to zero as $n\to \infty$.
Thus all of $\mathbb C^3$ is the basin of attraction of $(0,0,0)$.
This gives an example of a polynomial contraction which is
{\em not} strongly algebraic.

\subsection{The GAGA theorem for cone varieties}

\theorem\label{_algebraic_str_on_cone_uniqueness_Theorem_}
Let $(X_0, \gamma)$ be an open cone variety.
Then $X_0$ admits a unique algebraic structure 
such that the map $\gamma$ is strongly algebraic.
Moreover, in this algebraic structure, $\gamma$-finite functions
are regular (that is, algebraic), 
and, conversely, all regular sections of $H^0(X, \calo_X)$
are $\gamma$-finite.

\hfill

%

\proof
Clearly, the uniqueness of the algebraic structure
follows from the last assertion, which was proven
in \ref{_regular_G_finite_Claim_}.

By \ref{_cone_embedded_Theorem_}, $X_0$ admits an embedding
to $\C^n \backslash 0$, in such a way that $\gamma$
is obtained by the restriction of a linear contraction 
$A\in \GL(\C^n)$. Let $X$ be its closure in $\C^n$;
this closure is complex analytic by Remmert-Stein 
(\ref{_Remmert_Stein_Theorem_}). Using \ref{_RS_surj_Lemma_},
we obtain that every $\gamma$-finite function 
on $X$ can be extended to a $\gamma$-finite function
on $\C^n$. By \cite[Lemma 6.1]{_OV_Algebraic_Cones_},
a function on $\C^n$ is $\gamma$-finite 
if and only if it is polynomial. Therefore,
in this algebraic structure, $\gamma$-finite
functions coincide with the regular ones.
\endproof

\hfill

\remark\label{_normalization_closed_cone_algebraic_Remark_}
Let $(X, \gamma)$ be a closed cone variety, and 
$(X_1,\gamma)$ its normalization. 
Applying the same argument as used
in \ref{_algebraic_str_on_cone_uniqueness_Theorem_}
to  $(X_1, \gamma)$, we obtain that it is an affine variety,
and the algebraic structure is uniquely determined by
$(X, \gamma)$.

\hfill

Let $(X, \gamma,x)$ be a closed cone variety,  $\goth m$ the
maximal ideal of $x$, 
and $\hat \calo_X:= \lim\limits_{\leftarrow} \calo_X/{\goth m}^n$ 
the completion of $H^0(X, \calo_X)$ in $x$.
The variety $X$ is Stein (see \ref{_closed_cone_Stein_Theorem_} above)
and the natural map $\calo_X \to \hat \calo_X$
is injective by Krull's intersection theorem.
Therefore, the group $\Aut(X,x)$ of complex automorphisms
of $X$ preserving $x$ injects in $\Aut(\hat \calo_X)$.
The group $\Aut(\hat \calo_X)$ is clearly an inverse limit: 
$\Aut(\hat \calo_X)= \lim\limits_{\leftarrow} \Aut(\calo_X/{\goth m}^n)$.

\hfill

\remark
Consider an automorphism $\gamma\in \Aut(X)$,
and let $\Gamma:= \langle \gamma\rangle$ be the corresponding
cyclic group. Clearly, $\gamma$ is strongly algebraic if and only if
the proalgebraic Zariski closure $\hat \Gamma$ is finite-dimensional.
For a polynomial automorphism $F\in \Aut(\C^n)$, this is equivalent
to $\deg F^n$ being bounded (\ref{_strong_alg_degree_Remark_}).


\section{Cone varieties and affine cones}
\label{_affine_cone_Section_}


In this section, we prove that any
open cone variety is isomorphic to the space of non-zero vectors in 
the total space of an ample
line bundle over a projective variety.

In Subsection
\ref{_cone_example_Subsection_},
 we explore the cone structure on the space $\Tot^\circ(L)$ of 
non-zero vectors in  the total space of an ample
line bundle over a projective variety $Z$.
The corresponding ``affine cone'' is the closure
of $\Tot^\circ(L)$ in $H^0(Z, L)^*$. We identify
it with the spectrum of the ring 
$\bigoplus_i H^0(Z, L^{\otimes i})$.
When $Z$ is a smooth orbifold, these spaces
were studied in \cite{_OV_Algebraic_Cones_}; we called them
``an open algebraic cone'' and ``a closed algebraic cone''.

\subsection{Jordan-Chevalley decomposition}

For the sequel, we need a few basic results from
the algebraic group theory (\cite{_Humphreys_}).
We introduce the  Jordan-Chevalley decomposition
generalizing the Jordan normal form to an arbitrary
algebraic group.

\hfill

\definition 
An element of an algebraic group $G$ is called
{\bf semisimple} if its image is semisimple
for some faithful
algebraic representation of $G$, and is called 
{\bf unipotent} if its image is unipotent
(that is, the exponential of a nilpotent) 
for some exact algebraic representation of $G$.

\hfill

\remark 
For any algebraic representation of an algebraic group
$G$, the image of any semisimple element is a 
semisimple operator, and the image of any unipotent
element is a unipotent operator (\cite[\S 15.3]{_Humphreys_}).

\hfill

\theorem \label{_JC_Theorem_}
 ({Jordan-Chevalley  decomposition}, 
\cite[\S 15.3]{_Humphreys_}) \\
Let $G$ be an algebraic group, and $A\in G$.
Then there exists a unique decomposition $A= S U$
	of $A$ in a product of commuting elements $S$ and $U$,
	where $U$ is unipotent and $S$ semisimple.
\endproof

\subsection{Cone varieties and open cones}
\label{_cone_example_Subsection_}

In this subsection we construct an example of an open cone variety,
which we call ``an open quasi-affine cone''. It is an affine
cone without the origin. In Subsection
\ref{_all_are_cones_Subsection_} we are going to
show that all open cone varieties are obtained this way.

\hfill

\definition
Let $Z$ be a projective variety,
and $L$ an ample line bundle on $Z$.
Consider the total space
$\Tot^\circ(L)$ of non-zero vectors
in $L$; this is an algebraic
variety $\C^*$-fibered over $Z$.
Fix $k$ such that $L^k$ is very ample, and
$Z$ is embedded to $\C P^n= {\Bbb P}H^0(Z, L^k)^*$.
Let $C(Z, L^k)\subset H^0(Z, L^k)^*$ be the union 
of all lines in the corresponding vector space
passing through the image of $Z$.
Clearly, $\Tot^\circ(L)$ is isomorphic
to a $\Z/k$-cover of $\Tot^\circ(L^k)$.
We denote by $C(Z,L)$ or by $C(Z)$
the corresponding $\Z/k$ ramified cover of $C(Z, L^k)$;
it is a closed cone variety, equipped
with an algebraic structure induced
from $\Tot^\circ(L)$. Regular functions
on $C(Z,L)$ are the same as 
elements of $\bigoplus_i H^0(Z, L^i)$.
The variety $C(Z,L)$ is called 
{\bf the affine cone of $Z$}.

\hfill

\definition\label{_equiv_on_line_bundle_Definition_}
Let $L$ be an ample bundle on a projective
variety $Z$, and $\phi\in \Aut(Z)$ an automorphism.
{\bf A $\phi$-equivariant structure} on $L$
is a fiberwise automorphism of the total space
$\Tot(L)$ acting linearly on fibers
and coinciding with $\phi$ on $Z$.

\hfill

\remark
Let $\phi$ be an automorphis of $Z$ 
acting on $L$ (and, hence, on $C(Z,L)$) equivariantly.
The aim of the present section is to show that
the composition of $\phi$ and a homothety
$x \to \lambda x$ with sufficiently small
$\lambda\in \C^*$ is a contraction. This statement
is intuitively obvious, but it takes a bit of work
to prove it formally. We prove this statement
for $C(Z,L)\subset H^0(X,L)^*$ when $L$ is very ample.
Since any $C(Z,L)$ is obtained as a finite ramified covering
of such a variety, the argument carries over to the general case
when $L$ is ample.

\hfill

\lemma\label{_uniform_implies_contraction_Lemma_}
Let $X$ be an affine cone variety, 
and $F:\; X \to X$ an automorphism
preserving the origin $x\in X$. Assume that
for any regular function $f$ on $X$,
the sequence $\{(F^*)^i(f)\}$ converges to
a constant $f(x)$ uniformly on compacts.
Then $F$ is a contraction.

\hfill

\proof
Let $K(f,n) \subset X$ be a subset determined
as $K(f,n)= \{z\in X\ \ |\ \ |f(z)|\leq n\}$,
where $f$ is a regular function, vanishing at $x$, and $C$ a constant.
Taking an appropriate homogeneous polynomial function on the cone,
we can assume that $U(f, \epsilon) :=\{z\in X\ \ |\ \ |f(z)| <\epsilon\}$
(where $\epsilon$ varies) is a base of open 
neighbourhoods of $x$, the sets $K(f,n)$ are compact, and 
$\bigcup_n K(f,n)=X$.
To show that $F$ is a contraction,
it will suffice then to prove that for
any $n$, $\epsilon$  there exists 
a sufficiently big $N$ such that
$F^N(K(f,n)) \subset U(f, \epsilon)$.
This assertion clearly follows from $f(x)=0$ and
the uniform convergence $\lim_i (F^*)^i(f)=f(0)=0$.
\endproof

\hfill

\proposition\label{_examples_via_action_on_Zariski_tangent_Proposition_}
Let $L$ be a very ample bundle over a projective variety $Z$,
$\phi\in \Aut(Z)$ an automorphism, and 
$\tilde \phi:\; \Tot(L)\to \Tot(L)$ a $\phi$-equivariant
structure. By construction, $\tilde \phi$
can be extended to the corresponding 
algebraic cone variety, giving an automorphism  $\gamma:\; C(Z,L) \to C(Z,L)$.
Let $\goth m$ denote the maximal ideal of the origin, 
identified with the zero in $H^0(X, L)^*\supset C(Z,L)$. Then $\gamma$ is a contraction if
and only if its action on the
Zariski tangent space $\left(\frac{\goth m}{\goth  m^2}\right)^*$ 
is a linear contraction.

\hfill

\pstep
Assume that $\gamma$ is a contraction with apex $x$.
Then $\lim_i (\gamma^*)^i(f)=f(x)$.
Take $f\in {\goth m}$, and let $df$
denote its class in $\frac{\goth m}{\goth  m^2}$.
Since the sequence  $\{\gamma^i(f)\}$ converges to zero for all
$f\in \goth m$, the sequence $\{\gamma^i(df)\}$ converges
to zero, hence all eigenvalues of the $\gamma$-action
on $\frac{\goth m}{\goth  m^2}$ are less than 1 in absolute value.

\hfill

{\bf Step 2:} 
Conversely, assume that 
$\gamma^*$ acts on $\frac{\goth m}{\goth  m^2}$ as a linear contraction.
By construction, $\frac{\goth m}{\goth  m^2}$ 
is the space of all sections of $\bigoplus_k H^0(X, L^k)$,
up to products of sections non-zero degree:
\[
\frac{\goth m}{\goth  m^2}= 
\frac{\bigoplus_{k>0} H^0(X, L^k)}
{\bigoplus_{k>0} H^0(X, L^k)\cdot \bigoplus_{k>0} H^0(X, L^k)}.
\]
Therefore, $\gamma^*$ acts on $\frac{\goth m}{\goth  m^2}$ 
with eigenvalues $|\alpha_j|<1$ if and only if
it acts on each space $H^0(X, L^k)$ with
eigenvalues $|\alpha_j|<1$.
This gives $\lim_i (\gamma^*)^i(f)=f(x)$
for any $f\in \bigoplus_{k=0}^n H^0(X, L^k)$,
hence $\lim_i (\gamma^*)^i(f)=f(x)$ for any
regular function on $C(Z,L)$.
By \ref{_uniform_implies_contraction_Lemma_},
$\gamma$ is a contraction.
\endproof

\hfill

\remark
Consider an affine cone $C(Z,L)$,
and let $\tilde \phi$ be an automorphism of $Z$
extended to $C(Z,L)$ using an equivariant structure.
Taking the composition with a homothety
action $v\mapsto \lambda v$, for $|\lambda|$
sufficiently small, we obtain a contraction
of $C(Z,L)$ (\ref{_examples_via_action_on_Zariski_tangent_Proposition_}).

\subsection{All closed cone varieties are affine cones over projective varieties}
\label{_all_are_cones_Subsection_}

\theorem\label{_cone_str_on_open_cone_variety_Theorem_}
Let $(X, \gamma)$ be a closed cone variety, and
$X_0$ the corresponding open cone variety.
Then $X_0$ is equipped with 
a locally free $\C^*$-action, commuting with the
action of $\gamma$, such that $Z:=X_0/\C^*$ is a projective
variety. 

\hfill

\pstep
The proof repeats {\em mutandis mutatis} the proof
of \cite[Theorem 7.9]{_OV_Algebraic_Cones_}.
Using \ref{_cone_embedded_Theorem_}, we can 
embed $H_X:=X_0/\langle \gamma\rangle$ to a linear Hopf manifold.
This gives a holomorphic map
$X_0 \to \C^n\backslash 0$ taking
$\gamma$ to a linear contraction $A\in \End(\C^n)$.
Let $X'$ be the closure of $X_0$
in $\C^n$; it is complex analytic by Remmert-Stein theorem
(\ref{_Remmert_Stein_Theorem_}). By construction, $(X', \gamma)$ is
a closed cone variety, possibly different from $X$
(but with the same normalization, as explained in
\cite[Section 8]{_OV_Algebraic_Cones_}).

By \cite[Lemma 6.1]{_OV_Algebraic_Cones_}), the 
$\gamma$-finite functions on $\C^n\backslash 0$
are polynomials. Since 
the image of $X_0$ is $\gamma$-invariant, it
is algebraic by \cite[Theorem 6.3]{_OV_Algebraic_Cones_}
(see also \ref{_algebraic_str_on_cone_uniqueness_Theorem_}).

\hfill

{\bf Step 2:}
To prove \ref{_cone_str_on_open_cone_variety_Theorem_}
we are going to construct a holomorphic 
$\C^*$-action on $\C^n$, locally free
outside of the origin $0$, and preserving
$X'\subset \C^n$. Then the quotient
$W\CP^{n-1}= \frac{\C^n\backslash 0}{\C^*}$ is a weighted
projective space (\cite[Example 9.13]{_OV_book_}),
which is well known to be a projective orbifold, hence
$Z:= X_0/\C^*\subset W\CP^{n-1}$ is a projective
orbivariety (\ref{_orbivariety_Definition_}). On the other hand, 
$X_0$ is the total space of the corresponding
$\C^*$-bundle, or, what is the same,
the total space of non-zero vectors in a
holomorphic orbifold line bundle.

\hfill

{\bf Step 3:}
Let ${\cal G}_A$ be the Zariski closure of
$\langle A \rangle$ in $\GL(\C^n)$. This is a commutative
algebraic group, acting on the variety $X'\subset \C^n$.
To prove \ref{_cone_str_on_open_cone_variety_Theorem_}, we shall find
a subgroup $\C^*\subset {\cal G}_A$ containing a contraction and
apply Step 2. 

Let $A=SU$ be the Jordan-Chevalley decomposition 
(\ref{_JC_Theorem_}) for
$A$, with $S, U\in {\cal G}_A$. The algebraic group
${\cal G}_A$ preserves $X_0$, hence
the endomorphisms $S$ and $U\in \End(\C^n)$ also act on $X_0$.
Since the eigenvalues of $S$ are the same
as the eigenvalues of $A$, the map $S$ is a linear contraction.
Let $G_S\subset {\cal G}_A$, $G_S= e^{\C \log S}$ be the one-parametric subgroup
containing $S$. We prove that $G_S$ can be approximated by
subgroups of ${\cal G}_A$ isomorphic to $\C^*$. Then some of these subgroups
also contain a contraction, and we are done.

\hfill

{\bf Step 4:}
Consider the map taking any $A_1\in {\cal G}_A$ to its unipotent component
$U_1$. Since ${\cal G}_A$ is commutative, this map is a group
homomorphism. Therefore, its kernel ${\cal G}_s$ (that is,
the set of all semisimple elements in ${\cal G}_A$) is an
algebraic subgroup of ${\cal G}_A$. A connected commutative
algebraic subgroup of $\GL(\C^n)$ consisting of semisimple elements
is always isomorphic to $(\C^*)^k$ 
(\cite[Proposition 1.5]{_Borel_Tits:Groupes_Reductifs_}).
 The one-parametric subgroups $\C^*\subset (\C^*)^k$ 
are dense in $(\C^*)^k$ because one-parametric
complex subgroups $\C^*\subset (\C^*)^k$ can be obtained as 
complexification of subgroups $S^1\subset \mathrm{U}(1)^k \subset (\C^*)^k$,
and those are dense in $\mathrm{U}(1)^k$. Therefore,
the contraction $S\in {\cal G}_s=(\C^*)^k$ can be approximated
by an element of $\C^*$ acting on $X_0$.
\endproof 

\hfill

\remark
As shown in 
\ref{_C^*_action_orbispace_Theorem_} 
below, the variety $Z= X_0/\C^*$ defined in
\ref{_cone_str_on_open_cone_variety_Theorem_} is equipped
with a natural structure of an 
orbivariety\footnote{The orbivarieties are discussed in detail in
Section \ref{_orbivarieties_equiv_sheaves_Section_}, 
see, in particular, \ref{_orbivariety_Definition_}.}
and an orbifold line bundle $L$, such that
$X_0$ is biholomorphic to $\Tot^0(L)$.
Later we will show that all
cone varieties are obtained this way (\ref{_C^*_action_orbispace_Theorem_}).
We postpone this argument because to state it properly,
we need to define the orbifolds, orbispaces, and
line bundles over orbispaces.


\section{Fr\'echet structure on coherent sheaves}


In this section, we recall the analytic side of 
the theory of coherent sheaves on complex analytic varieties.
This makes a digression from
our narrative. To continue the arguments which 
lead to \ref{_cone_str_on_open_cone_variety_Theorem_},
we need to define the coherent sheaves on orbispaces. Before we
do that, we prove a GAGA-type theorem for equivariant
sheaves on cone varieties. This requires us to introduce
a number of analytic concepts, used in the
proof of the GAGA theorem.

\subsection{Montel theorem for coherent sheaves}

Montel's theorem can be extended to the space
of sections of a coherent sheaf, as follows. 
We define the Fr\'echet structure
on coherent sheaves as in
\cite[Chapter VIII.A]{_Gunning_Rossi_}, who refer
to Cartan and Serre for their definition.
For Montel sheaves and their application, see 
\cite{_Grothendieck_Montel_} (Grothendieck used
the term ``compact sheaves'' for the notion that, following
\cite{_Gunning_Rossi_}, we call ``Montel sheaves'').

\hfill

\definition
Let $U\subset V$ be open subsets of a Hausdorff 
topological space. Recall that $U$ is called
{\bf relatively compact}, or {\bf precompact} if its closure in $V$ is compact.
We denote the relative compactness by $U\Subset V$.

\hfill

\definition
Let $V$ be a topological vector space.
A subset $K\subset V$ is {\bf bounded}
if for any neighbourhood $U$ of zero,
there exists a number $\lambda \neq 0$
such that $\lambda K \subset U$.

\hfill

\definition
A {\bf Fr\'echet sheaf} is a sheaf of Fr\'echet spaces,
with continuous restriction maps. A Fr\'echet sheaf
${\cal F}$ is called {\bf Montel} if for any open sets
$U\Subset V$, the restriction map of topological
vector spaces ${\cal F}(V) \to {\cal F}(U)$ is compact
(that is, takes bounded sets to precompact sets).

\hfill

\remark\label{_Frchet_on_free_Remark_}
The standard Fr\'echet topology on the sheaf
of holomorphic functions $\calo_X$ is defined as above, by
uniform convergence on compact sets.
Then, \ref{_Montel_Theorem_} can be reformulated
by saying that the sheaf of holomorphic functions
on a complex analytic variety is Montel. 

\hfill

\remark  
Consider a coherent subsheaf  ${\cal F}\subset \calo_X^n$.
Clearly, $H^0(U, {\cal F})$ is closed in $H^0(U,\calo_X^n)$
for any open set $U\subset X$. This defines the structure
of a Fr\'echet sheaf on ${\cal F}$. Also, 
for any quotient sheaf  ${\cal H}$ of ${\cal F}$, the corresponding
projection maps $H^0(U, {\cal F})\to H^0(U, {\cal H})$
have closed kernels, hence ${\cal H}$ is also a Fr\'echet sheaf.

\hfill

\proposition\label{_coherent_Frechet_Proposition_}
Let  $X$ be a complex analytic variety. 
Using \ref{_Frchet_on_free_Remark_}, we
equip a free sheaf $\calo_U^n$ on an open subset $U\subset X$
with a structure of Fr\'echet sheaf.
Consider a coherent sheaf ${\cal F}$ on $X$.
Then there 
exists a unique structure of a Fr\'echet sheaf on 
${\cal F}$ such that for any open subset $U\subset X$
and any surjective coherent sheaf morphism 
$\phi:\; \calo_U^n \to {\cal F}\restrict U$,
the map $\phi$ defines a morphism of Fr\'echet sheaves
(that is, $\phi$ is continuous on the spaces of sections).
Moreover, the topology on $H^0(U, {\cal F})$
is defined by a collection of pseudonorms
$\|\cdot \|_K$ where $K$ runs through all compact 
subsets $K\subset U$, and satisfies the following 
conditions.
\begin{description}
\item[(i)] For any $f\in H^0(U, {\cal F})$, 
the value $\|f\|$
is uniquely determined by the germ of $f$ in $K$,  and
\item[(ii)]
 $\|f\|_K\geq \|f \|_{K'}$ 
for any $K' \subset K$.
\end{description}

\noindent
\proof \cite[Theorem VIII.A.8]{_Gunning_Rossi_}.
\endproof

\hfill

\remark
Using \ref{_coherent_Frechet_Proposition_}, we can
always consider a coherent sheaf as a Fr\'echet sheaf;
the Fr\'echet topology is canonically defined by
 \ref{_coherent_Frechet_Proposition_}. It is not
hard to see that any coherent sheaf morphism is
continuous in the Fr\'echet topology as well
(\cite[\S VIII.A.9]{_Gunning_Rossi_}).

\hfill

\theorem\label{_coherent_sheaves_Montel_Theorem_}
 (\cite[Theorem VIII.A.9]{_Gunning_Rossi_})\\
All coherent sheaves on complex analytic spaces are
Montel.
\endproof

\subsection{Banach structure on the space of sections}

\remark\label{_norm_on_quotient_Remark_}
For a coherent sheaf ${\cal F}$ on $M$ and $x\in M$,
let $r_x({\cal F}):= \dim_\C \frac{\cal F}{\goth m_x{\cal F}}$,
where $\goth m_x$ is the maximal ideal of $x$.
On a Stein variety $M$, any coherent sheaf ${\cal F}$, for which
$r_x({\cal F})$ is bounded, is obtained
as a quotient of a free sheaf, ${\cal F}=\calo_M^n/{\cal H}$
(\cite[Theorem 4.3]{_Forster:Steinischen_}).
Let $\phi:\; \calo_M^n\to {\cal F}$ be the quotient map.
For such a quotient, the seminorms $\|\cdot \|_K$ 
of \ref{_coherent_Frechet_Proposition_}
can be constructed explicitly: for any
$f\in H^0({\cal F})$, we write 
$\|f\|_K:= \inf_{z_1, ..., z_n}\sum_{i=1}^n \|z_i\|_K$,
where the infimum is taken over the set of all representatives
$z_1, ..., z_n\in H^0(\calo_M)$
such that $\phi(z_1, ..., z_n)=f$.

\hfill

\remark\label{_seminorm_is_norm_Remark_}
In general, $\|\cdot\|_K$
is only a seminorm. Indeed, let ${\cal F}$
be the sum of two skyscraper sheaves, ${\cal F}=k_x\oplus k_y$, 
and $K\subset M$ a subset containing $x$ but not $y$.
Then, clearly, $ \|\cdot\|_K$ vanishes on the summand
$k_y$, and is non-zero on $k_x$.
However, $\|\cdot\|_K$ is a norm whenever
the support of any non-zero section $f\in H^0({\cal F})$ 
intersects the interior of $K$ (for example,
when ${\cal F}$ is torsion-free).

\hfill

\definition
Let ${\cal F}$ be a coherent sheaf on a Stein variety,
obtained as a quotient ${\cal F}=\calo_M^n/{\cal H}$
of a free sheaf. Fix an open set $U\Subset M$,
and assume that  the support of any section of ${\cal F}$
intersects with $U$. Define
$\|\cdot \|_U:= \inf_{z_1, ..., z_n}\sum_{i=1}^n \sup_U |z_i|$,
where the infimum is taken, as in \ref{_norm_on_quotient_Remark_}, 
over all representatives $z_1, ..., z_n\in H^0(\calo_M)$
such that $\phi(z_1, ..., z_n)=f$. A section of ${\cal F}$
is called {\bf bounded} if $\|f \|_U$ is finite.

\hfill

\remark
Since the space of sections of ${\cal H}$
is closed in $H^0({\cal F})$, the map
$f\to \|f \|_U$ defines a norm on the space
$H^0_b({\cal F})$ of bounded sections
of ${\cal F}$.\footnote{Generally speaking, $\|\cdot \|_K$
is a seminorm, but \ref{_seminorm_is_norm_Remark_} implies
that $\|\cdot \|_U$ is zero only on sections with support
outside of the interior of $M$, which means zero sections.}

\hfill

The following assertion is true for all coherent sheaves,
but it is hard to state it clearly without the restrictive
assumptions which we used. Note that 
the sheaf ${\cal F}$ is a quotient of $\calo_M^n$
if for some $C>0$, one has
$r_x({\cal F}) <C$ for all $x\in M$
(\ref{_norm_on_quotient_Remark_}).

\hfill

\proposition\label{_H^0_b(F)_Banach_Proposition_}
Let ${\cal F}$ be a coherent sheaf on a Stein variety $M$.
Assume that ${\cal F}=\calo_M^n/{\cal H}$,
where ${\cal H}\subset \calo_M^n$ is a 
coherent subsheaf. Consider a relatively compact
subset $U\subset M$, and let $H^0_b(U,{\cal F})$
denote the space of sections which can be
extended to a neighbourhood of $\bar U$.
Then the norm $\|\cdot \|_{\bar U}$ is 
complete on $H^0_b(U,{\cal F})$.

\hfill

\proof Since $M$ is Stein, 
$H^0_b(U,{\cal F})= \frac{H^0_b(U, \calo_U^n)}{H^0_b(U, {\cal H})}$.
By construction, $H^0_b(U,{\cal F})$ is a 
quotient of a Banach space $H^0_b(U,\calo_U)^n$ 
(\ref{_bounded_functions_Banach_Theorem_}), 
and $\|\cdot \|_{\bar U}$  is the quotient norm.
\endproof

\subsection{The Fr\'echet topology and the adic topology}

\definition
Let ${\cal F}$ be a coherent sheaf on $M$, $U\subset M$
an open set, and $x\in U$ a point. 
Denote by ${\cal F}_x$ the space of germs of ${\cal F}$ in
$x$, and let ${\goth m}_x$ denote the maximal ideal.
Consider the map
$\phi_r:\; H^0(U, {\cal F})\to \frac{{\cal F}_x}{\goth m_x^r{\cal    F}}$
associating to a section its $r$-jet in $x$.
The {\bf ${\goth m}_x$-adic topology}, or simply
{\bf adic topology} on $H^0(U, {\cal F})$ is the weakest topology such that
all maps $\phi_r$ are continuous.

\hfill

\remark\label{_proj_to_adic_continuous_Remark_}
The following proposition implies that 
this topology is weaker than the Fr\'echet topology
defined previously.

\hfill

\proposition\label{_jet_continuous_Proposition_}
Let ${\cal F}:= \calo_M^n/{\cal H}$ be a coherent sheaf on
a Stein manifold $M$, and $H^0_b(U, {\cal F})$ the space of 
bounded sections, equipped with the Banach structure as in
\ref{_H^0_b(F)_Banach_Proposition_}.
Then the map 
$\phi_r:\; H^0_b(U, {\cal F})\to \frac{{\cal F}_x}{\goth m_x^r{\cal    F}}$
defined above is continuous.

\hfill

\proof
Using the Cauchy formula, we can express the partial
derivatives of a function at a point through a Cauchy integral over
a contour. This implies that the map
$\phi_r:\; H^0_b(U, {\cal F})\to \frac{{\cal F}_x}{\goth  m_x^r{\cal    F}}$
is continuous when ${\cal F}=\calo_M^n$.
When ${\cal F}:= \calo_M^n/{\cal H}$, the space
$H^0_b(U, {\cal H})$ is by construction closed
in $H^0_b(U, \calo_M^n)$; this implies that
$\phi_r:\; H^0_b(U, {\cal F})\to \frac{{\cal F}_x}{\goth  m_x^r{\cal    F}}$
is continuous as well.
\endproof


\section{GAGA for coherent sheaves over cone varieties}
\label{_GAGA_coherent_Section_}


\subsection{Coherent sheaves over cone varieties: main result}

In \cite{_OV_Algebraic_Cones_}, 
we proved a version of GAGA theorem for 
smooth open cone varieties, showing that the algebraic
structure on a smooth open cone is uniquely determined
by the analytic structure and the holomorphic contraction map
$\gamma$. This theorem was extended to  all cone varieties
in \ref{_algebraic_str_on_cone_uniqueness_Theorem_}.

In this section, we extend  \ref{_algebraic_str_on_cone_uniqueness_Theorem_} to
$\langle \gamma\rangle$-equivariant reflexive coherent sheaves
on cone varieties. The notion of an equivariant coherent
sheaf is quite complex. We already defined equivariant
structures on line bundles in \ref{_equiv_on_line_bundle_Definition_};
the same definition works nicely for holomorphic
vector bundles. Most general definitions will be introduced
in Subsection \ref{_equivariant_coherent_Subsection_}.
For the present purposes, it is enough to define
the $\Z$-equivariant structure on a coherent sheaf.
\ref{_equiv_Z_Definition_}
is compatible with \ref{_equiv_on_line_bundle_Definition_}
and with the more general definitions introduced in 
Subsection \ref{_equivariant_coherent_Subsection_}.

\hfill

\definition\label{_equiv_Z_Definition_}
Let $\gamma$ be an automorphism of a complex variety $X$,
and $\langle\gamma\rangle\cong\Z$ the corresponding cyclic group.
{\bf A $\langle\gamma\rangle$-equivariant structure} on a coherent
sheaf ${\cal F}$ is an isomorphism ${\cal F} \to \gamma^* {\cal F}$,
that is, an isomorphism 
$\phi:\; H^0(U, {\cal F})\to H^0(\gamma^{-1}(U), {\cal F})$
which satisfies $\phi(a f) = \gamma^*(a) \phi(f)$
for any open $U\subset M$ and any 
 $a\in H^0(\calo_U)$, $f\in H^0(U, {\cal F})$. We also  want $\phi$ to be
compatible with restrictions to smaller open subsets $U'\subset U$.
In the future, we will usually write $\gamma^*$
in place of $\phi$, to make the formulas for coherent sheaves
 consistent with those for holomorphic functions.

\hfill

\theorem\label{_coherent_GAGA_Theorem_}
Let $(X, \gamma)$ be a closed cone variety.
Using \ref{_algebraic_str_on_cone_uniqueness_Theorem_}, 
we consider $X$ as an affine variety.
Consider a $\langle \gamma\rangle$-equivariant analytic reflexive coherent 
sheaf ${\cal F}$ on $X$. Then ${\cal F}$
admits an algebraic structure,\footnote{We 
give a formal definition of an algebraic structure
on an analytic coherent sheaf in 
\ref{_alg_structure_Definition_}.} which is 
uniquely determined by the 
$\gamma$-equivariant structure.

\hfill

We prove \ref{_coherent_GAGA_Theorem_}
in Subsection \ref{_GAGA_coherent_Subsection_}.

%

\subsection{$\gamma$-finite sections of a $\gamma$-equivariant coherent sheaf}

As in \ref{_coherent_Frechet_Proposition_}, 
we define the Fr\'echet structure on
the sheaf ${\cal F}$, using the seminorms $\|\cdot\|_K$
on $H^0(X, {\cal F})$ associated with compact subsets
$K\subset X$. To simplify the exposition,
we assume that each $K$ is obtained as 
the closure of an open subset.  When the sheaf ${\cal F}$ 
is torsion-free, and $X$ is irreducible,
these seminorms are in fact norms
(\ref{_seminorm_is_norm_Remark_}).

Given $U\Subset X$, we denote by $H^0_b(U, {\cal F})$
the space of $\|\cdot\|_{U}$-bounded sections of ${\cal F}$ on $U$;
by \ref{_H^0_b(F)_Banach_Proposition_}, this space is Banach.

\hfill

\theorem\label{_pullback_on_coherent_Theorem_}
Let $X$ be a 
closed cone variety, and $\gamma:\; X \to X$
an invertible holomorphic contraction of $X$ to the origin.
Consider a $\langle\gamma\rangle$-equivariant
analytic coherent sheaf ${\cal F}$ on $X$, obtained as
a quotient of a free sheaf,  and let
$U\subset X$ be an open set such that $\gamma(U)\Subset U$. 
Using the equivariant structure, we get an isomorphism
$H^0_b(U, {\cal F})\to H^0_b(\gamma^{-1}(U), {\cal F})$.
Then, composing this map with the restriction to $U\subset \gamma^{-1}(U)$,
we obtain an operator 
\begin{equation}\label{_gamma^*_on_sec_Equation_}
\gamma^*:\; H^0_b(U, {\cal F})\to H^0_b(U, {\cal F})
\end{equation}
on the corresponding Banach spaces 
which is compact and satisfies the assumptions of 
\ref{_bounds_by_C^n_Claim_}.

\hfill

\pstep
Let $U\Subset X$ be an open, relatively
compact subset such that $\gamma(U) \Subset U$,
and $\|\cdot\|_{\bar U}$ the corresponding norm.
By \ref{_coherent_sheaves_Montel_Theorem_}, the restriction map 
\[ 
(H^0(U, {\cal F}), \|\cdot \|_{\bar U})\longrightarrow 
(H^0(\gamma(U), {\cal F}), \|\cdot \|_{\gamma(\bar U)})
\]
is compact (as a morphism of topological vector spaces).

\hfill

{\bf Step 2:} 
The pullback map
$\gamma^*:\; H^0_b(U, {\cal F})\to H^0_b(U, {\cal F})$ 
can be obtained as a composition of the isomorphism
$\gamma^*: H^0_b(U, {\cal F})\to H^0(\gamma^{-1}(U), {\cal F})$
and the restriction $H^0(\gamma^{-1}(U), {\cal F})\to H^0_b(U, {\cal F})$,
which is compact by Step 1. Therefore, 
$\gamma^*:\; H^0_b(U, {\cal F})\to H^0_b(U, {\cal F})$
is compact. 

This proves the first claim of 
\ref{_pullback_on_coherent_Theorem_};
it remains to show that this map satisfies
the assumptions of \ref{_bounds_by_C^n_Claim_}.
We follow the same logic as used in the
proof of \ref{_Banach_no_quasinilpotent_Proposition_}.

\hfill

{\bf Step 3:}
Let $\goth m$ be the maximal ideal of $x$.
 {\bf The $k$-jet of a section of ${\cal F}$ in $x$} is its class
in $\frac{\cal F}{{\goth m}^{k+1}{\cal F}}$. Since the ring of 
germs of holomorphic functions is Noetherian,
Krull lemma implies that $\bigcap_n {\goth m}^n{\cal F}=0$.
Therefore, any non-zero section of ${\cal F}$
has non-zero $k$-jet for some $k>0$.

\hfill

{\bf Step 4:}
Consider the map which associates to a section $f\in H^0_b(U, {\cal F})$
its $k$-jet $f_k$ in the origin $x$ of the cone variety 
$M$. By \ref{_jet_continuous_Proposition_}, this map
is continuous.
We are going to show that $\lim\sup_n |C^n (\gamma^*)^nf_k|> 0$
for appropriate $C$ and $k$. 

\hfill

{\bf Step 5:} 
Let $\mu$ be the eigenvalue with the smallest absolute value of $\gamma^*$
on the Zariski cotangent space ${\goth m}/{\goth m}^2$
and $\lambda$ the eigenvalue with the smallest absolute value of $\gamma^*$
on $\frac{{\cal F}}{{\goth m}{\cal F}}$.
The natural map 
$\frac{{\cal F}}{{\goth m}{\cal F}}\otimes {\goth m}^k/{\goth m}^{k+1}
\to \frac{{\goth m}^k{\cal F}}{{\goth m^{k+1}}{\cal F}}$
is surjective by Nakayama's lemma.
Therefore, the eigenvalue with the smallest absolute value of $\gamma^*$
acting on the space of $k$-jets is bounded from below by $\lambda \mu^k$.
This implies that $\lim\sup C^n|(\gamma^*)^nf_k| >0$ for 
$C> \lambda^{-1}\mu^{-k}$,
and the sequence $C^n (\gamma^*)^nf$ cannot converge to 0
in $H^0_b(U, \calo_U)$, when $k$ is the smallest integer
such that $f_k \neq 0$.
\endproof

\hfill

Applying the definition of the equivariant action,
we obtain an isomorphism 
$H^0(X, {\cal F})\to H^0(X, \gamma^*{\cal F})=
H^0(X, {\cal F})$, also denoted $\gamma^*$. 
Restricted to a smaller open set $U$, this map
is equal to the arrow in \eqref{_gamma^*_on_sec_Equation_}. 
By definition, a section $f\in H^0(X, {\cal F})$  
is {\bf $\gamma$-finite} if 
the sequence $f, \gamma^*(f), \gamma^*(\gamma^*(f)), ...$
generates a finite-dimensional subspace in $H^0(X, {\cal F})$.
Applying \ref{_Riesz_Schauder_Theorem_}, we immediately
obtain the following result.

\hfill

\corollary\label{_finite_dense_in_bounded_Corollary_}
Let ${\cal F}$ be a $\gamma$-equivariant coherent sheaf
on a cone variety $(X, \gamma)$. Assume that ${\cal F}$
is a quotient of a free sheaf. Then the space
of $\gamma$-finite sections of ${\cal F}$ is
dense in $H^0_b(U, {\cal F})$.
\endproof

\hfill

\corollary\label{_finite_dense_in_H^0_Corollary_}
Let ${\cal F}$ be a $\gamma$-equivariant coherent sheaf
on a cone variety $(X, \gamma)$. Assume that ${\cal F}$
is a quotient of a free sheaf. Then the space
of $\gamma$-finite sections of ${\cal F}$ is
dense in $H^0(X, {\cal F})$ with the natural Fr\'echet topology.

\hfill

\proof
The Fr\'echet topology is obtained by taking the 
restriction to $U\Subset X$ and using the
Banach norm on $H^0_b(U, {\cal F})$.
By \ref{_finite_dense_in_bounded_Corollary_}, the $\gamma$-finite sections
are dense in every open set of $H^0(X, {\cal F})$.
\endproof

\subsection{The completion map}

The following theorem is implied by the same
argument as used to prove \ref{_finite_function_extends_Proposition_}.

\hfill

\theorem\label{_completion_extension_finite_Theorem_}
Let ${\cal F}$ be a $\gamma$-equivariant
coherent sheaf on a closed cone variety $(X, \gamma)$,
and ${\cal F}_x$ the space of germs of sections of ${\cal F}$
in the origin $x\in X$. Denote by $\widehat{\cal F}_x$
the completion of ${\cal F}_x$ in $x$,
\[\widehat{\cal F}_x:= 
\lim\limits_{\longleftarrow}\frac{{\cal F}_x}{\goth m^k{\cal F}_x},
\]
where $\goth m\subset \calo_X$ is the maximal ideal of $x\in X$.
 Then:
\begin{description}
\item[(i)] Any $\gamma$-finite germ $f\in {\cal F}_x$ 
can be represented by a $\gamma$-finite section of ${\cal F}$ on $X$.
\item[(ii)] The completion map defines a bijection
between the spaces of $\gamma$-finite vectors in ${\cal F}_x$
and $\widehat{\cal F}_x$.
\end{description}

\pstep
The argument in this step repeats the one
used to prove \ref{_finite_function_extends_Proposition_}.
Let $f\in  {\cal F}_x$ be a $\gamma$-finite germ. Then
the space 
$W:=\langle f, \gamma^* f, (\gamma^*)^2 f,..., (\gamma^*)^{r-1} f\rangle$
is $\gamma^*$-invariant, for some $r\in \Z^{>0}$. 
We denote the germ $(\gamma^*)^k f$ by $f_{(k)}$.
Let $U\subset X$ be an open set
containing the origin, such that all $f_{(i)}$ are holomorphic
in $U$ and $\gamma(U) \Subset U$; such an open set 
clearly exists. The $\gamma^*$-invariance of $W$ is then equivalent
to an equation 
$((\gamma^{-1})^*)^r(f)=\sum_{i=0}^{r-1} \alpha_i f_{(i)}$,
for some $\alpha_i \in \C$, or, equivalently,
\[
f=  (\gamma^*)^r\sum_{i=0}^{r-1} \alpha_i f_{(i)}.
\]
Since all $f_{(i)}$ are sections of $H^0(U, {\cal F})$,
the expression $(\gamma^*)^r\sum_{i=0}^{r-1} \alpha_i f_{(i)}$
is a section of $H^0(\gamma^{-r}(U), {\cal F})$.
Therefore, $f$ is the restriction of a section in 
$H^0(\gamma^{-r}(U), {\cal F})$. Since $\gamma$ is a
contraction, $\bigcup_i \gamma^{-ri}(U) =X$,
and $f$ is obtained as the restriction of a section in 
$H^0(X, {\cal F})$. This proves 
\ref{_completion_extension_finite_Theorem_} (i).

\hfill

{\bf Step 2:}
It remains to prove \ref{_completion_extension_finite_Theorem_} (ii).
The map  $H^0_\gamma(X, {\cal F})\to \widehat{\cal F}_x$
is injective by Krull's intersection theorem (its kernel is 
$\bigcap_k \goth m_x^k{\cal F}$). To prove its surjectivity
on $\gamma$-finite vectors, fix an eigenvalue $\alpha$
of $\gamma^*$ acting on $\widehat{\cal F}_x$, and 
$V_\alpha\subset H^0_\gamma(X, {\cal F})$
the corresponding eigenspace, which is clearly 
finite-dimensional. 

The rest of the argument is similar to the one which proves
\cite[Proposition 5.6]{_Verbitsky:HC_}.
Denote by $H^0_\gamma(X, {\cal F})$ the space of $\gamma$-finite sections.
By \ref{_finite_dense_in_H^0_Corollary_},
$H^0_\gamma(X, {\cal F})$ is dense in $H^0(X, {\cal F})$.
Since the map $H^0(X, {\cal F})\to \frac{\cal F}{\goth m_x^k}$
is continuous (\ref{_jet_continuous_Proposition_}),
and a dense subspace of a finite-dimensional space $W$ coincides with $W$,
we obtain that the map 
$H^0_\gamma(X, {\cal F})\to \frac{\cal F}{\goth m_x^k{\cal F}}$
is also surjective. 
Since $H^0_\gamma(X, {\cal F})\to \frac{\cal F}{\goth m_x^k{\cal F}}$
is surjective, it is surjective on the eigenspace corresponding
to $\alpha$, which is identified with $V_\alpha$.
\endproof

\subsection{GAGA-type theorem for equivariant 
coherent sheaves on closed cone varieties}
\label{_GAGA_coherent_Subsection_}

We are going to prove that any $\langle \gamma\rangle$-equivariant reflexive 
analytic coherent sheaf on a cone variety is in fact algebraic;
this statement is a coherent sheaf analogue of
\ref{_algebraic_str_on_cone_uniqueness_Theorem_},
where the same result is proven for the cone varieties.

\hfill

\definition\label{_alg_structure_Definition_}
Let ${}^a\!X$ be an algebraic variety, $X$ the
same variety considered as a complex analytic variety,
and ${\cal F}$ an analytic coherent sheaf on $X$.
Let ${}^a\!\calo_X\subset \calo_X$ denote the sheaf of regular 
holomorphic functions on $X$. 
An {\bf algebraic structure} on ${\cal F}$
is a subsheaf ${}^a\!{\cal F}\subset {\cal F}$ 
with the following properties.
\begin{enumerate}[(a)]
\item
${}^a\!{\cal F}$ is a sheaf of
${}^a\!\calo_X$-modules, that is, a section of ${}^a\!{\cal F}$
multiplied by a section of ${}^a\!\calo_X$ belongs to ${}^a\!{\cal F}$.
\item
${}^a\!{\cal F}$ is locally finitely generated
over ${}^a\!\calo_X$ in the Zariski topology, that is, 
${}^a\!{\cal F}$ is coherent
as a sheaf on ${}^a\!X$ (in algebraic category).
\item ${\cal F}= {}^a\!{\cal F}\otimes_{{}^a\!\calo_X} \calo_X$.
\end{enumerate}

\hfill

\theorem\label{_GAGA_reflexive_Theorem_}
Let $X$ be a closed cone variety, and $\gamma:\; X \to X$
an invertible holomorphic contraction of $X$ to the origin
commuting with the natural $\C^*$-action on the cone.
By \ref{_algebraic_str_on_cone_uniqueness_Theorem_}, 
$\gamma$ determines the algebraic structure on $X$.
Consider a reflexive $\langle\gamma\rangle$-equivariant
analytic coherent sheaf ${\cal F}$ on $X$. Then
${\cal F}$ admits a natural algebraic structure,
such that the corresponding algebraic coherent sheaf
${}^a\!{\cal F} \subset {\cal F}$ is also
$\langle\gamma\rangle$-equivariant. 
Moreover, $H^0(X,{}^a\!{\cal F})$ is the 
space of $\gamma^*$-finite sections of ${\cal F}$.

\hfill

\pstep
Let ${}^a\!\calo_X$ be the sheaf of regular
(that is, algebraic) functions on $X$.
By \ref{_algebraic_str_on_cone_uniqueness_Theorem_},
the sections of ${}^a\!\calo_X$
are the same as $\gamma^*$-finite functions.

Let ${}^a\!{\cal F}\subset {\cal F}$ be 
the sheaf of ${}^a\!\calo_X$-modules
generated by all $\gamma^*$-finite sections of ${\cal F}$.
By \ref{_finite_dense_in_H^0_Corollary_}, 
the space of sections of ${}^a\!{\cal F}$ is dense in 
the corresponding space of sections of
${\cal F}$. It remains to show that ${}^a\!{\cal F}$
is finitely generated over ${}^a\!\calo_X$.

\hfill

{\bf Step 2:}
Since $X$ is affine (\ref{_cone_str_on_open_cone_variety_Theorem_}), 
to show that ${}^a\!{\cal F}$ is finitely generated over ${}^a\!\calo_X$
it suffices to prove that $H^0(X, {}^a\!{\cal F})$
is finitely generated as an $H^0(X, {}^a\!\calo_X)$-module.
However, $H^0(X, {}^a\!{\cal F})$ is identified
with the space of $\gamma^*$-finite elements
in the completion $\widehat {\cal F}_x$ 
(\ref{_completion_extension_finite_Theorem_}),
and $H^0(X, {}^a\!\calo_X)$ with the space
of $\gamma^*$-finite elements of the
adic completion of $H^0(X, \calo_X)$.
By \ref{_algebraization_R_modules_Theorem_},
the former is finitely generated over the latter.
\endproof 

\hfill

\remark
Clearly, \ref{_GAGA_reflexive_Theorem_}
defines an equivalence of categories 
between $\gamma$-equivariant complex analytic and 
 algebraic coherent reflexive sheaves on $X$.

\hfill

\remark\label{_formal_used_here_Remark_}
Most of this paper deals with the
topology on functions induced by the 
uniform convergence, obtaining an
algebraization of the (complex analytic) cone varieties
and coherent sheaves on the cone varieties.
In Sections
\ref{_algebraiz_for_rings_Section_} 
and \ref{_modules_over_complete_Section_} we deal with the adic topology,
obtaining an algebraization of a complete local ring
$R$ equipped with a contraction $\gamma$ and of 
any $\gamma$-equivariant $R$-module. 
In the essense, in Sections
\ref{_algebraiz_for_rings_Section_} 
and \ref{_modules_over_complete_Section_}
we prove the same results as in the
rest of this paper, but in the category
of formal varieties instead of the
complex analytic category. These
two sections are separate from the rest of this
paper; however, the algebraization theorem
from Section \ref{_modules_over_complete_Section_}
is necessary in the second step of the
proof of \ref{_GAGA_reflexive_Theorem_}.
This is the only place in this paper
where we use the formal completions or
any result from Sections \ref{_algebraiz_for_rings_Section_} 
and \ref{_modules_over_complete_Section_}.


\section{Equivariant coherent sheaves on orbivarieties}
\label{_orbivarieties_equiv_sheaves_Section_}


\subsection{Equivariant coherent sheaves}
\label{_equivariant_coherent_Subsection_}

In this section we define the notion of an equivariant structure on a
coherent sheaf. Let $G$ be a
group acting on a variety $M$, and ${\cal F}$ a coherent sheaf,
algebraic or analytic depending on the context. 
When ${\cal F}$ is a bundle, we can just define
the equivariant structure on ${\cal F}$ as the
lifting of the action of $G$ to the total space
$\Tot({\cal F})$ acting linearly in the fibers.
For a more general situation, this definition will
not work. Here we have to follow a more traditional
route, in the spirit of \cite{_SGA:3_1_}, where this
notion is defined in the more general context
of $G$-objects in a category. However, the intuition
should follow the geometric route: a $G$-equivariant
structure on a sheaf is a fiberwise linear action
of $G$ on its ``total space''.

We have already defined equivariant coherent sheaves 
for the $\Z$-action (\ref{_equiv_Z_Definition_}); this definition is compatible
with the present one.

\hfill

\definition\label{_equivariant_sheaf_Definition_}
Let $\cac$ be a category, and $G$ a group acting on $\cac$ by endofunctors,
that is, by functors from $\cac$ to itself.  
{\bf A $G$-equivariant structure
on an object ${\cal F}$ of $\cac$} is a family of isomorphisms
$R_g:\; g^*({\cal F}) \to {\cal F}$, for all $g\in G$ that is  associative in the
following sense:
\[ g_2^*(R_{g_1}) \circ R_{g_2}= R_{g_1 g_2},\] where
$R_{g_1 g_2}$ is considered to be  an isomorphism
${(g_1g_2)}^*{\cal F} \to {\cal F}$, and $g_2^*(R_{g_1}) \circ R_{g_2}$
is the composition of the isomorphism
$g_2^*(R_{g_1}):\; g_2^*g_1^* ({\cal F})\to g_2^*({\cal F})$
and $R_{g_2}:\; g_2^*({\cal F}) \to {\cal F}$, for any $g_1, g_2\in G$.

\hfill

\definition\label{_G_equiv_sheaf_Definition_}
Let $G$ be a group acting on a  complex variety
$X$ by automorphisms. Then $G$ acts on the category of coherent sheaves
on $X$ by endofunctors. {\bf A $G$-equivariant coherent sheaf}
on $X$ is a $G$-equivariant object in the category of 
coherent sheaves.

\hfill

\remark
Note that in this definition we consider $G$ as an abstract group.
In many applications, one needs to consider $G$ as a group object
in $\cac$, and then a $G$-equivariant structure is a
morphism $G\times {\cal F}\to {\cal F}$, where $\times$ denotes the
fibred product in $\cac$. This is the definition used in 
\cite{_SGA:3_1_}.

\hfill

\definition (\cite[\S 39.12]{_Stacks_Project_})
Let $G$ be an algebraic group and $a:\; G\times X\to X$
its action on an algebraic variety. Denote by $\pi_2:\;  G\times X\to X$
the projection map. Consider a coherent sheaf
${\cal F}$ on $X$. {\bf An equivariant structure} on ${\cal F}$
is an isomorphism $a^* {\cal F}\stackrel \phi \to\pi_2^*{\cal F}$
such that the following diagram is commutative
\[ 
\xymatrix{ (1_ G \times a)^*\pi_1^*\mathcal{F} \ar[r]_-{\pi_{23}^*\phi } &  \pi_2^*\mathcal{F} \\  (1_ G \times a)^*a^*\mathcal{F} \ar[u]^{(1_ G \times a)^*\phi } \ar@{=}[r] &  (m \times 1_ X)^*a^*\mathcal{F} \ar[u]_{(m \times 1_ X)^*\phi } }.
\]
In this diagram, $1_G$ and $1_X$ are identity maps, $m:\; G \times G\to G$
is the multiplication map, and $\pi_{23}:\; G\times G\times X \to G\times X$ 
the projection.


\hfill

In the sequel, we will be interested in $\C^*$-equivariant
coherent sheaves on a space $X$ equipped with a $\C^*$-action.
When this action is free and proper, there exists
a quotient space $Y$, and the category of $\C^*$-equivariant
coherent sheaves on $X$ is identified with the category
of coherent sheaves on $Y$, by 
\ref{_Roberts_invariant_sections_Theorem_} below
(see also \cite{_Roberts:G_sheaves_,_Kollar:quotient_}).
However, when the action is not free in a special point, 
such equivalence does not hold, even for reflexive sheaves.
This is where we need to use the orbispaces,
introduced below in Subsection \ref{_orbispaces_Subsection_}. 

This phenomenon
is apparent from the example given in Subsection
\ref{_orbifold_sheaf_example_Subsection_}.

\subsection{Equivariant coherent sheaves on $\C^2\backslash 0$: an example}
\label{_orbifold_sheaf_example_Subsection_}

The last sections of this paper are based on the following
observation, originally due to M. Roberts
\cite{_Roberts:G_sheaves_}. Let $M$ be a complex
variety equipped with a proper, free action of a complex
reductive Lie group $G$. Then $G$-equivariant
coherent sheaves on $M$ are the same as coherent
sheaves on the quotient $M/G$. We apply this
result only to $G=\C^*$. For algebraic groups
(not necessarily reductive) acting on quasiprojective varieties,
a stronger framework exists, due to J. Koll\'ar
(\cite{_Kollar:quotient_}). For most applications
in the present paper, the earlier complex analytic
framework suffices, though the algebraic version
is more powerful (especially when the action is not free).

Let $G$ be a reductive algebraic group and
$X$ a complex algebraic variety equipped with a $G$-action.  Following the
conventions of \cite{_Kollar:quotient_}, we say that
the action $\rho:\; G\times X \to X$ is {\bf proper} if 
$\rho\times \Id:\; (G, X)\to (X,X)$ is
proper. In \cite[Theorem 1.5]{_Kollar:quotient_} 
it is shown that the geometric quotient $X/G$ is a complex
variety when $G$ acts properly. 

Given a $G$-equivariant coherent sheaf $F$ on $X$, we
produce the sheaf $F^G$ of $G$-invariant sections on
$X/G$. The  pushforward of $F$ to $X/G$ is coherent, as Koll\'ar shows in 
\cite[Theorem 3.12]{_Kollar:quotient_}.
Clearly, $\pi_*(\pi^*(H)^G) =H$ for any $H$ on $X/G$. However, if
$F$ is a sheaf on $X$, the natural morphism 
$\pi^*((\pi_*F)^G)\to F$ is not necessarily an isomorphism, 
even if $F$ is reflexive, $X$
is normal and $G=\C^*$, unless the singular locus of $G$ has
codimension $\geq 2$. We give an explicit counterexample,
related to the constructions used further in this section.

\hfill

Consider the ${\Bbb C}^*$-action $\rho$ on $X= {\Bbb
  C}^2\backslash 0$ such that $\lambda\in {\Bbb C}^*$
takes $(x, y)$ to $(\lambda x,\lambda^2 y)$. This action
is the composition of the ${\Bbb  Z}/2$-action taking $(x, y)$ 
to $(-x, y)$ and of the $({\Bbb C}^*/\pm 1)$-action  
$(x, y)\to (\lambda^2 x,\lambda^2 y)$, from the
exact sequence 
\[0\to {\Bbb Z}/2 \to {\Bbb C}^* \to {\Bbb C}^*\to 0.
\] 
The $({\Bbb C}^*/\pm 1)$-action $(x, y)\to (\lambda^2 x,\lambda^2 y)$ 
is free, hence the category of
$\rho$-invariant coherent sheaves on $X$ is equivalent to
the category of ${\Bbb Z}/2$-invariant sheaves on 
${\Bbb  C} P^1$, where ${\Bbb Z}/2$ acts via 
$(x, y) \to (-x,y)$. At this point we can 
switch from the group $G={\Bbb C}^*$ to 
$G={\Bbb Z}/2$ and from $X={\Bbb C}^2\backslash 0$ 
to $X={\Bbb C} P^1$, and show that the
functor $F\mapsto F^G$ is not an equivalence in this case;
it won't be an equivalence for the original case as well.

It remains to show that there are $G$-equivariant
reflexive sheaves on ${\Bbb C} P^1$ which are not
pullbacks of reflexive sheaves on 
${\Bbb  C} P^1/G$. Consider the trivial sheaf $\cal O$, with the
$G$-equivariant action $f \mapsto - f$. The corresponding
sheaf of $G$-invariant sections on ${\Bbb C} P^1/G$ is the sheaf of 
functions which have zeros in two fixed points of the 
$G$-action, and hence it is isomorphic to ${\cal
  O}(-2)$. Then its pullback to ${\Bbb C} P^1$ is ${\cal
  O}(-4)$, and the map $\pi^*(\pi_*(F^G))\to F$ is a map from
${\cal O}(-4)$ to ${\cal O}$, manifestly not an
isomorphism.

\subsection{The orbispaces, a. k. a. Deligne-Mumford stacks}
\label{_orbispaces_Subsection_}

We start by introducing the orbifolds, following
\cite{_Lerman:Orbifolds_}. 
Then we explain how to generalize this definition 
to (possibly singular) complex varieties, arriving at the notion
of the orbispace, or the (topological) Deligne-Mumford stack.

\hfill

\definition
{\bf A groupoid} is a category with all
morphisms invertible. A groupoid is {\bf finite}
if $\Mor(X,Y)$ is finite for all objects $X, Y$.

\hfill

\example A group is a groupoid with just one object.

\hfill

\example Let $M$ be a manifold or a CW-space.
Its {\bf Poincar\'e grou\-po\-id} is the category
with objects the points of $M$, and morphisms $\Mor(x,y)$
the homotopy classes of paths from $x$ to $y$.

\hfill

\definition
Let $G_1$ be the set of all morphisms
in a category, and $G_0$ the set of all objects.
The {\bf source  map} 
$s:\; G_1 \to G_0$ takes $x\in \Mor(A,B)\subset G_1$
to $A\in G_0$, and the {\bf target map} 
$t:\; G_1 \to G_0$ takes $x$ to $B\in G_0$.

\hfill

The notion of groupoid already contains the
group operations (the composition of morphisms and
taking the inverse), also called {\bf the structure maps.}
However, it is possible to give another definition
of groupoid, defining the groupoid axiomatically in terms
of the structure maps. Presently we will formulate the groupoid axioms
using this language. This gives an alternative definition,
better suited for the geometric constructions ahead.

Denote by $G_1 {}_s\times_t G_1$ the set of all
$(\phi, \psi)\in G_1\times G_1$ such that
$t(\phi)=s(\psi)$.
There are three structure maps, {\bf the composition}
$\mu:\; G_1 {}_s\times_t G_1\to G_1$, {\bf the inverse}
$\alpha:\; G_1 \to G_1$, and {\bf the unit} $\iota:\; G_0\to G_1$ which
are almost self-evident.

The axioms of the unit: $s(\iota(A))=t(\iota(A)=A$, 
$\mu(\iota(A), \phi)=\phi$ when $A=s(\phi)$, and
$\mu(\phi, \iota(B))=\phi$ when $B=t(\phi)$.

The axioms of inverse:
$s(\alpha(\psi))=t(\psi)$,  $t(\alpha(\psi))=s(\psi)$, 
and $\mu(\psi, \alpha(\psi))=\iota(s(\psi))$,
 $\mu(\alpha(\psi), \psi)=\iota(t(\psi))$.

Associativity: the maps 
$(\Id\times \mu)\circ \mu$ and 
$(\mu\times \Id)\circ \mu:\; G_1 {}_s\times_t G_1 {}_s\times_t G_1\to G_1$
are equal.

\hfill

\definition
A {\bf topological groupoid} is
a groupoid where $G_0$ and $G_1$ are topological
spaces, and the structure maps are continuous.
The {\bf coarse  space} or {\bf the orbit space} 
of a topological groupoid is the quotient of $G_0$ by the
equivalence relation generated by all pairs $(x, y)\in G_0\times G_0$
such that $\Mor(x,y)$ is non-empty.

\hfill

\definition
A {\bf Lie groupoid} is
a groupoid where $G_0$ and $G_1$ are smooth manifolds, 
the source and target maps 
are smooth submersions, and the rest of the
structure maps are also smooth. 

\hfill

\definition
A Lie groupoid is called
{\bf an \'etale groupoid} if the source 
and the target maps are local diffeomorphisms.
It is called {\bf proper} if the map
$G_1 \to G_0\times G_0$ sending a morphism
to a pair (source, target) is proper.

\hfill

\example\label{_action_groupoid_Example_}
Let $G$ be a Lie group acting on a manifold $M$.
Define its {\bf action groupoid} as follows:
the objects of the action groupoid are the points of $M$,
and the set of morphisms $\Mor(x,y)$ is the set
of all $g\in G$ such that $g(x)=y$. Formally,
$G_1 = G\times M$, $G_0=M$, the source map takes $(g,m)$ to $m$
and the target map takes $(g,m)$ to $g(m)$.
We denote the action groupoid by $M/G$.

\hfill

\example
Let $G$ be a Lie group acting on $M$ with compact
stabilizers. Then its action groupoid is proper.
If $G$ is finite, then its action groupoid is 
\'etale and proper.

\hfill

\definition
Let $F: \ (G_0, G_1)\to (H_0, H_1)$ be
a morphism of Lie groupoids (that is, a functor
from one category to another, compatible with
the topology and smooth structure on $G_i, H_i$).
Consider the space ${\cal X}:= (G_0\times G_0)\times_{(F,F)=(s,t)} H_1$
obtained as the set of all triples $(x, y, h)\in G_0\times G_0\times H_1$,
such that $(F(x), F(y))=(s(h), t(h))$.
Let  ${\cal Y} \subset G_0 \times H_1$
be the space of all $(x, h) \in G_0 \times H_1$ such 
that $s(h) = F(x)$.
The functor $F$ is called {\bf an equivalence}
if the map $G_1 \to {\cal X}$ taking $g$ to $(s(g), t(g), F(g))$
is a diffeomorphism, and the map ${\cal Y} \to H_0$ taking $(x, h)$ 
to $s(h)$ is a surjective submersion.

\hfill

\remark
Essentially, the equivalence $(G_0, G_1)\to (H_0, H_1)$  means that the
coarse spaces of these groupoids are identified, and 
this identification is compatible (at least locally)
with the natural group action on the local charts, see 
\ref{_orbifold_charts_Remark_} below.

\hfill

\definition
Let $(G_0, G_1)$ be a Lie groupoid, $s,t:\; G_1\to G_0$ the source 
and target maps, and $U\subset G_0$ an open set. The {\bf restriction groupoid}
is the groupoid whose set of objects is $U$ and set of morphisms
is $s^{-1}(U)\cap t^{-1}(U)$, that is, all arrows connecting points
of $U$ to points of $U$. Clearly, the restriction groupoid is
a Lie groupoid as well.

\hfill

\definition
A groupoid $(G_0, G_1)$ is called {\bf effective} if for general 
$x\in G_0$, the set $\Mor(x,x)$ has only one element (the identity). 
{\bf An orbifold} is an equivalence class of
 proper \'etale effective groupoids. Abusing the language,
we call any proper \'etale effective groupoid ``an orbifold''.

\hfill

This definition
is justified by the following proposition.

\hfill

\proposition\label{_proper_etale_groupoid_Proposition_}
Let $(G_1, G_0)$ be an orbifold, 
that is, a proper \'etale effective groupoid.
Then for any point $x\in G_0$ there exists a neighbourhood
$U$ such that the restriction of $(G_1, G_0)$ to $U$ is 
isomorphic to an action groupoid $U/G$, where $G$ is
a finite group acting on the manifold $U$ by diffeomorphisms, freely
at a general point.

\proof 
This is a special case of \cite[Theorem 2.3]{_Zung:groupoids_};
see also \cite{_Moerdjik_Pronk_}. \endproof

\hfill

\remark\label{_orbifold_charts_Remark_}
In the traditional definition of an orbifold,
the orbifold is a topological space, locally covered
by charts identified with $\R^n/G$, where $G$ is a finite
group, and the information about the $G$-action is kept when
passing to transition maps.
From \ref{_proper_etale_groupoid_Proposition_} it is
clear that these two definitions are equivalent.
Following Lerman, we use groupoids to define the orbifolds;
this approach is pioneered in symplectic and Poisson geometry. 
For us this definition is essential,
because we need to deal with coherent sheaves on complex
varieties equipped with an orbispace
(that is, a proper \'etale effective groupoid) structure.

\hfill

\remark
From \ref{_proper_etale_groupoid_Proposition_},
it is clear how to obtain an orbifold atlas from 
the groupoid. To construct the converse, we use the 
notion of a {\bf cover \'etale groupoid}, defined as follows.

\hfill

\definition
Consider a manifold $M$ with an open cover $\{U_i\}$.
Let $G_0:={\cal U}:= \coprod_i U_i$ be the disjoint union of $U_i$,
and $G_1:= {\cal U}\times_M {\cal U}= \coprod_{i,j}U_i \times U_j$. 
The source and the target maps $s, t: G_1 \to G_0$ are projections
to the first and second component, ${\cal U}\times_M {\cal U}\to {\cal U}$. 
Clearly, $(G_0, G_1)$ is a proper \'etale groupoid
(all sets $\Mor(x, y)$ are either empty or singletons),
and its coarse space is $M$. The groupoid
$(G_0, G_1)$ is called  {\bf a cover \'etale groupoid}
associated with $M$ and its cover $\{U_i\}$.

\hfill

\remark
If we use the same construction with an ``orbifold cover'',
of an orbifold $M$ and enlarge the space of morphisms by adding the group
action, we obtain the proper \'etale groupoid corresponding to $M$.
Equivalent groupoids correspond to equivalent atlases.

\hfill

\remark
At this point we start considering ``orbivarieties'',
the possibly singular complex-analytic version of 
the notion of the orbifold. The Lie groupoids
do not work in this situation, because $G_0$
is manifestly non-smooth. We resolve this by using
the \'etale topological
groupoids instead of \'etale groupoids.

\hfill

\definition
Let $(G_0, G_1)$ be a groupoid
such that $G_0$ and $G_1$ are topological spaces.
We say that $(G_1, G_0)$ is {\bf an \'etale topological
groupoid} if the source and target maps are local
homeomorphisms, and the rest of the structure maps are
continuous. Also, $(G_0, G_1)$ is {\bf proper}
if the map $G_1 \to G_0\times G_0$ sending a morphism
to a pair (source, target) is proper.

\hfill

\example
If $(G_0, G_1)$ is an action groupoid of a group $G$ acting on $M$,
its coarse space is $M/G$. Also, the coarse space of an orbifold
(that is, of an effective proper \'etale groupoid) is the
variety underlying this orbifold. 

\hfill

\definition
Let $(G_0, G_1)$ be 
a proper, effective \'etale topological groupoid. Assume that
every point in $G_0$ has a neighbourhood $U$
such that the restriction groupoid is an action
groupoid of a finite group acting on $U$ by
homeomorphisms. Then $(G_0, G_1)$ is called
{\bf an orbispace.} The {\bf structure of
an orbispace} is the groupoid data attached
to the coarse space of $(G_0, G_1)$.

\hfill

\definition\label{_orbivariety_Definition_}
{\bf A complex orbivariety} is a complex variety
$X$ e\-quip\-ped with the structure of an orbispace
compatible with the complex structure.
More precisely: a structure of a complex
orbivariety on $X$ is a proper, effective \'etale 
topological groupoid $(G_0, G_1)$ such that
$G_0, G_1$ are complex varieties, the 
structure maps are complex analytic, 
the source and target maps are locally
biholomorphic, and the complex variety
underlying $X$ is identified with 
the coarse space of $(G_0, G_1)$.

\hfill

Let us explain this definition. An orbivariety
is a more abstract version of the notion
of an orbifold, adapted to singular varieties.
An orbifold is a topological space
covered by an atlas of open sets homeomorphic to
$\R^n/G$, where $G$ is a finite group, and this
finite group action is carried over to refinements
of the atlas. Similarly, for complex orbivarieties,
the {\bf orbispace cover} of $M$ is an open cover
$(U_i, G_i)$ by open subvarieties $U_i\subset M$ 
obtained as a quotient $U_i=W_i/G_i$
by a finite group $G_i$. We assume that the transition maps are
compatible with the group actions on the intersections
$U_i\cap U_j$. This is the intuitive, though less rigorous,
approach to the orbifolds, and we will use it for
complex orbivarieties as well.

\hfill

\remark
\ref{_orbivariety_Definition_} is not
enough: for a complete picture,
we need to define the equivalence of
groupoid structures, and define
an orbivariety as the equivalence
class of the groupoid structures.
This would parallel the definition
of a smooth manifold as an equivalence
class of smooth atlases. We don't follow
this route to the end, referring the reader
to \cite{_Lerman:Orbifolds_} instead, or to \cite{_Calle:Stacks_}
for the Deligne-Mumford stack flavoured definition.

\hfill

\remark
In algebraic geometry, the orbivarieties are called
{\bf Deligne-Mumford stacks} 
(\cite[Definition 101.4.1]{_Stacks_Project_}; see also
\cite{_Calle:Stacks_}).

\subsection{Coherent sheaves on orbivarieties}

\definition\label{_coherent_on_orbi_Definition_}
Consider an orbivariety $M$,
represented by its orbispace cover $\{U_i\}$,
with the local groups $G_i$ acting on each $U_i$.
A coherent sheaf on $M$
(see also \cite[Definition 9.15]{_OV_book_},
\cite[Definition 4.2.7]{_Boyer_Galicki_})
is a collection of coherent sheaves over each open 
chart $U_i$, equivariant under $G_i$,
equipped with $G_i$-invariant gluing
maps satisfying the cocycle condition.

\hfill

\remark
It is possible to define the coherent sheaves
over orbivarieties in the language of groupoids
(Subsection \ref{_orbispaces_Subsection_}).
Consider the groupoid $(G_0, G_1)$ associated with the
orbivariety $M$; here $G_0$ is the disconnected union of all
the charts in a given orbifold atlas $\{U_i\}$, and $G_1$ the
graphs of the transition maps, defined on appropriate
open subsets of $U_i$, with $s, t:\; G_1 \to G_0$
the source and the target maps. We can define
a coherent sheaf on $M$ as a coherent sheaf
${\cal F}$ on $G_0$, equipped with an isomorphism
$\psi:\; s^*{\cal F}\tilde\to t^*{\cal F}$, satisfying the
usual equivariancy axioms, which can be formulated
in the same way as in \ref{_equivariant_sheaf_Definition_}.

\hfill

\remark
From the example described in Subsection 
\ref{_orbifold_sheaf_example_Subsection_},
it is apparent that the category of coherent sheaves
on an orbivariety is not, generally speaking,
equivalent to the category of coherent sheaves
on its quotient variety (that is, the coarse space
of the corresponding groupoid). This distinction
becomes important even when we consider the
line bundles. In the literature, the line bundle
over an orbifold (defined as in \ref{_coherent_on_orbi_Definition_})
is often called ``an orbifold bundle'', to distinguish
from vector bundles on the quotient space.
Orbifolds were first defined by I. Satake, who
called them ``V-manifolds''. In his earliest papers
\cite{_Baily:decomposition_,_Baily:embedding_},
W. L. Baily has defined what he called ``V-bundles'' and
proved the Kodaira embedding theorem for V-bundles
of rank 1 with positive curvature; his definition
of a V-bundle is the same as our definition of
an orbifold bundle.
In the present paper, we will always consider
coherent sheaves and vector bundles in the
sense of \ref{_coherent_on_orbi_Definition_}. 
When the distinction between the orbifold sheaves
and the coherent sheaf on the underlying space becomes
important, we will occasionally call the
coherent sheaves (defined as in 
\ref{_coherent_on_orbi_Definition_})
{\bf the orbisheaves,} to distinguish
them from the coherent sheaves on the
corresponding coarse space.

\hfill

\example
Let $Z$ be a complex variety equipped with
an action of a finite group $\Gamma$, and
$B$ a $\Gamma$-equivariant vector bundle on $Z$.
Then $Z/\Gamma$ is an orbivariety,
and $\Tot(B)/\Gamma$ is the total space
of a holomorphic orbifold bundle over the orbivariety
$Z/\Gamma$.

\subsection{A $\C^*$-action on a cone variety defines an
  orbispace structure on its quotient}
\label{_C^*_cone_Subsection_}

For the rest of this section,
we consider a cone variety $(X, \gamma)$.
By \ref{_cone_str_on_open_cone_variety_Theorem_}, 
there is a $\C^*$-action on $X$ containing 
a contraction $\gamma_0$. 
This $\C^*$-action is obtained
by embedding $(X, \gamma)$ to $(\C^n,A)$,
where $A\in \GL(n, \C)$ is a linear contraction,
and applying the Jordan-Chevalley decomposition;
we can define $\gamma_0$ as an element of the 
Zariski closure of the semisimple part of $A$,
close enough to the semisimple part of $A$, such that the Zariski closure of
$\langle \gamma_0\rangle$ is isomorphic to $\C^*$.

\hfill

\theorem\label{_C^*_action_orbispace_Theorem_}
Let $(X, x, \gamma_0)$ be a closed cone variety such that
the Zariski closure of $\gamma_0$ is isomorphic
to $\C^*$, and $X_0:= X\backslash x$ the
corresponding open cone variety. 
In this situation we denote the 
Zariski closure of $\gamma_0$ by $\C^*_Z$. Then 
the quotient $Z:=X_0/\C^*_Z$ is is equipped
with a natural structure of a projective
orbivariety and an ample orbifold line bundle $L$ on $Z$, such that
$X_0$ is biholomorphic to $\Tot^0(L)$.

\hfill

\pstep
By \ref{_cone_embedded_Theorem_},
$X_0$ admits a holomorphic embedding to $V:=\C^n$
in such a way that $\gamma_0$ is the restriction
of a contraction $A\in \End(V)$,  and the
Zariski closure of 
$\langle A \rangle \subset \GL(V)$ is $\C^*_Z\cong \C^*$. By 
\ref{_C^*_on_cone_embedding_Theorem_},
this embedding can be chosen $\C^*_Z$-equivariant.
Then the orbits of the $\C^*_Z$-action on $X_0$
correspond to the orbits of a $\C^*$-action
on $V\backslash 0$, which gives a holomorphic
embedding from $Z=X_0/\C^*_Z$ to $\frac{\C^n\backslash 0}{\C^*}$,
identified with a weighted projective space.

\hfill

{\bf Step 2:} 
Consider the weighted projected space 
$W\CP^{n-1}:=\frac{\C^n\backslash 0}{\C^*}$ 
as an orbifold, and induce the orbifold structure
on $X_0/\C^*_Z$ from the embedding
$X_0/\C^*_Z\hookrightarrow \frac{\C^n\backslash 0}{\C^*}$.
Since the total space of the orbifold line bundle
on $W\CP^{n-1}$ is $\C^n\backslash 0$, this construction
identifies the space $X_0\subset \C^n\backslash 0$
with the total space of an orbifold line bundle
over the orbifold quotient $X_0/\C^*_Z\subset W\CP^{n-1}$.
\endproof

\hfill

An  explicit
description of the orbivariety structure on $X_0/\C^*_Z$
 is given in \ref{_local_orbi_Proposition_}.

\hfill

\claim\label{_F_max_subgroup_Claim_}
Let $X_0$ be a cone variety equipped with a
$\C^*_Z$-action as in \ref{_C^*_action_orbispace_Theorem_}.
Then there exists a finite subgroup $F_{\max} \subset \C^*_Z$
such that the natural $\C^*$-action on $X_0/F_{\max}$
is free.

\hfill

\proof
It is possible to prove \ref{_F_max_subgroup_Claim_}
by embedding $X_0$ to $\C^n\backslash 0$ and then proving 
\ref{_F_max_subgroup_Claim_} in the case when
the orbifold under question is $W \C P^{n-1}$.
Instead, we give a direct proof, which might be
instructive as well.

Given $x\in X_0$, let $F(x)$ be the stabilizer of
$x$ in $\C^*_Z$. Clearly, the set of all $x\in X$ such that
$F(x) \ni \zeta$ is Zariski closed for any $\zeta\in \C^*_Z$.
Recall that any descending chain of subvarieties in
a compact complex variety terminates; this property
is called ``Noetherianity''. The variety $X_0$
is not compact, but the space of $\C^*_Z$-orbits is
compact, hence the Noetherianty is also true
for $\C^*$-invariant subvarieties: any descending
chain of $\C^*$-invariant subvarieties terminates.
We can apply this to the varieties 
$X_{F}:=\{x\in X_0 \ \ | \ F(x)\supset F\}$, where
$F\subset \C^*_Z$ is a finite subgroup. Take
a descending chain of subvarieties
$X_0 \supset X_{F_1} \supset ... \supset X_{F_n}$,
where $F_1 \subset F_2\subset ... \subset F_n \subset \C^*_Z$.
Since it terminates, $F(x)$ is maximal on
a finite collection of subvarieties.
This implies that the set 
$\bigcup_{x\in X_0} F(x)$  generates a finite subgroup 
$F_{\max}\subset \C^*_Z$.
\endproof

\hfill

In the assumptions of \ref{_F_max_subgroup_Claim_},
let $X_0':= X_0/F_{\max}$, considered as an orbivariety.
The group $\C^*/F_{\max}=\C^*$ acts freely
on $X_0'$, hence $Z=X_0'/\C^*$ 
inherits this orbifold structure.

\hfill

To make this construction even more explicit, we
would like to exchange the order of taking quotients: first we take
a quotient by the free $\C^*$-action, obtaining a variety,
and then the finite group quotient, which gives
an action groupoid orbifold. This is possible to
do only locally. Indeed, even the simplest example
of our construction, the weighted projective space,
is not isomorphic to an action groupoid orbifold.

To reverse the order of taking quotients in 
the statement of \ref{_F_max_subgroup_Claim_},
we identify the cyclic group\footnote{Every finite
subgroup of $\C^*$ is a cyclic group.} $F_{\max}$
with $\Z/r\Z$, where $r=|F_{\max}|$. 

\hfill

\remark
In \ref{_cone_embedded_Theorem_}, we approximate
the $\Z$-action $\gamma$ on a cone variety by $\gamma_0$
which can be extended to an action of $\C^*$.
Any action of $\C^*$ which contains a contraction
has finite stabilizers, contained in $S^1\subset \C^*$.
Otherwise this action restricted to the corresponding 
orbit $\C^*\cdot x$ would factorize through
a compact group $\C^*/F(x)$, where $F(x)\subset \C^*$ 
is its stabilizer, and a compact group cannot
contain contractions.

\hfill

\remark
Let $\pi:\; X_0 \to X_0/\C_Z$ be the fibration constructed earlier.
For proper complex analytic fibrations, it is known that
a neighbourhood of a central fiber retracts to the central
fiber (\cite{_Persson:degene_,_Clemens:degene_}).
Below we show that this statement is  also true for proper $\C^*$-fibrations.

\hfill

\claim\label{_pi_1_locally_Z_Claim_}
Let $X_0$ be a cone variety equipped with a locally free
and proper $\C^*_Z$-action, as in \ref{_C^*_action_orbispace_Theorem_},
and $U\subset X_0$ be a $\C^*_Z$-invariant open subset.
Let $A$ be an orbit of the $\C^*_Z$-action.
Then the natural map $A\to U$
is a homotopy equivalence, when $U$ is a sufficiently small
$\C^*_Z$-invariant neighbourhood of $A$.

\hfill

\pstep
Consider a point $a\in A$, and let
$U_a$ be a small neighbourhood of $a$.
Shrinking $U_a$ if necessary, we
may assume that each $\C^*$-orbit
intersects $U_a$ in a contractible
open set. The orbits of $\C^*$ define a foliation
on $U_a$, and the projection ${\cal V}\: :\: U_a\to W_a$
to the corresponding leaf space is a locally
trivial dimension 1 fibration. To trivialize
this fibration means to choose a complex analytic
section $\Delta$ which intersects each leaf
of ${\cal V}$ in only one point. Consider the
map $\Psi:\; \C^*\times \Delta \to U$
taking a point to its image under the $\C^*$-action.
Shrinking $U$ if necessary, we can assume that
$\Psi$ is a finite covering. Replacing $\Delta$
by an open subvariety if necessary, and shrinking $U$ again, 
we can assume that $\Delta$ is contractible 
(any complex variety is locally contractible),
and the map $\Psi:\; \C^*\times \Delta \to U$
is still a covering.

\hfill

{\bf Step 2:} 
Let $z\in \Delta$ be the point which 
satisfies $\Psi(\C^* \times \{z\})=A$.
The space $U$ is homotopy equivalent to $\C^*$,
because its covering is $\Delta \times \C^*$,
and $A$ is also indentified with $\C^*$. We need to show that
the corresponding map $\pi_1(A) \to \pi_1(U)$
is an isomorphism. 

Let $\phi_t:\; \Delta\times [0,1] \to \Delta$ be the 
continuous contraction of $\Delta$ to $z$.
Consider the homotopy
map $h_t:\; U\times [0,1] \to U$ contracting $U$ to the central fiber $A$,
taking a point $\Psi(x \times \lambda)\in U$ to 
$\Psi(\phi_t(x),lambda)$,
where $x\in \Delta$ and $\lambda\in \C^*$.
This map is well defined because
$\Psi$ is a finite covering, compatible
with the $\C^*$-action, hence $h_t(\Psi(x, t))$
is independent from the choice of the preimage
$(x, \lambda)\in\Delta\times \C^*$.
This map is a deformational retract 
of $U$ to $A$, which concludes the proof.
\endproof

\hfill

We are in position to give an explicit local description
of the orbivariety structure on $X_0/\C^*_Z$ 
(\ref{_C^*_action_orbispace_Theorem_}).

\hfill

\proposition\label{_local_orbi_Proposition_}
Let $X_0$ be a cone variety equipped with a
$\C^*_Z$-action as in \ref{_C^*_action_orbispace_Theorem_}, 
and $A\subset X_0$ an orbit of $\C^*_Z$.
Consider a $\C^*_Z$-invariant neighbourhood $U\subset X_0$ of $A$
such that $\pi_1(A)=\pi_1(U)$ (\ref{_pi_1_locally_Z_Claim_}).
Let $F_A\subset \C^*_Z$ be the stabilizer of the $\C^*_Z$-action on $A$,
considered as a cyclic group of order $d$,
and $\tilde U$ a $d$-sheeted cover of $U$ 
associated with the group $d\Z\subset \Z=\pi_1(U)$.
Then the action of $\C^*$ is lifted to the free
action on $\tilde U$, and the quotient orbifold
$U/\C^*$ is identified with the action groupoid
(\ref{_action_groupoid_Example_})
of the group $F_A$ acting on the variety 
$\tilde U/\C^*$.\footnote{Let $\Gamma$ be a finite 
group acting on a variety $X$. The corresponding action groupoid
is just $X/\Gamma$ with the induced orbivariety structure.}

\hfill

\proof
Clearly, $\C^*/F_A$ acts freely on the 
quotient $U/F_A$. Passing to the order $d$ covering $\tilde U$,
we obtain that $\C^*$ acts freely on $\tilde U$.
The identification $\tilde U/\C^*= \frac{U/F_A}{\C^*/F_A}$
follows directly from the construction.
\endproof

\hfill


\remark
The orbispace $Z:=X_0/\C^*$ is projective, 
because the line bundle $L$ constructed in 
\ref{_cone_str_on_open_cone_variety_Theorem_}
is ample. A similar result is true for the closed
cone varieties: they can be obtained as
affine cones over the polarized orbivarieties.
We give the necessary definitions and
prove this result in the next subsection.

\subsection{Polarized orbifolds and orbivarieties}

\definition
Let ${\cal F}$ be a coherent orbisheaf over an orbifold.
This is the same as a collection of $G_i$-equivariant 
coherent sheaves over each orbifold chart $\{U_i\cong B_i/G_i\}$,
where $B_i$ is an open ball and $G_i$ a finite group,
compatible with the transition maps.
Such a sheaf induces a coherent sheaf ${\cal F}^{inv}$ 
of $G_i$-invariant sections,
called {\bf the coherent sheaf underlying the orbisheaf ${\cal F}$}.

\hfill

\remark
Let $L$ be an orbibundle of rank 1 over an orbispace $X$.
The underlying coherent sheaf $L^{inv}$
also has rank 1, but it is not necessarily locally free.
Over each point $x\in U_i\cong B_i/G_i$ which has
a non-trivial stabilizer $\St_x(G_i)$,
any section of $L^{inv}$ would vanish.
However, the action of $\St_x(G_i)$ is trivial
on the tensor power $L^{\otimes k}$ of order $k:=|\St_x(G_i)|$.
This implies that $(L^{\otimes u})^{inv}=L^{\otimes u}$
for any $u$ which is divisible by $|\St_x(G_i)|$ for all $x\in M$.

\hfill

In \cite{_Baily:embedding_}, W. L. Bailey has proved
a version of the Kodaira embedding theorem for orbifolds.
He has shown that a line bundle $L$ with positive curvature
over a compact orbifold is ample, that is, induces a complex
analytic embedding to $\C P^n$. This embedding does not
pay any attention to the orbifold structure; indeed,
replacing $L$ by a tensor power, we can assume 
that the orbibundle structure on $L$ is trivial
(that is, $L$ is isomorphic to its underlying sheaf 
$L^{inv}$).

In \cite[Proposition 2.11]{_Ross_Thomas:Orbifolds_}, J. Ross and R. Thomas
have proved a Kodaira-type embedding theorem
which takes care of the orbifold singularities.
The same theorem was proven by D. Abramovich and
B. Hassett using the language of stacks and
algebraic spaces (\cite[Proposition 2.4.2]{_Abramovich_Hassett_}).
We are going to explain their results now.

An orbibundle of rank 1 on an orbifold $X$ is called
{\bf a polarization} if it 
admits a Chern connection with positive curvature
and the action of $\St_x(G_i)$ on the fiber $L\restrict x$
is faithful for all $x\in X$. In 
\cite{_Ross_Thomas:Orbifolds_},
Ross and Thomas have shown 
that the graded ring $R:=\bigoplus_i H^0(X, L^{\otimes i})$
satisfies $\Proj R= W\C P^n$, that is, its Proj space
is a weighted projectiive space, and the embedding
$X \hookrightarrow W\C P^n$ is compatible with the 
orbifold structure.

Note that the existence of a polarization is a restrictive
condition; for example, if $\St_x(G_i)$ is non-abelian
for some $x\in X$, the orbifold $X$ cannot admit a
polarization. In \cite{_Abramovich_Hassett_}, Abramovich and Hasset call
an orbispace with abelian $\St_x(G_i)$ {\bf cyclotomic stack}.

\hfill

An orbivariety version of Ross-Thomas theorem is
proven in the same way as 
\cite[Proposition 2.11]{_Ross_Thomas:Orbifolds_}. 
We define an ample orbibundle
over a compact orbivariety $Z$ as an orbibundle $L$
of rank 1 such that a sufficiently big tensor power of $L$
is globally generated on $Z$ and induces a projective
embedding of the underlying variety to $\C P^n$.

\hfill

\theorem
Let $Z$ be an orbivariety, and
$L$ an ample orbibundle such that
the action of $\St_x(G_i)$ on the fiber $L\restrict x$
is faithful for all $x\in X$. Then $R:=\bigoplus_i H^0(Z, L^{\otimes i})$
satisfies $\Proj R= W\C P^n$, and the corresponding embedding
is compatible with the orbivariety structure.

\hfill

\proof
The embedding $Z\hookrightarrow {\Bbb P}H^0(Z, L^{\otimes i})^*$ is 
given by the definition of an ample bundle. Therefore,
the natural map $Z \to \Proj R$ is a holomorphic embedding.
The argument which proves that $\Proj R$ is a weighted
projective space (\cite[Proposition 2.11]{_Ross_Thomas:Orbifolds_}) 
works for orbivarieties as 
well as for orbifolds, 
see also \cite[Section 2.4]{_Abramovich_Hassett_}.
\endproof

\hfill

By definition, a weighted projective
space is a quotient of $\C^n\backslash 0$ by an
action of $\C^*$, which has weights
$w_1, ..., w_n \in \Z^{>0}$.
Therefore, $\C^{n+1}$ is $\C^*$-fibered
over $W \C P^n$.
Let $(Z, L)$ be a polarized orbivariety,
and $\phi:\; Z \to W \C P^n$ the embedding
constructed above.
Define {\bf the affine cone}  $C(Z,L)$ of $(Z,L)$ 
as the closure of the preimage of $Z$ in $\C^{n+1}$.
By Remmert-Stein theorem, $C(Z,L)$ 
is a $\C^*$-invariant complex analytic subvariety
in $\C^{n+1}$.

\hfill

Let $(X, \gamma)$ be a closed cone variety, 
$X_0$ the corresponding open cone variety, 
and $\C^*_Z\subset \Aut(X)$ the $\C^*$-action
constructed in \ref{_C^*_action_orbispace_Theorem_}.
Then $Z:=X_0/\C^*_Z$ is a projective orbivariety.
Denote by $L$ the corresponding ample bundle on $Z$,
such that $\Tot^\circ(L)=X_0$.

\hfill

\theorem\label{_closed_cone_alg_cone_Theorem_}
In these assumptions, $X$ is naturally
isomorphic to the affine cone $C(Z,L)$.

\hfill

\proof
By construction, the orbifold structure on $Z$ is induced
from the weighted projective space
by the way of $\C^*$-equivariant
embedding 
$\Tot^\circ(L)\hookrightarrow \C^{n+1}\backslash 0$. 
This implies that $L$ defines a polarization on the
orbivariety $Z$.  The closed cone $X$ is identified
with the closure of $\Tot^\circ(L)$ in $\C^{n+1}$
by construction.
\endproof

\subsection{$\C^*_Z$-equivariant sheaves on cone varieties and orbisheaves}

\theorem\label{_Roberts_invariant_sections_Theorem_}
Let $G$ be a reductive complex Lie group
acting with finite stabilizers on 
a complex analytic variety $X$, and ${\cal F}$
a $G$-equivariant coherent sheaf.
Assume that the orbit space is 
a complex analytic variety, all orbits are closed,
and the projection to the space of the
orbits $\pi:\; X \to Z$ is a morphism of complex 
varieties. We assume that every $z\in Z$
has a neighbourhood such that its preimage
in $X_0$ is Stein. Let $\pi_*({\cal F}^G)$
be the pushforward of the sheaf of $G$-invariant sections
of ${\cal F}$ to $Z$. Then $\pi_*({\cal F}^G)$ is 
a coherent sheaf on $Z$.

\proof \cite[Theorem 3.1]{_Roberts:G_sheaves_}.
\endproof

\hfill

\remark
Using  \cite[Theorem 3.12]{_Kollar:quotient_}
instead of \cite[Theorem 3.1]{_Roberts:G_sheaves_},
the assumptions of \ref{_Roberts_invariant_sections_Theorem_}
can be relaxed a lot:
the group $G$ can be assumed algebraic,
instead of reductive, and the Stein
assumption is not needed. However,
\cite[Theorem 3.12]{_Kollar:quotient_}
is proven only in the algebraic category.

\hfill

\remark
We are going to use \ref{_Roberts_invariant_sections_Theorem_}
only when $G=\C^*$.

\hfill

\theorem\label{_C^*_equiv_sheaves_and_orbispace_Theorem_}
Let $X_0$ be a complex analytic variety equipped with an 
action of $\C^*$ with finite stabilizers, free in a general point,
and $Z$ the corresponding projective orbivariety
(\ref{_C^*_action_orbispace_Theorem_}).
Then the category of $\C^*$-equivariant coherent 
sheaves on $X_0$ is equivalent to the
category of coherent (orbi)-sheaves on the 
orbivariety $Z$.

\hfill

\pstep
Let ${\cal F}$ be a $\C^*$-equivariant coherent sheaf of $X_0$.
By \ref{_Roberts_invariant_sections_Theorem_},
the sheaf $\pi_*({\cal F}^{\C^*})$ of ${\C^*}$-invariant sections on 
the coarse space of $Z$ is always
coherent. It is also clear that for any coherent sheaf
${\cal H}$ on the coarse space of $Z$,\footnote{As usual, 
``the coarse space'' of an orbivariety denotes the
underlying complex analytic variety.}
the natural morphism ${\cal H} \to \pi_*((\pi^*({\cal H})^{\C^*})$
is an isomorphism, hence the functor $\pi^*$ 
admits a right inverse.
 However, the natural
morphism $\pi_*(\pi^*({\cal F})^{\C^*}) \to {\cal F}$
is not an isomorphism (Subsection 
\ref{_orbifold_sheaf_example_Subsection_}), 
and to produce the equivalence of 
categories, we need to use the 
orbisheaf structure on ${\cal F}^{\C^*}$.

\hfill

{\bf Step 2:} When the action
of $G$ is free, the natural monomorphism
$\Psi:\; \pi^*(\pi_*({\cal F})^G) \to {\cal F}$ is an isomorhism.
To see this, consider a subvariety $D$ transversal
to the orbits of $G$. Then, locally in a 
$G$-invariant neighbourhood $U$ of an orbit, $U=G\times D$,
and any section of ${\cal F}\restrict D$ extends
uniquely to a $G$-invariant section of ${\cal F}\restrict U$,
hence $\Psi$ is surjective.

\hfill

{\bf Step 3:}
We are going to produce the functor ${\cal F}\mapsto \pi_*({\cal F}^{\C^*})$ 
from sheaves on 
$X_0$ to the orbisheaves on the orbivariety $X_0/\C^*$
and ${\cal H} \mapsto \pi^*({\cal H})$
from orbisheaves on $X_0/\C^*$ to $X_0$, and show
that they are inverse to each other. 
We use the local description of the
orbivariety structure on $X_0/\C^*$
given in \ref{_local_orbi_Proposition_}.

Let $A\subset X_0$ be an orbit of $\C^*$, and
 $U\subset X_0$ a $\C^*$-invariant
neighbourhood of $A$ satisfying the assumptions of 
\ref{_local_orbi_Proposition_}. As in \ref{_local_orbi_Proposition_},
$F_A$ denotes the stabilizer of the $\C^*$-action on $A$,
$F_A=\Z/d$. Then $\C^*/F_A$ acts freely on $U/F_A$.
In \ref{_local_orbi_Proposition_} we constructed
a $d$-sheeted cover $\tilde U\stackrel \sigma \to U$
such that the pullback of the $\C^*$-action is free on $\tilde U$,
and $\tilde U/\C^*= \frac{U/F_A}{\C^*/F_A}$.

Let ${\cal F}$ be a $\C^*$-equivariant coherent sheaf
on $X_0$. Pulling back ${\cal F}$ to $\tilde U\stackrel \sigma \to U$,
we obtain a sheaf $\sigma^*({\cal F})$ which is equivariant with respect
to the free action of $\C^*$.
By \ref{_Roberts_invariant_sections_Theorem_},
the sheaf of $\C^*$-invariant sections of 
$\sigma^*({\cal F})$ is coherent on
$\tilde U/\C^*$; it is $F_A$-equivariant by 
construction. This defines an equivalence
of categories between $\C^*$-invariant 
coherent sheaves on $U$ and $F_A$-equivariant
sheaves on $\tilde U/\C^*$, which is the 
same as orbisheaves on the orbivariety $U/\C^*$.

\hfill

{\bf Step 4:} In Step 3, we produced a local equivalence
of categories; to make it global, we notice that the
pullback of an orbifold sheaf to $X_0$ is by construction
$\C^*$-equivariant. This defines a functor from
the category of orbisheaves on $Z$ to the category
of $\C^*$-equivariant sheaves on $X_0$.
Since this functor is an eqiuvalence on any 
open subset of $X_0$, and the definition of a sheaf
is local, it is an equivalence.
\endproof


\section{Equivariant reflexive sheaves over cone varieties}


\subsection{Siu's extension theorems}

We are going to prove the extension of coherent sheaves
result that we use later in this paper, using the 
classical theorems of Y.-T. Siu and J.-P. Serre.

\hfill

\theorem\label{_Siu_extending_Theorem_}
(Siu, \cite{_Siu:Extending_})\\
Let $S\subset X$ be a subvariety of a complex space $X$, $\dim S \leq n$,
and ${\cal F}$ a coherent analytic sheaf defined on $X\backslash S$.
Assume that for any $R\subset X\backslash S$ such that $\dim R\leq n+1$, any open 
subset $U\subset X\backslash S$, and any 
 $f\in H^0(U\backslash R,{\cal F})$, the section
$f$ can be extended to a section of ${\cal F}$ over $U$. 
Then ${\cal F}$ can be extended to a coherent sheaf on $X$.

\proof \cite[Theorem 1]{_Siu:Extending_}. \endproof

\hfill

\definition
A coherent sheaf ${\cal F}$ on a variety $X$
is called {\bf normal} if for any complex analytic
subset $S\subset X$, $\codim S \geq 2$, and any open $U\subset X$,
the restriction map $H^0(U, {\cal F}) \to H^0(U\backslash S, {\cal F})$
is an isomorphism.

\hfill

\remark
An algebraic variety $X$ is normal if and only if it is non-singular
in codimension 1 and its structure
sheaf $\calo_X$ is normal (\cite[Theorem 5.8.6]{_EGA:IV_2_}); 
this explains the term.

\hfill

\theorem\label{_normal_is_reflexive_Theorem_}
Let $X$ be a normal complex analytic variety, and
${\cal F}$ a coherent sheaf over $X$.
Then ${\cal F}$ is reflexive if and only if
it is normal.

\hfill

\proof
This result is classical and well known,
but we did not find a good reference for
the situation when $X$ is a complex analytic variety.
When $X$ is a complex manifold, it is 
\ref{_normal_is_reflexive_Theorem_} is 
\cite[Lemma 1.1.2]{_OSS_} or 
\cite[Proposition 1.30]{_Brinzanescu:LMN_1624_}.
When $X$ is a normal integral scheme, 
\ref{_normal_is_reflexive_Theorem_} is
\cite[Proposition 1.6]{_Hartshorne:stable_reflexive_}.

Let $S\subset X$ be a subvariety, 
and $j:\; X\backslash S\to X$ the standard embedding.
Consider the natural sheaf morphism ${\cal F} \to j_* j^* {\cal F}$.
Clearly, ${\cal F}$ is normal if and only if this map
is an isomorphism whenever $\codim S \geq 2.$
As explained in the paragraph after 
\cite[Proposition 7]{_Serre:Prolongement_},
$j_* j^* {\cal F}={\cal F}$ for any reflexive sheaf.
Therefore, reflexive sheaves are normal.
Conversely, let ${\cal F}$ be a normal
sheaf. Then ${\cal F}$ is non-singular
in codimension 1 (\cite[Cor.  Lem. 1.1.8]{_OSS_}). Let $S\subset X$ be its
singular locus, and $j:\; X\backslash S\to X$ the
standard embedding. Then ${\cal F}=j_*j^*{\cal F}$ by normality.
However, the sheaf $j^*{\cal F}$ is locally free, hence
reflexive, and a pushforward of a reflexive sheaf
is reflexive by \cite[Proposition 7]{_Serre:Prolongement_}.
\endproof

\hfill

Together with \ref{_normal_is_reflexive_Theorem_},
Siu's extension theorem (\ref{_Siu_extending_Theorem_})
immediately implies the following.

\hfill

\corollary\label{_Siu_over_point_Corollary_}
Let $X$ be a normal variety, $\dim X \geq 3$,
and $X_0:= X \backslash x$ the complement to a
point $x\in X$. Denote the natural embedding $X_0\hookrightarrow X$
by $j$. Then for any reflexive sheaf ${\cal F}$ on $X_0$,
its pushforward $j_* {\cal F}$ is coherent on $X$.
\endproof

%
%
%
%
%
%
%
%
%
%
%
%
%
%

\subsection{Equivariant line bundles on cone varieties}

Let $(X_0, \gamma)$ be an open cone variety.
Consider the algebraic structure on $(X_0, \gamma)$ constructed in 
\ref{_algebraic_str_on_cone_uniqueness_Theorem_}.
By construction, $\gamma$ is contained in
an algebraic group acting on $X_0$; let
${\cal G}$ be the Zariski closure of $\langle \gamma\rangle$
in this algebraic group. As explained in the proof of
\ref{_cone_str_on_open_cone_variety_Theorem_},
${\cal G}$ contains a subgroup $\C^*_Z$ isomorphic to $\C^*$
acting on $X_0$ with finite stabilizers, and the quotient
$Z:=X_0/\C^*_Z$ is a projective variety. Moreover, $X_0$ is isomorphic
to the space of non-zero vectors in an ample
line bundle on $Z$. Throughout this
section, we keep these assumptions
and this notation.

\hfill

\remark\label{_gamma_acts_on_Z_Remark_}
By construction, $\gamma$ commutes with $\C^*_Z$ in ${\cal G}$,
hence $\gamma$ induces a $\Z$-action on $Z= X_0/\C^*_Z$.
We denote this action by $\langle \gamma\rangle$ as well.

\hfill

We are interested in the category of coherent  sheaves
on the quotient
$X_0/\langle \gamma\rangle$, where $\gamma$ is a contraction.
Clearly, this category is equivalent
to the category of $\langle \gamma\rangle$-equivariant
analytic coherent sheaves on $X_0$.

\hfill

\proposition\label{_categories_on_cone_Proposition_}
Let $(X_0, \gamma)$ be an open cone variety, $\dim_\C X_0> 2$,
considered as a quasiprojective variety 
(\ref{_algebraic_str_on_cone_uniqueness_Theorem_}),
with the induced $\langle \gamma\rangle$-action.
Then the following three categories are naturally
equivalent: 
\begin{description}
\item[(i)] the category of $\langle \gamma\rangle$-equivariant
reflexive analytic coherent sheaves on $X_0$,
\item[(ii)] 
the category of  $\langle \gamma\rangle$-equivariant
reflexive {\em algebraic} coherent sheaves on $X_0$
\item[(iii)]
the category of reflexive coherent sheaves on the complex analytic variety
$X_0/\langle \gamma\rangle$.
\end{description}

\proof
The equivalence of (i) and (iii) is clear.
The category of reflexive analytic sheaves on $X_0$ 
is equivalent to the one on its (normal) Stein completion $X$
by \ref{_Siu_over_point_Corollary_}.
By \ref{_closed_cone_from_open_Theorem_}, $X$ is a closed cone variety.
The algebraic structure on  $\langle \gamma\rangle$-equivariant
coherent sheaves on $X$ is determined by the
$\langle \gamma\rangle$-equivariant structure, as follows from
\ref{_GAGA_reflexive_Theorem_}. By \ref{_GAGA_reflexive_Theorem_},
the algebraic sections are $\gamma$-finite, and
all $\gamma$-finite sections are algebraic. This implies 
the equivalence between (i) and (ii). 
\endproof

\hfill

\proposition\label{_G_F_equivariant_on_F_Proposition_}
Let ${\cal F}$ be a $\gamma$-equivariant coherent sheaf
on a cone variety $(X_0, \gamma)$, and
${\cal G}$ the Zariski closure of $\langle \gamma\rangle$ in
$\Aut(X_0)$, obtained as in \ref{_cone_str_on_open_cone_variety_Theorem_}.
Then there exists a complex algebraic Lie group ${\cal G}_{\cal F}$
equipped with a surjective map to ${\cal G}$ (and hence,
acting on $X_0$) such that the equivariant action of
$\gamma$ on ${\cal F}$ is extended to a
${\cal G}_{\cal F}$-action. 

\hfill

\proof
Let $(X, x, \gamma)$ be the closed cone associated wih $X_0$.
Since $X$ is an affine variety, $X=\Spec(R)$, where $R$ 
is the ring of algebraic functions on $X$
(\ref{_algebraic_str_on_cone_uniqueness_Theorem_}), and coherent
sheaves on $X$ are the same as finitely generated $R$-modules.
Let $W$ be the $R$-module of algebraic sections of  ${\cal F}$.
By construction, $R$ is the ring of $\gamma$-finite 
holomorphic functions on $X$, and $W$ the space
of $\gamma$-finite sections of ${\cal F}$
(\ref{_GAGA_reflexive_Theorem_}).
Let $\check W\subset W$ be a finite-dimensional
$\gamma$-finite space generating $W$ over $R$.
Denote by ${\cal G}_{\cal F}$ the Zariski closure of
$\langle\gamma\rangle\subset \GL(\check W)\times {\cal G}$.
The group ${\cal G}_{\cal F}$ is abelian,
because it is the closure of a $\Z$-action,
and projects to ${\cal G}$ surjectively
by construction. Also, the equivariant action
of $\gamma$ factorizes through ${\cal G}_{\cal F}$.
\endproof

\hfill

\remark
The group ${\cal G}_{\cal F}$
is constructed using the space of generators
$\check W\subset W$. However, different choices
of $\check W\subset W$ lead to the same group,
because ${\cal G}_{\cal F}$ preserves any finite-dimensional
$\gamma$-invariant space in $W$. Therefore,
${\cal G}_{\cal F}$ is defined canonically by ${\cal F}$.

\hfill

It would be optimal 
to use \ref{_C^*_equiv_sheaves_and_orbispace_Theorem_}
to establish an equivalence between
$\gamma$-equivariant reflexive sheaves on $X_0$
and $\gamma$-equivariant (orbi)-sheaves on $Z$.
Unfortunately, it does not work, because
the subgroup $\C^*_Z\subset {\cal G}$ does not
necessarily lift to $\C^*\subset {\cal G}_{\cal F}$.
This happens even in the most elementary 
case: for the line bundles on a Hopf manifold.
In the next subsection, we are going 
to compute ${\cal G}_{\cal F}$ for a line bundle ${\cal F}$
over Hopf manifold $H$, to illustrate the phenomena which are
abound in the general case.

\subsection{Line bundles on Hopf manifolds}
\label{_line_on_Hopf_Subsection_}

\remark
For a classical Hopf manifold, obtained as 
$\frac{\C^n\backslash 0}{\langle A\rangle}$,
where $A=\const\cdot \Id$, we have ${\cal G}=\C^*=\C^* \! \cdot\!\Id$,
which is clear, because $\C^* \! \cdot\!\Id$
is the smallest algebraic group containing $A$.
The same happens for elliptic Hopf manifolds;
for all other Hopf manifolds, $\dim_\C {\cal G} >1$.

\hfill

\remark\label{_Hopf_Picard_Remark_}
The Picard group of a Hopf manifold $H$ is not hard to
compute using the exponential exact sequence 
(\cite[Exercise 33.20]{_OV_book_}):
\[
0\to H^1(H, \Z) \to H^1(H, \calo_H)\to \Pic(H) \to H^2(H, \Z)=0.
\]
It is easy to see that $H^1(H, \calo_H)=\C$ 
(\cite[Theorem 27.9]{_OV_book_}),
which gives $\Pic(H) = \C^*$.
Each line bundle on a Hopf manifold is lifted to 
a line bundle on $\C^n\backslash 0$, which can
be extended to $\C^n$ by \ref{_Siu_over_point_Corollary_};
a rank 1 reflexive sheaf is a line bundle (\cite{_OSS_}),
and a line bundle on $\C^n$ is trivial, which can be
seen from the same exponential exact sequence.
Different $\gamma$-equivariant structures on a line bundle
on $\C^n\backslash 0$ differ by a $\gamma$-invariant automorphism
of this bundle, which has to be constant because
$\frac{\C^n\backslash 0}{\langle \gamma\rangle}$ is compact.
This implies that $\gamma$-equivariant structures on a line bundle
on $\C^n\backslash 0$ are
parametrized by $\C^*$. Since $\gamma$-equivariant
line bundles on $\C^n\backslash 0$ are identified
with line bundles on $H$, this gives another 
proof of the isomorphism $\Pic(H)=\C^*$.

\hfill

The upshot of this argument is that any holomorphic line
bundle $L$ on a Hopf manifold admits a holomorphic connection,
and its monodromy determines the bundle uniquely
(\cite[Subsection 27.5.3]{_OV_book_}).

Since the extension of the pullback of $L$ to $\C^n$
is trivial, the space $\check W$ of generators of the
corresponding $R$-module $W$ (introduced in the proof
of \ref{_G_F_equivariant_on_F_Proposition_}) can be
chosen 1-dimensional, hence ${\cal G}_{\cal F}$
is a subgroup of $\GL(\check W)\times {\cal G}=\C^*\times {\cal G}$.

Using the holomorphic connection on $L$, we lift
the Lie algebra of ${\cal G}$  
to its horizontal lift in $H^0(T\Tot(L))$, defined by the connection.
This horizontal lift is holomorphic, and defines
a $\tilde {\cal G}$-equivariant structure on $L$,
where $\tilde {\cal G}=\C$ is the universal covering
of ${\cal G}$. The bundle $L$ is ${\cal G}$-equivariant
if and only if the kernel of $\tilde {\cal G}\to {\cal G}$
acts on $L$ trivially. 

\hfill

We summarise these arguments in the following proposition.

\hfill

\proposition
Let $L$ be a line bundle on a Hopf manifold.
Then $L$ admits a natural holomorphic connection $\nabla$,
and the group ${\cal G}_L$ is obtained as the algebraic
closure of $\Gamma\times {\cal G}$ in $\C^*\times {\cal G}$,
where $\Gamma$ is the monodromy group of $\nabla$.
\endproof

\hfill

\claim\label{_Z_in_Pic_Hopf_Claim_}
Let $H$ be a classical Hopf manifold,
and $\pi:\; H \to \C P^{n-1}$ be the standard
elliptic projection. Then the pullback map
defines a closed embedding $\pi^*:\; \Pic(\C P^{n-1})\to \Pic(H)$,
from $\Z=\Pic(\C P^{n-1})$ to $\C^*=\Pic(H)$.

\hfill

\proof We need to show that  $\Pic(\C P^{n-1})=\Z$
lifted to $\Pic(H)$ does not
belong to a compact subgroup of $\Pic(H)=\C^*$.
Consider a Gauduchon metric $\omega$ on $H$,
and let $\deg_\omega:\; \Pic(H)\to \R$ denote
the corresponding degree map (\cite[\S 27.3.1]{_OV_book_}).
Since $\pi^*(\calo(1))$ has holomorphic sections,
its degree is positive. Indeed, the degree $\deg L$
is the integral $\int_H c_1(L)\wedge \omega^{n-1}$.
For $L=\pi^*(\calo(i))$, $i>0$, the class $c_1(L)$
is represented by the fundamental class $D$ of 
a zero divisor of its holomorphic section, 
giving $\int_H c_1(L)\wedge\omega^{n-1}=\int_D\omega^{n-1}>0$.
Since the degree of $\pi^*(\calo(1))$ is positive,
$\deg \pi^*(\calo(i))$, $i\in \Z$ is unbounded.
Therefore, the image $\pi^*(\Pic(\C P^{n-1}))$ does not
belong to a compact subgroup of $\Pic(H)$.
\endproof

\subsection{Equivariant structures on trivial line bundles
  on cone varieties}

Let $L$ be a line bundle on $\C^n$, $n>2$.
Then $L$ is trivial, per \ref{_Hopf_Picard_Remark_}.
For a general cone variety $X$, this is 
false: the group $\Pic(X)$ can be non-trivial.
However, the arguments used in Subsection
\ref{_line_on_Hopf_Subsection_} can be applied
to trivial line bundles on $X$ equipped with
a (possibly non-trivial) $\gamma$-equivariant structure.

\hfill

\proposition\label{_C_action_on_L_W_1-dim_Proposition_}
Let $(X_0, \gamma)$ be a cone variety,
and $L$ a trivial line bundle on $X_0$
equipped with a $\gamma$-equivariant action.
Then $L$ is trivialized by a section $f$
such that $\gamma^* f= \lambda f$
for some $\lambda\in \C^*$.

\hfill

\pstep
Applying \ref{_Siu_over_point_Corollary_},
we extend $L$ to the corresponding closed cone variety 
$X$, using the same letter to denote
the extension. We identify the space of sections of $L$ on $X_0$
with $H^0(X, L)$. Let $H^0(L)^\gamma$ be the space
of $\gamma$-finite sections of $L$.
Denote by $H^0(L)^\gamma_\alpha$ 
the root space of $\gamma$-action with the eigenvalue $\alpha \in \C$;
it is defined as $\bigcup_n \ker (\gamma-\alpha\Id)^n$.
Clearly, $H^0(L)^\gamma= \bigoplus_\alpha H^0(L)^\gamma_\alpha$.
Since $H^0(L)^\gamma$ is dense in $H^0_b(U, L)$
(\ref{_finite_dense_in_bounded_Corollary_}), there exists
$f\in H^0(L)^\gamma$ which is non-zero in the origin $x\in X$.
Projecting $f$ to the root spaces, we can assume
that $f\in H^0(L)^\gamma_\lambda$ for some $\lambda\in \C^*$.
In the next step we are going to show that
$|\lambda|$ is strictly smaller than $|\alpha'|$
for all other eigenvalues $\alpha'$, $|\alpha'|>1$, such that
$H^0(L)^\gamma_{\alpha'}\neq 0$.

\hfill

{\bf Step 2:} 
Consider the filtration on
$H^0(\calo_X)$ by the powers of the maximal ideal
${\goth m}_x$. The eigenvalues of $\gamma^*$ acting on
the space $H^0(\calo_X)^\gamma$
of $\gamma$-finite functions are
the same as its eigenvalues on the
associate graded space, which is identified
with $\bigoplus \frac{{\goth m}_x^k}{{\goth m}_x^{k+1}}$.
Since $\frac{{\goth m}_x^k}{{\goth m}_x^{k+1}}$ is a quotient of
$\Sym^k(\frac{{\goth m}_x}{{\goth m}_x^{2}})$,
these eigenvalues are products of the eigenvalues
of $\gamma^*$ acting on $\frac{{\goth m}_x}{{\goth m}_x^{2}}$,
which satisfy $|\alpha_i|>1$.
Therefore, the smallest eigenvalue of $\gamma$ acting
on $H^0(\calo_X)^\gamma$ has multiplicity 1.
The same is true for any free $H^0(\calo_X)^\gamma$-module
$W$ of rank 1, because all eigenvalues of the $\gamma$-action
on $W$ are obtained as $u v_i$, where $u$ is a fixed complex
number, and $v_i$ the eigenvalues of $\gamma^*$ on 
$H^0(\calo_X)^\gamma$. This implies, in particular,
that $H^0(L)^\gamma_\lambda$ has rank 1 when
$|\lambda|>1$ is the smallest possible, and
its generator $f$ is non-zero in a neighbourhood of $x\in X$.

\hfill

{\bf Step 3:} Since $\lambda$ is an eigenvalue 
of multiplicity 1 and $f\in H^0(L)^\gamma_\lambda$,
we have $\gamma^* f= \lambda f$. Since
$f\neq 0$ in a neighbourhood $U$ of $x$,
and $\gamma^n(K)\subset U$ for any compact
$K\subset X$ and $n$ sufficiently big,
this implies that $f\neq 0$ anywhere in $X$.
\endproof

\hfill

\proposition\label{_C_action_on_L_Proposition_}
Let $(X_0, \gamma)$ be a cone variety,
and $L_\lambda$ a trivial line bundle on $X_0$
with $\gamma$-equivariant action multiplying
a  non-vanishing section $f$ by $\lambda\in \C^*$.
Then the  corresponding group ${\cal G}_{L_\lambda}$
is isomorphic to  the Zariski closure of the group generated by
$(\lambda, \gamma)$ in 
$\C^*\times {\cal G}$.

\hfill

\proof 
Let $W$ denote the space of $\gamma$-finite sections of 
$L$, considered as an $R$-module.
Since $L_\lambda$ is generated by a single
non-degenerate section $f$ such that $\gamma^* f =\lambda f$,
the space of generators $\check W$
(introduced in the proof of \ref{_G_F_equivariant_on_F_Proposition_})
can be identified with 
$\C=\langle f \rangle$. By construction,
${\cal G}_{L_\lambda}$ is the closure of
$\langle \gamma^*\rangle$ in $\GL(\check W) \times {\cal G}$. \endproof

%
%

\subsection{Equivariant Picard group}

In light of \ref{_categories_on_cone_Proposition_},
the Picard group of the variety $\frac{X_0}{\langle\gamma\rangle}$
is the same as the group of $\gamma$-equivariant line
bundles on $X_0$. We denote this group $\Pic_\gamma(X_0)$.

\hfill

Let $(X_0, \gamma)$ be a cone variety,
and $Z:= X_0/\C^*_Z$ its projective quotient
obtained as in \ref{_cone_str_on_open_cone_variety_Theorem_},
equipped with an action of $\gamma$ per \ref{_gamma_acts_on_Z_Remark_}. 
Let $\pi:\; X_0\to Z$ be the natural projection. 
By construction, $X_0$ is identified
with the total space of non-zero vectors in a line bundle 
$L$ on $Z$. Clearly, $\pi^*(L)$ is trivial: the corresponding $\C^*$-bundle is 
a principal bundle admitting a tautological
section $\tau$ associating $y$ to every point 
$y\in X_0=\Tot^\circ(L)$.  We choose the equivariant
action on $\pi^*(L)$ which multiplies  
$\tau$ by a constant $\lambda\in\C^*$. Choosing different
$\lambda$, we obtain a natural embedding 
$\xi:\; \C^*\hookrightarrow \Pic_\gamma(X_0)$.
We denote the bundle $\xi(\lambda)$ by $L_\lambda\in \Pic_\lambda(X_0)$.

\hfill

Let $\C^*_Z=\C^*\subset {\cal G}$ be the subgroup used 
in \ref{_C^*_action_orbispace_Theorem_} to
define $Z$. Since ${\cal G}$ is commutative, 
the action of $\gamma$ commutes with $\C^*_Z$, 
hence it passes to $Z= X_0/\C^*_Z$.
We want to describe the $\gamma$-equivariant
line bundles on $X_0$ which are pullbacks
of $\gamma$-equivariant line bundles on $Z$
in terms of the $\gamma$-action. 
More precisely: we want to decompose any
$\gamma$-equivariant line bundle on $X$
into a tensor product of a trivial
line bundle $L_\lambda = \xi(\lambda)$ 
with a non-trivial equivariant
action and the line bundles lifted from $Z$.

Note that the map 
$\Pic_\gamma(Z)\to \Pic_\gamma(X_0)$
is not injective, hence $\Pic_\gamma(Z)$
cannot serve this purpose. Indeed, consider the cone
$(X_0, \gamma)$, where $X_0=\C^n \backslash 0$ and
$\gamma= \const \cdot \Id_{\C^n}$.
Then $\Pic_\gamma(X_0)$ is the Picard group
of a Hopf manifold, identified with $\C^*$ by
\ref{_Hopf_Picard_Remark_}. The action
of $\gamma$ on $Z=\C P^{n-1}$ is trivial
which gives $\Pic_\gamma(Z)=\C^*\times \Z$.
However, the subgroup of line bundles $\Pic_\gamma(X_0)$
which are pullbacks is $\Z$ (\ref{_Z_in_Pic_Hopf_Claim_}), 
and it is strictly
smaller than $\Pic_\gamma(Z)=\C^*\times \Z$. 

\hfill

Instead of $\Pic_\gamma(Z)$, we are forced to use
a smaller group $\Pic_\gamma(Z)_0$, defined
as follows. First, $Z$ is already equipped
with an ample line bundle $E$, with
$X_0=\Tot^\circ(E^*)$ (\ref{_C^*_action_orbispace_Theorem_}).
This bundle is by construction $\gamma$-equivariant.
Then $Z= \Proj\left(\bigoplus\limits_{i=0}^\infty H^0(Z, E^{\otimes i})\right)\!,$
and ${\cal G}$ is the Zariski closure of
$\langle\gamma\rangle$ acting on
$\bigoplus_{i=0}^d H^0(Z, E^{\otimes i})$
for $d$ sufficiently big to contain all generators
of the ring $\bigoplus_i H^0(Z, E^{\otimes i})$.

Let $L\in \Pic_\gamma(Z)$ be
a $\gamma$-equivariant bundle. Consider the group 
${\cal G}_L$ defined as in \ref{_G_F_equivariant_on_F_Proposition_};
it is obtained as the Zariski closure of the $\gamma$-action
on $\bigoplus_{i=0}^d H^0(Z, L\otimes E^{\otimes i})$ 
for $d$ sufficiently big to contain all generators
of $\bigoplus_i H^0(Z, L\otimes E^{\otimes i})$
considered as a module over 
$H^0(\calo_{X_0})=\bigoplus_i H^0(Z, E^{\otimes i})$.

Clearly, this action defines a $\gamma$-equivariant
structure on the pullback of $L$ to $X_0$,
and ${\cal G}_L$ is obtained as the Zariski closure of the
action of $\langle \gamma\rangle$.
Therefore, ${\cal G}_L$ is equipped
with a natural homomorphism to ${\cal G}$.

Define $\Pic_\gamma(Z)_0$ as the set of all
$L\in \Pic_\gamma(Z)_0$  such that the natural
map $P:\; {\cal G}_L\to {\cal G}$
is an isomorphism. 

\hfill

\remark\label{_C^*_Z_equiv_are_G_equiv_Remark_}
A $\gamma$-equivariant bundle $L\in \Pic_\gamma(X_0)$
can be obtained as a pullback of $L_0\in \Pic_\gamma(Z)$
only if $L_0\in \Pic_\gamma(Z)_0$. Indeed, the group
${\cal G}_L$ is by construction equal to ${\cal G}_{L_0}$.
Since $L= \pi^*(L_0)$, its sections are obtained
by multiplying sections of $L_0$ with holomorphic 
functions on $X_0$. Therefore, the preimage of
$\C^*_Z$ in ${\cal G}_{L}={\cal G}_{L_0}$
acts on $L$ in the same way as on $H^0(\calo_{X_0})$,
that is, as $\C^*$ (and not as a bigger group).
This implies that the projection
$P:\; {\cal G}_{L}\to {\cal G}$
is an isomorphism on $P^{-1}(\C^*_Z)$,
and therefore it has no kernel.

\hfill

\claim\label{_pullb_inj_on_Pic_Claim_}
The pullback map $\Pic_\gamma(Z)_0\to \Pic_\gamma(X_0)$
is injective.

\hfill

\proof Let $\pi:\; X_0\to Z$ be the natural projection.
Consider $L\in \Pic_\gamma(Z)_0$.
Its pullback $\pi^*(L)$ is $\gamma$-equivariant
by construction. However, the $\gamma$-equivariant
structure is induced by the ${\cal G}_L$-equivariant structure,
and the preimage of $\C^*_Z$ in ${\cal G}_L$ is $\C^*$,
by definition of $\Pic_\gamma(Z)_0$,
forcing $\pi^*(L)$ to be $\C^*_Z$-equivariant. Then 
$L$ is reconstructed from $\pi^*(L)$
as the sheaf of $\C^*_Z$-invariant sections.
\endproof

\hfill

The next theorem is the main result of this section,
giving an  explicit description of the Picard group of
$X_0/\langle\gamma\rangle$ in terms of the
$\gamma$-equivariant line bundles on $Z$.

\hfill

\theorem\label{_Picard_X_0_Theorem_}
Let $(X_0, \gamma)$ be a cone variety,
$Z:= X_0/\C^*_Z$ its projective quotient
obtained as in \ref{_cone_str_on_open_cone_variety_Theorem_},
and $\xi:\; \C^*\hookrightarrow \Pic_\gamma(X_0)$
the embedding constructed above.
Then the tensor product of line bundles map 
$\Pic_\gamma(Z)_0\times \xi(\C^*)\to \Pic_\gamma(X_0)$
is surjective, and injective on both components.

\hfill

\pstep
Injectivity on the second component is
inherent in the construction, and
injectivity on the first component is shown in
\ref{_pullb_inj_on_Pic_Claim_}. It remains to
prove the surjectivity.

Let $L_\lambda$ denote 
a trivial line bundle on $X_0$ equipped with a $\gamma$-equivariant
action multiplying a section by $\lambda\in \C^*$, as in 
\ref{_C_action_on_L_Proposition_}.
Let $E$ be a $\gamma$-equivariant line bundle on $X_0$,
and let $E_\lambda:= E\otimes L_\lambda$. 
Since $\C^*_Z$-equivariant bundles
belong to $\pi^*(\Pic_\gamma(Z)_0)$,
surjectivity would follow if we 
show that $E_\lambda$ is $\C^*_Z$-equivariant
for some $\lambda\in \C^*$.

\hfill

{\bf Step 2:}
Let $S\subset {\cal G}_E$, 
$S_\lambda\subset {\cal G}_{L_\lambda}$ be 1-parametric
complex Lie groups which surjectively project
to the group $\C^*_Z\subset {\cal G}$ used in the
definition of $Z=X_0/\C^*_Z$. The natural
projections $S \to \C^*_Z$ and $S_\lambda \to \C^*_Z$
are surjective homomorphisms of connected complex
Lie groups of dimension 1. If $S=\C$, the kernel
of this homomorphism is $\Z$, and if $S=\C^*$,
it is a finite cyclic group. Let
$u, u_\lambda$ be the generators of the
kernels of the projections $S\to \C^*_Z$, $S_\lambda \to \C^*_Z$.
Since $u$ and $u_\lambda$ are $\gamma$-invariant
automorphisms of a line bundle, and any $\gamma$-invariant
function on $X_0$ is constant (because it passes to the quotient $X_0/\langle\gamma\rangle$ which is compact), both $u$ and $u_\lambda$ act
on the line bundles $E, L_\lambda$ as constants; 
accordingly, we will treat $u, u_\lambda$ as elements of $\C^*$.
Clearly, $E_\lambda$ is $\C^*_Z$-invariant
if and only if $u u_\lambda=1$.

\hfill

{\bf Step 3:}
By definition, $L_\lambda\otimes L_\mu=L_{\lambda\mu}$.
Such bundles generate the subgroup
 $\xi(\C^*)\subset \Pic_\gamma(X_0)$,
isomorphic to $\C^*$. Identifying $\Pic_\gamma(X_0)$
with $\Pic(X_0/\langle\gamma\rangle)$,
we can treat $L_\lambda$ as a flat line bundle
on $X_0/\langle\gamma\rangle$. Note that 
$u_\lambda$ is equal to the monodromy of the corresponding
connection on an orbit of the $\C^*_Z$-action.
Therefore, the map $L_\lambda \mapsto u_\lambda$
defines a group homomorphism from $\xi(\C^*)$ to
$\C^*$. To show that $u u_\lambda=1$,
it sufficies to prove that this homomorphism
is surjective, that is, non-trivial.
Arguing ad absurdum, we can suppose that 
$u_\lambda=1$ for any $\lambda$.

\hfill

{\bf Step 4:} 
Let $\gamma_0\in \C^*_Z\subset {\cal G}$ be an element such that
the quotient $\C^*_Z/\langle \gamma_0\rangle$ is an elliptic curve.
The bundle $L_\lambda$, restricted to
the orbit of $\C^*_Z$,
has trivial monodromy, because this monodromy is equal to $u_\lambda$.
This implies, in particular, that the section
$f$ used in the definition of $L_\lambda$ is invariant under the action
of the Lie algebra of $\C^*_Z\subset {\cal G}_{L_\lambda}$.
Then $\gamma_0^*(f)=f$.

Since ${\cal G}_{L_\lambda}$ multiplies the section
$f$ by a constant, the metric which satisfies
$|f|=1$ satisfies $|a^* f|=\chi(a) |f|$,
for a fixed character $\chi:\; {\cal G}_{L_\lambda}\to \R^{>0}$
and any $a\in {\cal G}_{L_\lambda}$.
Since $f$ is $\gamma_0$-invariant,
it is obtained as the pullback of a section of
a line bundle on the variety $X_0/\langle \gamma_0\rangle$,
which is compact; then $|f|$ is bounded on $X_0$,
bringing a contradiction for $|\lambda| \neq 1$,
because $f$ satisfies $|f(x)| =|\lambda| |f(\gamma(x))|$.
We proved that the map $\lambda\to u_\lambda$
is surjective, finishing the proof of
\ref{_Picard_X_0_Theorem_}.
\endproof

\hfill

\remark
The kernel of the surjective map 
$\Pic_\gamma(Z)_0\times \xi(\C^*)\to \Pic_\gamma(X_0)$
defined in \ref{_Picard_X_0_Theorem_}
is identified with the intersection of these two
factors. By definition of $\xi$, this intersection is the set of all
$L\in \Pic_\gamma(Z)_0$ such that the pullback
of $L$ to $X_0$ is trivial as a holomorphic line bundle.

\hfill

\proposition\label{_intersection_in_Picard_Proposition_}
The intersection $\Lambda$ of the images $\Pic_\gamma(Z)_0$
and $\xi(\C^*)$ in $\Pic_\gamma(X_0)$ is a discrete 
subgroup in $\xi(\C^*)=\C^*$.

\hfill

\pstep
Let $j:\; \Pic_\gamma(Z)_0\to \Pic_\gamma(X_0)$
be the injective map defined above
(\ref{_Picard_X_0_Theorem_}). We need to show that 
$\xi(\C^*) \cap j(\Pic_\gamma(Z)_0)$ is discrete in $\xi(\C^*)$.
Otherwise $\xi(\C^*) \cap j(\Pic_\gamma(Z)_0)= \xi(\C^*)$
because it is an algebraic subgroup in $\xi(\C^*)$.
Arguing ad absurdum,
suppose that $j(\Pic_\gamma(Z)_0)\supset \xi(\C^*)$.
Then all elements of $j^{-1}(\xi(\C^*))$
are infinitely divisible in $\Pic(Z)$, hence their
first Chern class vanishes. Indeed, a non-zero
element $H^2(Z, \Z)$ cannot be infinitely divisible.
We are going to prove
that (some of) these bundles have non-zero first Chern class,
arriving at contradiction.

\hfill

{\bf Step 2:}
In Step 3, we prove that
there exists a bundle $L_\lambda\in \xi(\C^*)$
which admits a $\gamma$-invariant section $s$
with non-trivial zero divisor.
Since $\cal G$ is the smallest algebraic
group containing $\gamma$, the section
$s$ is also $\cal G$-invariant. By definition,
$\Pic_\gamma(Z)_0$ is the group of $\gamma$-equivariant
line bundles on $Z$ which are
lifted to $\C_Z$-equivariant bundles 
on $X_0$. Since $L_\lambda$ admits
a $\C_Z$-invariant section, it is 
$\C_Z$-equivariant; indeed, any bundle
$\calo(D)$ for any $\C^*_Z$-invariant divisor
$D$ is $\C_Z$-equivariant. Therefore, 
$L_\lambda$ is the pullback
of a line bundle $\underline L_\lambda \in \Pic_\gamma(Z)_0$
(\ref{_Roberts_invariant_sections_Theorem_}).
Let $\underline s$ be the corresponding section of 
$\underline L_\lambda$. As shown in Step 1,
$c_1(\underline L_\lambda)\neq 0$ implies 
\ref{_intersection_in_Picard_Proposition_},
and $c_1(\underline L_\lambda)\neq 0$ 
because $c_1(\underline L_\lambda)$
is cohomologous to the fundamental class
of the zero divisor of 
$\underline s\in H^0(Z,\underline L_\lambda)$.

\hfill

{\bf Step 3:} It remains to produce
$L_\lambda\in \xi(\C^*)$
which admits a $\gamma$-invariant section $s$
with non-trivial zero divisor.
Let $\psi:\; X_0 \to \C^n\backslash 0$
be the holomorphic embedding constructed
in \ref{_cone_embedded_Theorem_}, and $A \in \GL(\C^n)$ a linear contraction
which restricts to $\gamma\in \Aut(X_0)$.
Using the Jordan normal form of $A$,
we construct an $A$-invariant linear
subspace $W\subset \C^n$ such that
$\dim (W\cap X_0)=\dim X_0-1$.\footnote{This argument,
which proves a theorem 
originally due to Ma. Kato 
\cite{_Kato:subvarieties_}, was first given in 
\cite{_OV:surf_in_lck_pot_}.
The Jordan normal form of $A$ gives a full flag of linear
subspaces in $\C^n$, hence provides a full flag of Hopf submanifolds
in $H=\frac{\C^n\backslash 0}{\langle A\rangle}$.
Intersecting this flag with a subvariety 
$X_0/\langle \gamma\rangle\subset H$ gives a full flag
of subvarieties in $X_0/\langle \gamma\rangle$.}
The corresponding
line bundle $\calo(W)$ is by
construction $\langle A \rangle$-equivariant;
therefore, its restriction $\calo(W\cap X_0)$
is a $\gamma$-equivariant bundle.
Passing to $Z=X_0/\C_Z$, the divisor $W$
becomes a hyperplane section of $Z$, hence it is non-trivial.
Since $\Pic_\gamma(\C^n\backslash 0)=\C^*$,
this group coincides with $\xi(\C^*)$ (\ref{_Hopf_Picard_Remark_}), hence
$\calo(W)\in \xi(\C^*)$, which also implies
that $\calo(W\cap X_0)\in \xi(\C^*)$.
\endproof

\hfill

As a closing note to the proof, we observe the following.
Let $L_\lambda\in \xi(\C^*)$ be a $\gamma$-equivariant
bundle on a cone variety. Taking the algebraic
closure of $\langle \gamma\rangle$, we obtain
an algebraic group ${\cal G}_{L_\lambda}$
which is projected to ${\cal G}$, with the
fiber which is 0-dimensional or 1-dimensional.
This projection is an isomorphism
if and only if 
$L_\lambda \in  \xi(\C^*)\cap\Pic_\gamma(Z)_0$.

\subsection{Equivariant coherent sheaves on cone varieties}

By the Lie-Kolchin theorem  (\cite[Theorem 16.30]{_Milne:AG_}),
for any  collection of commuting operators $A_i\subset \End(V)$,
the finite-dimensional representation space $V$ admits an invariant full flag. 
We are going to describe the corresponding filtration more explicitly,
as follows.

\hfill

Let ${\goth S}=\{A_1, A_2, ...\}\in \End(V)$ be a collection of 
commuting linear operators acting on a finite-dimensional 
vector space $V$, and ${\goth F}$ the set of functions
$\alpha:\; {\goth S} \to \C$
taking each $A_i$ to one of its eigenvalues.
The corresponding root decomposition
$V =\bigoplus_{\alpha\in \goth F} \bigcap_{A_s\in \goth S}V_{\alpha(A_s)}$ is 
the decomposition of $V$ onto the union of the 
root spaces,
where $V_{\alpha(A_s)}$ is the root space of $A_s$
associated with the eigenvalue $\alpha(A_s)$.

\hfill

\definition\label{_generalized_weight_Definition_}
On each $V_\alpha$, each $A_s$ acts as a Jordan block
with an eigenvalue $\alpha(A_s)$. Thus, we have a 
Jordan cell filtration for each $A_s\in {\goth S}$,
$F_1^{\alpha, s} \subset F_2^{\alpha, s}\subset ...$,
where $F_i^{\alpha, s}$ denotes the kernel of
$(A_s -\alpha(s))^i$.
The decomposition $V =\bigoplus_{\goth F} V_\alpha$ 
is called {\bf the generalized weight decomposition},
and the corresponding filtration $F_i^{\alpha, s}$
{\bf the generalized Jordan filtration}.

\hfill

\remark\label{_gen_Jordan_scalars_Remark_}
On the (smallest) subquotients 
of the generalized Jordan filtration, all operators $A_s$ act
as scalars. In other words, all $A_s$ acts by scalars on
the associated graded space of the generalized Jordan filtration.

\hfill

\remark
\ref{_gen_Jordan_scalars_Remark_}
is a special case of Lie theorem, which
claims that any representation $V$ of a
solvable Lie algebra $\g$ admits a filtration
$V=V_n \supset V_{n-1} \supset ... \supset V_0=0$
such that every element of $\g$ acts on $V_i/V_{i-1}$ by scalars.
In our case $\g$ is the Lie algebra generated
by $A_1, A_2, ...$, which is commutative,
hence solvable.

\hfill

\proposition\label{_filtration_on_reflexive_Proposition_}
Let ${\cal H}$ be a commutative algebraic
group acting by automorphisms on a $\gamma$-equivariant 
reflexive coherent sheaf ${\cal F}$ over an open cone variety $X_0$.
Then ${\cal F}$ admits an ${\cal H}$-invariant, 
$\gamma$-equivariant filtration, such that
${\cal H}$ acts by scalars on its subquotients.

\hfill

\proof
Denote by $S\subset X_0$ the singular locus of ${\cal F}$;
it has codimension $\geq 2$, because ${\cal F}$ is reflexive
and $X_0$ is normal (\cite{_Hartshorne:stable_reflexive_}).
For each $h\in {\cal H}$, and $z\in X_0\backslash S$,
denote the eigenfunctions of $h$ on the fiber of
${\cal F}$ in $z$ by $\alpha_1, ..., \alpha_s$.
This defines a $\gamma$-invariant holomorphic map
from $X_0\backslash S$ to $\Sym^s(\C)$; since $X_0/\langle\gamma\rangle$
is compact, such a map is constant by maximum principle.
This implies that the eigenvalues of $h$ are constant, and
${\cal F}$ admits a root decomposition with
respect to $h$: 
\[ {\cal F}= \bigoplus {\cal F}_{\alpha_s}.
\]
Similarly, for a collection of operators
we have a generalized weight decomposition
${\cal F} =\bigoplus_{\alpha \in \goth F} {\cal F}_\alpha$ 
(\ref{_generalized_weight_Definition_}).
Taking the subquotients of the corresponding
generalized Jordan filtration and using
\ref{_gen_Jordan_scalars_Remark_}, we obtain
a filtration such that ${\cal H}$ acts
by (constant) scalars on its subquotients.
\endproof

\hfill

The main result of this section is 
\ref{_filtra_with_subquot_lifted_up_to_L_i_Theorem_}
below. However, before we formulate it, we need to 
explain how its assumptions fit with the rest of this section.

\hfill

\remark
\ref{_filtra_with_subquot_lifted_up_to_L_i_Theorem_} is stated for 
${\cal G}_{\cal F}$-equivariant sheaves on a closed cone variety
$(X,\gamma)$, but  we shall apply it (\ref{_filtrable_Corollary_}) to reflexive 
sheaves on $X_0$ (which extend naturally to $X$ by Siu's
extension theorem).
Let us recall how the ${\cal G}_{\cal F}$-equivariance
is proven for reflexive sheaves.
Let ${\cal G}$ be the algebraic closure of $\langle\gamma\rangle$
in $\Aut(X)$ (\ref{_cone_str_on_open_cone_variety_Theorem_}). 
As in \ref{_G_F_equivariant_on_F_Proposition_}, 
we extend the $\gamma$-equivariant action
on ${\cal F}$ to the ${\cal G}_{\cal F}$-equivariant action,
where ${\cal G}_{\cal F}$ is a commutative complex 
Lie group equipped with a surjective map to ${\cal G}$.
The proof of \ref{_filtra_with_subquot_lifted_up_to_L_i_Theorem_} 
involves replacing ${\cal F}$ by its 
${\cal  G}_{\cal  F}$-invariant subquotients, and we can no longer
guarantee that these subquotients are torsion-free. This
is why we drop the reflexivity assumption, and replace
it by the weaker assumption of existence of a ${\cal G}_{\cal F}$-action.

\hfill

\theorem\label{_filtra_with_subquot_lifted_up_to_L_i_Theorem_}
Let $(X, \gamma)$ be a closed cone variety, 
and ${\cal F}$ a $\gamma$-equivariant 
coherent sheaf on $X$. 
Assume that the $\gamma$-equivariant action
is extended to an action of 
an algebraic group ${\cal  G}_{\cal F}$ 
which projects to ${\cal G}$ surjectively. 
Assume that ${\cal F}$ admits a ${\cal K}_{\cal F}$-invariant, 
$\gamma$-equivariant filtration, such that the kernel
${\cal K}_{\cal F}$
of the map ${\cal  G}_{\cal F}\to {\cal G}$ acts by scalars on its
subquotients.\footnote{We apply 
\ref{_filtra_with_subquot_lifted_up_to_L_i_Theorem_} to
reflexive sheaves, where this assumption follows from
\ref{_filtration_on_reflexive_Proposition_}.}
Consider the $\gamma$-invariant
action of $\C^*_Z\subset {\cal G}$
constructed in \ref{_cone_str_on_open_cone_variety_Theorem_},
and let $\pi:\; X_0 \to Z=X_0/\C^*_Z$  
be the projection map. Since $\gamma$ commutes with $\C^*_Z$,
its action descends to the orbivariety $Z$, and we use the same
letter $\gamma$ to denote the corresponding automorphism of $Z$
(\ref{_gamma_acts_on_Z_Remark_}).  Then there exists
a collection of $\gamma$-equivariant coherent
sheaves ${\cal F}_1$, ..., ${\cal F}_k$ on $Z$
and a collection of $\gamma$-equivariant line
bundles
$L_1, ..., L_k\in \xi(\C^*)\subset \Pic_\gamma(X_0)$ such that
${\cal F}$ admits a $\gamma$-invariant filtration with subquotients
isomorphic, as $\gamma$-equivariant sheaves,
to $L_i \otimes \pi^*{\cal F}_i$.

\hfill

%

\pstep  We need to show that for each
subquotient of ${\cal F}$ with scalar action of
${\cal K}_{\cal F}$, there exists a $\gamma$-equivariant
line bundle $L$ such that ${\cal F}\otimes L$
is $\C^*_Z$-equivariant, where $\C^*_Z\subset {\cal G}$
is the group defined in \ref{_C^*_action_orbispace_Theorem_} 
and satisfying $Z=X_0/\C^*_Z$.

If the action of ${\cal K}_{\cal F}$ 
is trivial, then ${\cal G}={\cal G}_{\cal F}$,
and the existence of 
the $\C^*_Z$-equivariant action on ${\cal F}$
is automatic. To finish the proof of
\ref{_filtra_with_subquot_lifted_up_to_L_i_Theorem_}, 
it remains to show that for any 
$\gamma$-equivariant coherent sheaf ${\cal F}$
on $X_0$ with scalar action of ${\cal K}_{\cal F}$,
there exists a $\gamma$-equivariant 
line bundle $L\in \xi(\C^*)$ on $X_0$ such that ${\cal F}\otimes L$
is $\C^*_Z$-equivariant (that is, satisfies
${\cal K}_{{\cal F}\otimes L}=\{1\}$).

\hfill

{\bf Step 2:}
Recall how the group ${\cal G}_{\cal F}$ is defined  
(\ref{_G_F_equivariant_on_F_Proposition_}).\footnote{\ref{_G_F_equivariant_on_F_Proposition_}
assumes reflexivity to construct the
finite-dimensional space $\check W$ of generators;
however, since ${\cal F}$ is already assumed
${\cal G}_{\cal F}$-equivariant, we don't need
reflexivity.}
Consider a finite-dimensional $\gamma$-invariant
space $\check W\subset H^0(X, {\cal F})$ of sections
of ${\cal F}$ generating ${\cal F}$ over
$H^0(\calo_X)$. Then 
${\cal G}_{\cal F}$ is the algebraic closure of
$\gamma^*$ acting on $\check W\oplus H^0(\calo_X)$.
Then ${\cal K}_{\cal F}$ is the kernel
of the natural projection from ${\cal G}_{\cal F}$ 
to ${\cal G}\subset \Aut(H^0(\calo_X))$. 
The action of ${\cal G}_{\cal F}$ on $\check W$
defines a generalized Jordan filtration
(\ref{_generalized_weight_Definition_}),
and acts by diagonalizable operators on its
associated graded space.
To prove
\ref{_filtra_with_subquot_lifted_up_to_L_i_Theorem_},
we can always replace the sheaf ${\cal F}$
by its associated graded sheaf for an equivariant filtration.
Replacing ${\cal F}$ by an associated graded
quotient if necessarily,  we can always assume that
$\check W$ is generated by the eigenfunctions
of ${\cal G}_{\cal F}$. Without restricting the
generality, we can always replace ${\cal F}$ 
by one of these eigenspace sheaves. 

After this is done, we can assume that the sheaf ${\cal F}$ is 
generated by sections $v_1, ..., v_n$ such that 
$\gamma^* v_i=\alpha v_i$, for some $\alpha \in \C^*$.
This means that ${\cal F}$ is a quotient sheaf 
of $L_\alpha^{\oplus^n}$, where $L_\alpha:= \xi(\alpha)$ is the line
bundle constructed in \ref{_Picard_X_0_Theorem_}.
Taking $\lambda=\alpha^{-1}$, we obtain that
${\cal F}\otimes L_\lambda$ is a quotient of 
a free sheaf generated by $\gamma$-invariant
sections, hence it is $\C^*_Z$-equivariant.
\endproof

\subsection{Filtrable $\gamma$-equivariant sheaves on cone varieties}

\definition
Let $Z$ be a projective variety or orbivariety,
$G$ a group acting on $Z$, and ${\cal F}$
a $G$-equivariant coherent sheaf on $Z$.
We say that ${\cal F}$ is {\bf filtrable}
if it admits a $G$-equivariant filtration
$0={\cal F}_0 \subset {\cal F}_1\subset ... \subset {\cal F}_N={\cal F}$
such that the subquotients $\frac{{\cal F}_i}{{\cal F}_{i-1}}$
have rank $\leq 1$. 

\hfill

\remark
We have already discussed the filtrability in
Subsection \ref{_filtrability_Subsection_}.
For general $G$, the filtrability
is not always given even when $Z$ is projective.
For example, if $V$ is an irreducible
representation of $G$, which acts trivially on $Z$,
and ${\cal F}= V \otimes_\C \calo_Z$, all $G$-equivariant
subsheaves of ${\cal F}$ have rank 0 or rank $\dim V$.
However, when $G$ is abelian, filtrability 
of $G$-equivariant sheaves is easy to prove.

\hfill

\claim\label{_filtrability_for_G_sheaves_Claim_}
Let ${\cal G}$ be an abelian group acting on 
a projective (orbi)-variety $Z$,\footnote{As always,
when ${\cal G}$ acts on a projective variety, we
assume that its action preserves an ample bundle;
such an action is called {\bf linearizable}.} and ${\cal F}$
a ${\cal G}$-equivariant coherent sheaf on $Z$.
Then ${\cal F}$ is filtrable as a  ${\cal G}$-equivariant sheaf.

\hfill

\proof
Let $\calo(1)$ be an ample line bundle on $Z$.
Then ${\cal F}(n):={\cal F}\otimes_{\calo_Z} \calo(n)$ is globally
generated for $n \gg 0$. Consider the ${\cal G}$-invariant generalized
Jordan filtration $W_1\subset W_2 \subset ... \subset W_N=H^0({\cal F}(n))$; 
its subquotients are 1-dimensional vector spaces. 
The corresponding filtration on ${\cal F}$
is obtained as follows: ${\cal F}_i(n)$ is the subsheaf
of ${\cal F}(n)$ generated by the sections in $W_i$, and
${\cal F}_i:= {\cal F}_i(n)\otimes _{\calo_Z} \calo(-n)$.
Since ${\cal F}(n)$ is globally generated, ${\cal F}_N={\cal F}$.
Clearly, the subquotients $\frac{{\cal F}_i}{{\cal F}_{i-1}}$
have rank $\leq 1$.
\endproof

\hfill

\corollary\label{_filtrable_Corollary_}
Let $(X_0, \gamma)$ be a cone variety
and ${\cal F}$ a $\gamma$-equivariant
reflexive sheaf on $X_0$, or, equivalently,
a coherent sheaf on $X_0/\langle \gamma\rangle$.
Assume that $\dim_\C X_0 \geq 3$. Then 
${\cal F}$ admits a $\gamma$-equivariant 
filtration with all associated graded subquotients
$\gamma$-equivariant coherent sheaves 
of rank $\leq 1$.

\hfill

\proof
By \ref{_Siu_over_point_Corollary_},
${\cal F}$ is extended to
$\langle\gamma\rangle$-equivariant sheaf on the
closed cone $X$. By
\ref{_G_F_equivariant_on_F_Proposition_},
${\cal F}$ admits a ${\cal G}_{\cal F}$-equivariant
action, with the projection of ${\cal G}_{\cal F}$
to $X$ equal to ${\cal G}$, defined as the
Zariski closure of $\langle\gamma\rangle$
(\ref{_algebraic_str_on_cone_uniqueness_Theorem_}).
In this situation, we can apply
\ref{_filtra_with_subquot_lifted_up_to_L_i_Theorem_},
and obtain that ${\cal F}$ admits a filtration with the
subquotients ${\cal F}_i\otimes L_{\lambda_i}$
such that ${\cal F}_i$ is a ${\cal G}$-equivariant
pullback of a coherent sheaf $\underline {\cal F}_i$
on the projective orbivariety $Z:= X_0/\C^*_Z$ 
(\ref{_C^*_action_orbispace_Theorem_}).
Since a ${\cal G}$-equivariant
coherent sheaf on $Z$ admits a filtration 
with subquotients of rank $\leq 1$  
(\ref{_filtrability_for_G_sheaves_Claim_}), 
each sheaf ${\cal F}_i\otimes L_{\lambda_i}$
is filtrable as well.
\endproof

%

\hfill

\noindent{\bf Acknowledgements:} We thank Cezar Joi\c ta
for the reference \cite{_Kaup_}, Marian Aprodu for valuable suggestions, and Jason Starr for his
insightful comments on Mathoverflow. We are very grateful
to Dmitri Korshunov for his help (\cite{_Korshunov:MO_})
related to polynomial automorphisms of $\C^n$ used in 
Subsection \ref{_Henon_automo_Subsection_}.

\hfill

\hfill

{\scriptsize

\noindent {\sc Liviu Ornea\\
{\sc University of Bucharest, Faculty of Mathematics and Informatics, \\14
Academiei str., 70109 Bucharest, Romania}, \\
also:\\
Institute of Mathematics ``Simion Stoilow" of the Romanian
Academy\\
21, Calea Grivitei Str.
010702-Bucharest, Romania}\\
{\tt lornea@fmi.unibuc.ro,   liviu.ornea@imar.ro}

\hfill

\noindent {\sc Misha Verbitsky\\
 Instituto Nacional de Matem\'atica Pura e
Aplicada (IMPA) \\ Estrada Dona Castorina, 110\\
Jardim Bot\^anico, CEP 22460-320\\
Rio de Janeiro, RJ - Brasil }\\
{\tt verbit@impa.br}

}

\addcontentsline{toc}{section}{References}
{\small
\begin{thebibliography}{100}

\bibitem[AH]{_Abramovich_Hassett_}
D. Abramovich,  and B. Hassett, {\em Stable varieties
with a twist,} in {\em Classification of Algebraic Varieties} {\bf 3}, 
1-38, EMS Ser. Congr. Rep., EMS, Z\"urich, 2011.



\bibitem[AN]{_Andreotti_Narasimhan_} A. Andreotti, R. Narasimhan, {\em Oka's Heftungslemma and the Levi Problem for Complex Spaces}, Trans. Amer. Math. Soc. {\bf 111} (1964), 346-366.

\bibitem[AS]{_Andreotti_Siu:embeddings_} 
A. Andreotti,
  Y. T.  Siu,  {\em Projective embeddings of pseudoconcave
    spaces}, Ann. Scuola Norm. Sup. Pisa {\bf 24}, 231-278
  (1970).

\bibitem[AT]{_Aprodu_Toma_}
M. Aprodu, M. Toma, 
{\em Une note sur les fibr'es holomorphes non-filtrables}, 
C. R. Math. Acad. Sci. Paris {\bf 336} (2003), no. 7, 581-584.

\bibitem[Arn]{_Arnold:ODE+_} 
V. I. Arnol'd, {\em Geometrical
  Methods in the Theory of Ordinary Differential
  Equations,} Grundlehren der mathematischen Wissenschaften
  {\bf 250}, Springer, 1996.

\bibitem[Art]{_Artin:Approximation_}
M. Artin,  {\em 
Algebraic approximation of structures over complete local rings,} 
Publ. Math. IH\'ES, {\bf 36} (36) (1969). 23-58.


\bibitem[AM]{_Atiyah_MacDonald_} 
M. F. Atiyah, I. G. MacDonald, {\em Introduction to commutative
algebra,} Addison-Wesley, 1969.

\bibitem[Ba1]{_Baily:decomposition_}
W.L. Baily, 
{\em The decomposition theorem for V-manifolds},
Amer. J. Math. 78 (1956), 862-888.
32.00 (31.00)

\bibitem[Ba2]{_Baily:embedding_}
W.L. Baily, \emph{On the imbedding of V-manifolds
    in projective spaces}, Amer. J. Math. {\bf 79} (1957), 403-430.


\bibitem[BP]{_Banica_Le_Potier_} 	
C. B\u anic\u a, J. Le Potier
{\em Sur l'existence des fibr\'es vectoriels holomorphes 
sur les surfaces non-alg\'ebriques}, 
J. Reine Angew. Math. {\bf 378} (1987), 1-31.



\bibitem[BHPV]{_BHPV_} 
Barth, W., Hulek, K., Peters, C., van de Ven, A., 
{\em Compact complex surfaces,} 
Second enlarged edition. Springer Verlag, Berlin-Heidelberg, 2004. 




\bibitem[BeS]{_Behnke_Stein_} H. Behnke, K. Stein, {\em Konvergente Folgen Von Regularit\"atsbereichen und die Meromorphiekonvexit\"at}, Math.  Ann. {\bf 116} (1939) 204-216.

\bibitem[Ber]{_Berteloot:PD_}
F. Berteloot, 
{\em M\'ethodes de changement d'\'echelles en analyse complexe},
 Ann. Fac. Sci. Toulouse Math., Tome
{\bf XV}, 3 (2006), p.427-483


\bibitem[Bor]{_Borel:alg_groups_}
A. Borel, {\em Linear algebraic groups},
 GTM {\bf 176}, Springer-Verlag, 1991.

\bibitem[BT]{_Borel_Tits:Groupes_Reductifs_} A. Borel,
  J. Tits, {\em Groupes r\'eductives}, Publications
  Math\'ematiques de l'IH\'ES, {\bf 27} (1965),  55-151.

\bibitem[BG]{_Boyer_Galicki_} 
C. P.  Boyer, K.  Galicki, {\em Sasakian geometry}, 
Oxford Univ. Press, Oxford, 2008.

\bibitem[Bo]{_Boureau:Hopf_}
P. Boureau,
{\em Normal forms and geometric structures on Hopf
  manifolds},
arXiv:2501.10346

\bibitem[Br]{_Brinzanescu:LMN_1624_} 
V. Br\^inz\u{a}nescu
{\em Holomorphic Vector Bundles over Compact Complex Surfaces},
LNM {\bf 1624}, Springer Verlag, 1996.

\bibitem[BM]{_Brinzanescu_Moraru_}
V. Br\^inz\u anescu, R. Moraru {\em Holomorphic rank-2 vector bundles
on non-K"ahler elliptic surfaces,} Ann. Inst. 
Fourier  {\bf 55} No. 5 (2005), 1659-1683.

\bibitem[Ca]{_Calle:Stacks_}
M. E. Calle,
{\em Orbispaces as stacks: geometry and examples, talk notes},
\url{https://hiroleetanaka.com/workshop-2025/notes-05.pdf}


\bibitem[CD]{_Cantat_Dujardin_}
S. Cantat, R. Dujardin, {\em Holomorphically
conjugate polynomial automorphisms of $\mathbb {C}^2$ are
polynomially conjugate,}  Bull. London Math. Soc., {\bf 56} (2024), 
3745-3751.

\bibitem[CT]{_Clarke_Tipler_}
A. Clarke, C. Tipler, {\em Blowing-up hermitian Yang-Mills connections,}
 Ann. Scuola Norm. Sup. Pisa Cl. Sci., {\bf 23}. \url{https://doi.org/10.2422/2036-2145.202311_009} 

\bibitem[Cl]{_Clemens:degene_}
C. H. Clemens, 
{\em Degeneration of K\"ahler manifolds}, 
Duke Math. J. {\bf 44} (1977), no. 2, 215-290.



\bibitem[Del]{_Deligne:Tannakian_}
P. Deligne, {\em Cat\'egories tannakiennes}, 
The Grothendieck Festschrift, Vol.
II, Progr. Math. {\bf 87}, Birkh\"auser 
Boston, Boston, MA, 1990,  111-195.


\bibitem[Dem]{_Demailly:AG_} J.-P.  Demailly, {\em Complex analytic and differential geometry}, {\scriptsize\url{ https://www-fourier.ujf-grenoble.fr/~demailly/manuscripts/agbook.pdf}}


\bibitem[Dl]{_Dloussky:Kato_}
G. Dloussky, {\em Structure des surfaces de Kato,}
M\'em. Soc. Math. France (N.S.) No. 14 (1984).


\bibitem[Du]{_Dulac_} H.  Dulac,
{\em  Solutions d'un syst\`eme d'\'equations diff\'erentielles
 dans le voisinage des valeurs singuli\`eres,}
Bull. Soc. Math. France, 
{\bf 40} (1912), 324-383. 

\bibitem[EF]{_Elencwajg_Forster_}
G. Elencwajg and O. Forster,
{\em Vector bundles on manifolds without divisors and
a theorem on deformations,} Ann. Inst. Fourier 32.4 (1982), 25-51.


\bibitem[FR]{_Favre_Ruggiero_} 
C. Favre, M. Ruggiero, {\em Normal surface singularities
  admitting contracting automorphisms},
Ann. Fac. Sci. Toulouse Math. {\bf 23} (2014), 797--828.




\bibitem[Fo]{_Forster:Steinischen_}
O. Forster,
{\em Zur Theorie der Steinschen Algebren und Moduln}, 
Math. Z., 97:376-405, 1967.

\bibitem[FM]{_Friedland_Milnor:Henon_}
S. Friedland, J. Milnor,  
{\em Dynamical properties of plane polynomial automorphisms,} 
Ergodic Theory Dynam. Systems {\bf 9} (1) (1989), 67-99.

\bibitem[Fr]{friedman} A. Friedman, {\em Foundations of
  modern analysis}, Dover, 2010.






%
%
%
\bibitem[GO]{_GO_Hopf_surfaces_} P.  Gauduchon, L. Ornea, {\em Locally conformally K\"ahler metrics on Hopf surfaces}, Ann. Inst. Fourier (Grenoble) {\bf 48} (1998), 1107-1128.
 
 \bibitem[Gr]{_Grothendieck_Montel_} 
A. Grothendieck, {\em Th\'eor\`emes de finitude pour la cohomologie
 des faisceaux}, 
 Bull. S. M. F., tome {\bf 84} (1956), 1-7.

\bibitem[EGA4]{_EGA:IV_2_}
A. Grothendieck, 
{\em \'El\'ements de g\'eom\'etrie alg\'ebrique: IV. \'Etude locale des sch\'emas et des morphismes de sch\'emas, Seconde partie,} 
Publ. Math. Inst. Hautes Études Sci. {\bf 24} (1965), 5-231. 

\bibitem[SGA3]{_SGA:3_1_}
A. Grothendieck, {\em SGA 3, Exp. I,} in Demazure-Grothendieck (eds.),
Sch\'emas en groupes I, LNM {\bf 151}, Springer, 1970.


\bibitem[GR]{_Gunning_Rossi_} 
R. C. Gunning, H. Rossi,
 {\em Analytic functions of several complex variables}, Reprint
  of the 1965 original. AMS Chelsea Publishing,
  Providence, RI, 2009.



\bibitem[Ha]{_Hartshorne:stable_reflexive_}
R. Hartshorne,  {\em Stable reflexive sheaves}, 
Math. Ann. 254, 121-176 (1980).


 \bibitem[Hu]{_Humphreys_} 
J. E. Humphreys, {\em Linear algebraic groups}, GTM 21, 4th ed., Springer, 1998.

\bibitem[KV]{_KV1_} D. Kaledin, M. Verbitsky. {\em Hyperholomorphic sheaves and new examples of hyperk\"ahler manifolds}. alg-geom 9712012, a chapter in a book ``Hyperk\"ahler manifolds'' International Press, Boston, 2001.



\bibitem[Ka]{_Kato:subvarieties_}
Ma. Kato,
{\em Some Remarks on Subvarieties of Hopf Manifolds,}
    Tokyo J.  Math.
    {\bf 2}, Nr. 1 (1979), 47--61.

  
\bibitem[KK]{_Kaup_} L. Kaup, B. Kaup, {\em Holomorphic Functions of Several Complex Variables}, de Gruyter, 1983.

\bibitem[Kod]{_Kodaira_Structure_III_} 
K. Kodaira, 
{\em On the structure of compact complex surfaces, III},
  Amer. J. Math. {\bf 90} (1968), 55-83.



\bibitem[Kol]{_Kollar:quotient_}
J. Koll\'ar,
{\em Quotient Spaces Modulo Algebraic Groups},
Annals of Math. 
{\bf 145}, No. 1 (1997), 33-79.



\bibitem[Kor1]{_Korshunov:MO_}
D.~Korshunov,
``Polynomial contractions acting as automorphisms of $\mathbb{C}^n$'',
MathOverflow, 2026.
\url{https://mathoverflow.net/questions/507828/polynomial-contractions-acting-as-automorphisms-of-bbb-cn}.

\bibitem[Kor2]{_Korshunov:contractions_}
D.~Korshunov,
{\em Polynomial contractions, degree growth, and exotic algebraic $\C^n$},
arXiv:2605.29386.



%






\bibitem[LT]{_Lubke_Teleman:Book_}
M. L\"ubke, A. Teleman,  {\em 
The Kobayashi-Hitchin correspondence}, World
   Scientific Publishing Co., Inc., River Edge, 
NJ, 1995. x+254 pp

%


\bibitem[Le]{_Lerman:Orbifolds_}
E. Lerman, {\em Orbifolds as stacks?}
L'Enseign. Math. (2) {\bf 56} (2010), no. 3-4, 315--363

\bibitem[Ma]{_Madera:Hopf_}
M. Madera, {\em Holomorphic geometric structures on
  Hopf manifolds}, arXiv:2501.11364


\bibitem[Mi]{_Milne:AG_}
J. S. Milne, {\em Algebraic Groups. The theory of group schemes
	of finite type over a field}, Cambridge University Press; 2017. Also:  
\url{https://www.jmilne.org/math/Books/iAG2022.pdf}

\bibitem[MP]{_Moerdjik_Pronk_}
I. Moerdijk, D. A.  Pronk,  
{\em Orbifolds, sheaves and groupoids}, 
K-Theory, {\bf 12} (1) (1997), 3-21.


\bibitem[Mo]{_Morvan_} 
K. Morvan, {\em Singularities Admitting Contracting 
Automorphisms}, arXiv:2412.11583.

\bibitem[OSS]{_OSS_} 
 C. Okonek, M. Schneider, H.  Spindler,
{\it Vector bundles on complex projective spaces.}
 Progress in math., vol. {\bf 3},
 Birkhauser, 1980.

\bibitem[OW]{_Orlik_Wagreich_}
P. Orlik, P. Wagreich,
{\em Isolated Singularities of Algebraic Surfaces with $\C^*$-Action}
Annals of Math. {\bf 93}, No. 2 (1971), 205-228.




\bibitem[OV1]{_OV_lckpot_}
L. Ornea, M. Verbitsky, 
{\em Locally conformal K\"ahler manifolds with potential},
  Math. Ann. {\bf 348} (2010), 25-33.


\bibitem[OV2]{_OV_pams_} 
L. Ornea, M. Verbitsky, {\em
  Locally conformally K\"ahler metrics obtained from
  pseudoconvex shells}, Proc. Amer. Math. Soc. {\bf 144}
  (2016), 325-335.

%
%

\bibitem[OV3]{ov_indam} 
L. Ornea, M. Verbitsky, {\em Embedding of LCK manifolds with potential into Hopf manifolds using Riesz--Schauder theorem}, ``Complex and Symplectic Geometry'', Springer INdAM series, 2017, 137-148.

\bibitem[OV4]{_OV:surf_in_lck_pot_} 
L. Ornea, M. Verbitsky,
  {\em Hopf surfaces in locally conformally K\"ahler
    manifolds with potential}, Geom. Dedicata {\bf 207}
  (2020), 219-226.

 

\bibitem[OV5]{_OV_non_linear} L. Ornea, M. Verbitsky, {\em Non linear Hopf manifolds are locally conformally K\"ahler},  J. Geom. Anal. {\bf 33}, (2023), art. no 201.    arXiv:2202.12398.

\bibitem[OV6]{_OV_book_}  
L. Ornea, M. Verbitsky, {\em Principles of locally conformally K\"ahler geometry}, Progress in Math. {\bf 354}, Birkh\"auser, 2024, arXiv:2208.07188.



\bibitem[OV7]{_OV_Algebraic_Cones_} 
L. Ornea, M. Verbitsky, 
{\em Algebraic cones on LCK manifolds with
    potential}, J. Geom. Phys. {\bf 198} (2024),
  105103. arXiv:2208.05833.
  
 

\bibitem[OV8]{_OV:Mall_}
L. Ornea, M. Verbitsky, 
{\em Mall bundles and flat connections},
Ann. Inst. Fourier (Grenoble) {\bf 75} (2025), no. 1, 331-358.



\bibitem[P]{_Persson:degene_}
 U. Persson,  
{\em On degenerations of algebraic surfaces}, 
Mem. Amer. Math. Soc. {\bf 11} (1977), no. 189.



\bibitem[Po]{_Poincare:Thesis_} 
H. Poincar\'e, 
{\em Sur les propri\'et\'es des fonctions d\'efinies par les
\'equations aux diff\'erences partielles,}
 Paris, Gauthier-Villars, 1879.










%


%

\bibitem[Ro]{_Roberts:G_sheaves_}
M. Roberts, {\em A note on coherent G-sheaves,} 
Math. Ann. 275 (1986), 573-582.


\bibitem[RT]{_Ross_Thomas:Orbifolds_}
J. Ross, R. Thomas,  {\em 
Weighted projective embeddings, stability of orbifolds, and constant scalar curvature K\"ahler metrics},
J. Differential Geom. 88 (2011), no. 1, 109-159. 

\bibitem[Sa]{_Saito:quasihomo_}
K. Saito, 
{\em Quasihomogene isolierte Singularit\"aten von Hyperfl\"achen,} 
Invent. Math. {\bf 14} (1971), 123-142.


%
%


\bibitem[Se1]{_Serre_Faisceaux_} 
J.-P. Serre, {\em Faisceaux alg\'ebriques 
coh\'erents}, Annals of Math., {\bf 61} (2) (1955) 197-278.

\bibitem[Se2]{_Serre:GAGA_}
J.-P. Serre, {\em G\'eom\'etrie alg\'ebrique et
g\'eom\'etrie analytique}, Ann.Inst. Fourier
{\bf 6}  (1956), 1-42. 

\bibitem[Se3]{_Serre:Prolongement_} 
J.-P. Serre,
{\em Prolongement de faisceaux analytiques coh\'erents},
Ann. Inst. Fourier (Grenoble) 16 (1966), fasc. 1, 363-374.

\bibitem[Si]{_Siu:Extending_}
Y.-T. Siu,
{\em Extending Coherent Analytic Sheaves},
Annals of Math. {\bf 90}, No. 1 (1969), 108-143.


\bibitem[StPr]{_Stacks_Project_}
Aise Johan  de Jong et al.,
{\em The Stacks Project},
{\url{https://stacks.math.columbia.edu/}}

\bibitem[St]{_Stein:58_} K. Stein, {\em \"Uberlagerungen holomorphvollst\"andiger komplexer R\"aume}, Arch. Math. {\bf 7} (1956-7), 354-361.

\bibitem[Ste]{_Sternberg_contraction_} 
S.  Sternberg,
  {\em Local contractions and a theorem of Poincar\'e},
  Amer. J.  Math. {\bf 79} (1957), 809-824.



\bibitem[Sw]{_Swan_} R. G. Swan,   {\em Vector Bundles and Projective Modules}, Transactions  Amer. Math. Soc. {\bf 105} (2) (1962), 264-277.


%

\bibitem[Va]{_Vaisman_trans_} I. Vaisman, {\em On locally and globally conformal K\"ahler manifolds}, Trans. Amer. Math. Soc., {\bf 262} (1980), 533-542.



\bibitem[Ve1]{_Verbitsky:HC_} 
M. Verbitsky, {\em Hypercomplex varieties}, 
Comm. Anal. Geom. {\bf 7} (1999), no. 2, 355-396.



\bibitem[Ve2]{_Verbitsky:Stable_} 
 M. Verbitsky, {\em
Stable bundles on positive principal elliptic fibrations},
 Math. Res. Lett. {\bf 12} (2005), no. 2-3, 251--264.
 
 \bibitem[Ve3]{_Verbitsky:Filtrable_} 
 M. Verbitsky, 
 {\em Holomorphic bundles on diagonal Hopf manifolds}, Izv. Math. {\bf 70} (2006), no. 5, 13--30.
 
 \bibitem[VVO]{_ovv:surf_} 
M. Verbitsky, V. Vuletescu, L. Ornea, 
 {\em Classification of non-K\"ahler surfaces
 	and	locally conformally K\"ahler geometry}, 
 Russian Math. Surv. {\bf 76} (2021), 261-290. 

\bibitem[Vo]{_Vosegaard_}
H. Vosegaard, 
{\em A characterization of quasi-homogeneous complete intersection 
singularities,}
J. Algebraic Geom. {\bf 11} (2002), no. 3, 581-597.

\bibitem[Wa]{_Wang:Harder-Narasimhan_}
Z. Wang
{\em On the volume of a pseudo-effective class and semi-positive properties of the Harder-Narasimhan filtration on a compact Hermitian manifold},
Ann. Polon. Math. {\bf 117} (2016), 41-58

\bibitem[Wu]{_Wu:Montel_} H. Wu, 
{\em Normal families of holomorphic mappings},
Acta Math. {\bf 119} (1967), 193-233.


\bibitem[Zu]{_Zung:groupoids_}
N. T. Zung, {\em Proper groupoids and momentum maps: linearization, affinity, and convexity,} Ann. Sci. Ecole Norm. Sup. (4) {\bf 39} (2006), no. 5, 841-869.



\end{thebibliography}
}
\end{document}